\documentclass[preprint,authoryear,11pt]{elsarticle}
\usepackage{amssymb,amsmath,latexsym}
\usepackage[latin1]{inputenc}
\usepackage[T1]{fontenc}\usepackage{amsmath,amsfonts,amssymb,dsfont}
\usepackage[english]{babel}
\usepackage{latexsym}
\usepackage{epsfig}
\usepackage{graphics}

\usepackage{enumerate}
\usepackage{multirow}

\usepackage{tabularx}
\usepackage[left=2.5cm,right=2cm,top=2cm,bottom=2cm]{geometry}
\usepackage{graphicx}
 \usepackage{natbib}
\usepackage{here}
\usepackage{bm}
\usepackage[linesnumbered]{algorithm2e}
\usepackage{booktabs} 
\usepackage{lipsum}
\usepackage{multicol}
\usepackage{changepage}
\usepackage{dsfont}

\usepackage{accents}
\usepackage{caption}
\usepackage{xcolor} 
\usepackage[colorlinks=true, linkcolor=red, citecolor=red, urlcolor=blue]{hyperref}
\usepackage{cleveref}

\numberwithin{propo}{section}

\numberwithin{thm}{section}

\numberwithin{rem}{section}

\numberwithin{exemple}{section}

\numberwithin{lem}{section}

\numberwithin{cor}{section}
\newtheorem{l1}{Lemma}
\newtheorem{not1}{Notation}
\newtheorem{t1}{Theorem}
\newtheorem{r1}{Remark}
\newcommand{\Var}{\operatorname{Var}} 
\newcommand{\Cov}{\operatorname{Cov}}
\newcommand{\length}[1]{\operatorname{length}\left(#1\right)}

\makeatletter

\@addtoreset{equation}{section}
\makeatother

\begin{document}


\begin{frontmatter}

\title{Estimating weak Markov-switching AR$(1)$ models}

\author[uphf]{Yacouba Boubacar Ma\"{i}nassara \corref{cor1}} \ead{Yacouba.BoubacarMainassara@uphf.fr}
\cortext[cor1]{Corresponding author}
\address[uphf]{ Univ. Polytechnique Hauts-de-France, INSA Hauts-de-France, CERAMATHS - Laboratoire de
Mat\'{e}riaux C\'{e}ramiques et de Math\'{e}matiques, F-59313 Valenciennes, France}
\author[ufc]{Armel Bra}\ead{kja.bra@univ-fcomte.fr}
\address[ufc]{Universit\'{e} Marie et Louis Pasteur, Laboratoire de math\'{e}matiques de Besan\c{c}on,  UMR CNRS 6623, 16 route de Gray,  25030 Besan\c{c}on, France}

\author[ufc]{Landy Rabehasaina}\ead{landy.rabehasaina@univ-fcomte.fr}

\begin{abstract}
In this paper, we present the asymptotic properties of the moment estimator for autoregressive (AR for short) models subject to Markovian changes in regime under the assumption that the errors are uncorrelated but not necessarily independent. We relax the standard independence assumption on the innovation process to extend considerably the range of application of the Markov-switching AR models. We provide necessary conditions to prove the consistency and asymptotic normality of the moment estimator in a specific case. Particular attention is paid to the estimation of the asymptotic covariance matrix. Finally, some simulation studies and an application to the hourly meteorological data are presented to corroborate theoretical work.
\end{abstract}
\begin{keyword}
Weak AR models, Regime-switching models, Markov-switching models, Times series with changes in regime; Moment's method, Asymptotic normality, Asymptotic variance matrix.
\end{keyword}

\end{frontmatter}


\section{Introduction}\label{Introduction}
  Nonlinear models are becoming more and more employed because numerous real time series exhibit nonlinear dynamics. For example, consider a time series  that experiences regime changes at unknown times with a finite number of possible regimes. These models are commonly applied in financial time series, where regimes correspond to significant events that cause high volatility, followed by calmer periods. For instance, as illustrated in \citet[Fig 1.2, p. 7]{francq2010inconsistency}, high-volatility periods are often associated with notable events such as September 11, 2001, or the 2008 financial crisis.
 
In this paper, we investigate an autoregressive model with random coefficients, where the associated noise exhibits a multiplicative structure that depends on both a Markov chain and an exogenous noise. These models can be viewed as Markovian mixtures of dynamic systems, belonging to the class of Markov regime-switching models. More precisely, a Markov-switching model is a non-linear specification in which different states of the world affect the evolution of a time series (see, for examples, \cite{francq1997white, Hamilton1990, HS1994autoregressive}). Such models have attracted significant interest in the literature, with foundational contributions by \cite{hamilton1988rational}, \cite{hamilton1989new}, \cite{mcculloch1994bayesian} and \cite{chib1996calculating}. Their statistical properties have been extensively studied, for instance, by \cite{billio1999bayesian}. Recent research has further enriched this field. For instance, \cite{francq1997white} examined a time series model where the variance of the underlying process depends on the state of an unobserved Markov chain, considering only multiplicative noise. They proposed a maximum likelihood estimator and studied its asymptotic properties. This work was later extended by \cite{francq1998ergodicity} to encompass AR processes with random coefficients. The authors established conditions for the existence of a stationary and ergodic solution and proved the consistency of the maximum likelihood estimator. Another contribution is in \cite{xie2008general} who studied a general AR model with Markov regime-switching, allowing for AR with infinite order. Under some regular assumptions they demonstrated the consistency of the maximum likelihood estimators. We can also cite \cite{douc2004asymptotic} who studied the asymptotic properties of the maximum likelihood estimator for an AR process with Markov regime switching, potentially nonstationary, where the hidden state space is compact but not necessarily finite. They demonstrated consistency and asymptotic normality under the assumption of uniform exponential forgetting of the initial distribution of the hidden Markov chain given the observations. Additionally, \cite{FG2004jtsa}  investigated the estimation of time-varying Autoregressive Moving Average (ARMA for short) models with Markovian regime changes, where they  gave general conditions ensuring the consistency and asymptotic normality of least squares and quasi-generalized least squares estimators. \cite{francq2004estimation} provided explicit conditions ensuring the consistency and asymptotic normality of least squares and quasi-generalized least squares estimators and gave the asymptotic covariance matrix of the estimators when the changes between states are governed by the outcome of a Markov chain. Note also that in \cite{FG2004jtsa}  and \cite{francq2004estimation} the realization of the Markov chain is assumed to be observed. This body of work underscores the growing interest and ongoing advancements in Markov regime-switching models, as researchers continue to refine estimation techniques and broaden the applicability of these models. 
 All the works cited above have been conducted under the assumption that the noise is independent and identically distributed (i.i.d. for short).

 As above-mentioned, the works on the statistical inference of AR processes with Markov regime switching are generally performed under the assumption that the errors are independent. This independence assumption is often considered too restrictive by practitioners. It precludes conditional heteroscedasticity and/or other forms of nonlinearity (see \cite{fz05} for a review on ARMA models under the assumption that the errors are uncorrelated but not necessarily independent) which can not be generated by Markov regime switching models  with i.i.d. noises. Relaxing this independence assumption  allows to extend the range of application of the class of Markov regime switching models.

In this paper we focus on an AR$(1)$ process modulated by a hidden Markov chain with multiplicative noise, under the assumption that the errors are uncorrelated but not necessarily independent. For brevity, we refer to this as a weak AutoRegressive Hidden Markov Chain  (ARHMC for short) model. Conversely, when the noise is assumed to be i.i.d., we call it the strong ARHMC model. The term \textit{hidden} reflects the fact that the states of the Markov chain are not directly observable, yet they play a significant role in shaping the behavior of the time series. By relaxing the independence assumption, we extend the applicability of ARHMC models to encompass more complex nonlinear processes, including those with multiplicative noise structures akin to generalized autoregressive conditional heteroscedastic (GARCH for short) models introduced by \cite{engle1982autoregressive} and extended by \cite{bollerslev1986generalized} (see also \cite{francq2019garch}, for a reference book on GARCH models). These distinctions are critical for understanding the specific challenges and nuances of the methodology developed in this study. It is worth noting that very few studies have considered time series models with Markovian regime changes involving such noise. A notable exception \cite{francq2001stationarity} who investigated the stationarity conditions of such models in a multivariate framework. The authors demonstrated that local stationarity of these processes is neither sufficient nor necessary to ensure global stationarity. And finally \cite{boubacar2020estimation} who examined the asymptotic properties of the least squares estimator for a weak ARMA model with regime switching under the assumptions that the  realization of the Markov chain is observed. However, it is important to note also that in the model considered by \cite{boubacar2020estimation} the structure of the noise is not multiplicative. Instead, they assumed that the volatility was uniform within each segment of the time series.

 To our knowledge, it does not exist any estimation methodology for weak ARHMC models when the (possibly dependent) error is subject to known or unknown conditional heteroscedasticity. This paper is devoted to the problem of the estimation  of weak ARHMC processes. We propose the moment estimation procedure to estimate the parameters of a weak ARHMC model. We show that a strongly mixing property and the existence of moments are sufficient to obtain a consistent and asymptotically normally distributed of the proposed estimator.

 In our opinion there are two major contributions in this work. The first one is to show that the moment estimation procedure can be extended to weak ARHMC models. This goal is achieved thanks to Theorems \ref{thm:2} and \ref{thm:3} in which the consistency and the asymptotic normality are stated. The second one is to provide a weakly consistent estimator of the asymptotic variance matrix (see Theorem \ref{thm:4}). Thanks to this estimation of the asymptotic variance matrix, we can
construct a confidence region for the estimation of the parameters. Finally we extend the existing results on the statistical analysis of ARHMC models by addressing the estimation problem under more general error structures.

 The structure of the paper is as follows. Section \ref{sect2} introduces the weak ARHMC model that we consider here and outlines the underlying assumptions. Our methodology based on moments method is given in Section \ref{sect3} and the main results are given in Section \ref{sect4}. We provide a consistency analysis, showing that the moments estimator converges almost surely to the true parameter, along with the asymptotic normality of the moment estimator under certain mixing conditions for the linear innovation process. Notably, the asymptotic covariance of the moments estimator differs significantly between the weak and strong cases. Section \ref{sect5} is devoted to the estimation of this covariance matrix. The simulation studies and  illustrative applications on real data are presented and discussed in In Section \ref{sect6}.
The proofs of the main results are collected in Section \ref{proofs}.




\section{Model and assumptions}\label{sect2}

Let $(\Delta_t)_{t\in \mathbb{Z}}$ be an unobserved homogeneous Markov chain with a stationary distribution $\pi$ taking values in a discrete set $\mathcal{S}:=\{1,\dots,K\}$ and transition matrix $P=(p_{ij})_{i,j=1,\dots,K}$. We consider a stationary weak ARHMC$(1)$ process $(X_t)_{t \in \mathbb{Z}}$ defined as: 
\begin{align}
\label{eq:1}
   X _t = a(\Delta_t)X_{t-1}+f(\Delta_t)\eta_t, \  \ \forall t\in \mathbb{Z},
\end{align}
where the process $(\eta_t)_{t\in \mathbb{Z}}$ is a weak white noise satisfying $\mathbb{E}(\eta_t) = 0$, $\mathbb{E}(\eta_t \eta_{t'}) = \sigma^2\mathds{1}_{[t = t']}$ with $\sigma^2>0$, $f : \mathcal{S}\rightarrow \mathbb{R}\setminus\{0\}$, and $a : \mathcal{S}\rightarrow \mathbb{R}$. Without loss of generality, we will assume that \( \sigma^2 = 1 \). An example of weak white noise is the GARCH model (see \cite{francq2019garch}). It is customary to say that $(X_t)_{t\in \mathbb{Z}}$ is a strong ARHMC$(1)$ representation and we will do this henceforth if in (\ref{eq:1}) $(\eta_t)_{t\in \mathbb{Z}}$ is a strong white noise, namely an i.i.d. sequence of random variables with mean 0 and common variance 1. A strong white noise is obviously a weak white noise because independence entails uncorrelatedness. Of course, the converse is not true. It is clear from these definitions that the following inclusion hold: 
\begin{align*}
    \{\text{Strong ARHMC(1)}\} \subset \{\text{Weak ARHMC(1)}\}.
\end{align*}
In the rest of the paper, we will denote by \( M' \)  the transpose of the matrix \(M\).
 The unknown parameter of interest  is denoted $\theta_0:=(a(s),(p_{ij}),(f(s)), 1\leq i\leq K, 1\leq j\leq K-1,s\in \mathcal{S})'$ and belongs to the parameter space \begin{equation*}
\begin{split}
\Theta := \biggl\{
&\theta = (\theta_{11}, \dots, \theta_{KK}, \overline{\theta}_{11}, \overline{\theta}_{12}, \dots, \overline{\theta}_{1K-1}, \overline{\theta}_{21}, \overline{\theta}_{22}, \dots, \overline{\theta}_{2K-1}, \dots, \overline{\theta}_{K1}, \cdots, \overline{\theta}_{KK-1}, \\
&\tilde{\theta}_{11}, \dots, \tilde{\theta}_{KK})'\in \mathbb{R}^{K} \times {(0,1)}^{K(K-1)} \times \mathbb{R}^{K}\setminus\{0_{\mathbb{R}^K}\};
P_{\theta} : = (\overline{\theta}_{ij})_{1\leq i,j\leq K} \ \text{is an irreducible }  \\ &\text{transition matrix, }\; \overline{\theta}_{j,K} = 1-\displaystyle\sum_{i=1}^{K-1} \overline{\theta}_{j,i}\in (0,1), j=1,\dots,K;  A_{\theta}: = \text{diag}(\theta_{11}, \dots, \theta_{KK})  \\ &\hspace{4cm} \text{and the spectral radius of}\ A_{\theta}^2P_{\theta}' \ \text{is less than 1}
\biggr\}.
\end{split}
\end{equation*}
 In order to measure the temporal dependence of the processes \((\eta_t)_{t\in\mathbb{Z}}\) and \((\Delta_t)_{t\in\mathbb{Z}}\), we define the strong mixing coefficients \( (\alpha_\mathcal{Z}(h))_{h\in\mathbb{N}^\star} \), which are independent of \(t \in \mathbb{Z}\), for a stationary process \((\mathcal{Z}_t)_{t\in\mathbb{Z}}\) as follows:  
\begin{align}\label{eq:12}
\alpha_{\mathcal{Z}}(h) := \sup_{A\in \mathcal{F}_{-\infty}^t, B\in \mathcal{F}_{t+h}^\infty} \big|\mathbb{P}(A \cap B) - \mathbb{P}(A)\mathbb{P}(B)\big|,
\end{align}  
where \(\mathcal{F}_{-\infty}^t\) and \(\mathcal{F}_{t+h}^\infty\) denote the \(\sigma\)-fields generated by \(\{{\cal Z}_u : u \leq t\}\) and \(\{{\cal Z}_u : u \geq t + h\}\), respectively.

Our main results are proven under the following assumptions:
\begin{enumerate}
\item[$(\mathbf{A_1})$] The processes \((\eta_t)_{t\in\mathbb{Z}}\) and \((\Delta_t)_{t\in\mathbb{Z}}\) are  stationary and  $(\eta_t)_{t\in\mathbb{Z}}$ is ergodic.
\end{enumerate}
Note also that the process \((\Delta_t)_{t\in\mathbb{Z}}\) is  ergodic since the matrix $P=P_{\theta_0}$ is irreducible. 
In the sequel we suppose that there exists some constant $\nu>0$  
such that:
\begin{enumerate}
\item[$(\mathbf{A}_2)$] The spectral radii of \( A_{\theta_0}^{\beta}P'_{\theta_0} \) are each strictly less than 1, for \( \beta \leq \max\{8, 4 + 2\nu\} \) where \( A_{\theta_0} = \text{diag}(a(s), s \in \mathcal{S}) \). 

\item[$(\mathbf{A_3})$] The processes $(\Delta_t)_{t\in\mathbb{Z}}$ and $(\eta_t)_{t\in\mathbb{Z}}$ are independent.
\item[$(\mathbf{A_4})$] $\displaystyle\sum_{h = 0}^{\infty}\alpha_{\eta}(h)^{\frac{\nu}{2+\nu}}<\infty$, \ $\displaystyle\sum_{h = 0}^{\infty}\alpha_{\Delta}(h)^{\frac{\nu}{2+\nu}}<\infty$ and $\mathbb{E}(|\eta_t|^{4+2\nu})<\infty$.
\item[$(\mathbf{A_5})$] We have $\theta_0\in\mathring{\Theta}$, where $\mathring{\Theta}$ denotes the interior of $\Theta$. 
\end{enumerate}
 Under Assumptions $(\mathbf{A_1})$ and $(\mathbf{A_3})$, the process \((\Delta_t, \eta_t)_{t\in\mathbb{Z}}\) is ergodic. Consequently, the sequence \((a(\Delta_t), \epsilon_t)_{t \in \mathbb{Z}}\) is also strictly stationary and ergodic. Furthermore \(\mathbb{E}\log^{+}|a(\Delta_0)|\) and \(\mathbb{E}\log^{+}|\epsilon_0|\) are finite (where $\log^+(x) = \max\{\log(x),0\}, x > 0$). Additionally we assume that: 
\begin{align*}
    (\mathbf{A_6})  \text{ The Lyapunov exponent defined by }\gamma := \inf_{t \in \mathbb{N}^\star} \left\{ \mathbb{E} \frac{1}{t} \log |a(\Delta_t) a(\Delta_{t-1}) \ldots a(\Delta_1)| \right\} \;\text{is negative}.
\end{align*}
Then using  \citet[Theorem 1.1, page 1715]{bougerol1992strict}  (see also \citet[Theorem 1, page 212]{brandt1986stochastic}), the series
\begin{equation}
\label{eq:4}
X_t = \sum_{k=0}^{\infty} \left( \prod_{j=0}^{k-1} a(\Delta_{t-j}) \right) f(\Delta_{t-k}) \eta_{t-k}, \ \forall t \in \mathbb{Z},
\end{equation}
converges almost surely and is the unique strictly stationary solution of (\ref{eq:1}) with the usual convention that \(\prod_{j=\ell}^{\ell'} = 1\) for \(\ell, \ell' \in \mathbb{Z}\) and \(\ell' < \ell\).

\begin{r1}
Under the assumptions of ergodicity for the processes \((\Delta_t)_{t \in \mathbb{Z}}\) and \((\eta_t)_{t \in \mathbb{Z}}\) and the independence between \((\Delta_t)_{t \in \mathbb{Z}}\) and \((\eta_t)_{t \in \mathbb{Z}}\), the key to ensuring the strict stationarity of the model (\ref{eq:1}) lies in the hypothesis that \(\gamma < 0\). It is worth noting that \(\gamma < 0\) for strict stationarity is a more general condition, applicable to AR$(p)$ and ARMA$(p,q)$ models. However in the specific context of our study, a sufficient condition for strict stationarity as stated in \emph{\cite{francq2001stationarity}} is that
\begin{align}
\label{eq:p1234}
    \sum_{i \in \mathcal{S}} \pi(i) \log|a(i)| < 0,
\end{align}
where $\pi(i)$ is the $i$-th component of the stationary distribution $\pi$.
\end{r1}

\section{Estimation of the ARHMC$(1)$ model parameters}\label{sect3}
We state by the following theorem which provides an explicit expression of the autocovariance function of order $k$ of the centered process $(X_t)_{t\in\mathbb{Z}}$.
\begin{t1}\label{thm:1}
Under Assumptions $(\mathbf{A_1})$, $(\mathbf{A_2})$, $(\mathbf{A_3})$ and $(\mathbf{A_6})$, the joint moments of the process $(X_t)_{t\in\mathbb{Z}}$ defined in (\ref{eq:1}) satisfy for all $k\in\mathbb{N}$:
\begin{equation}
\label{eq:2}
    c_{k,0}(\theta_0) : = \mathbb{E}(X_kX_0) = \bm{1}'(A_{\theta_0}P_{\theta_0}')^k(I_K-A_{\theta_0}^2P_{\theta_0}')^{-1}\bm{\pi}_{f^2},
\end{equation}
where $I_K$ denotes the identity matrix of size $K$, $\bm{1}' := (1, \dots, 1)$ is a row vector of dimension $K$ and $\bm{\pi}_{f^2} := \{f^2(1)\pi(1), \dots, f^2(K)\pi(K)\}'$ is a column vector of dimension $K$.
\end{t1}
The proof of this theorem is given in Section \ref{proof_thm1}.

Thanks to Theorem \ref{thm:1} and in order to state our asymptotic normality result, we will explain our estimation procedure. 
 In the following we denote by \(\pi_\theta\) the stationary distribution associated to the Markov chain with transition matrix $P_\theta$, parametrized by some \(\theta \in \Theta\). 
 Since \(\pi_\theta\) is the unique solution to \(\pi_\theta P_\theta = \pi_\theta\) and \(\pi_\theta \bm{1} = 1\), there exists an invertible matrix \(B_\theta \in \mathbb{R}^{(K^2 + K) \times (K^2 + K)}\) and a vector \(v \in \mathbb{R}^{K^2 + K}\) such that \(\pi_\theta = v' {B_\theta'}^{-1}\). For instance we can take 
\begin{align}
B_\theta : = 
\begin{pmatrix}
    \overline{\theta}_{11}-1 & \overline{\theta}_{21} & \cdots & \overline{\theta}_{K1} \\
   \overline{\theta}_{12} & \overline{\theta}_{22}-1 & \cdots & \overline{\theta}_{K2}\\
     \vdots &  \ddots & \vdots & \vdots  \\
     \overline{\theta}_{1K-1} & \cdots & \overline{\theta}_{K-1K-1}-1 & \overline{\theta}_{KK-1}\\
  1 & \cdots & 1 & 1
\end{pmatrix} \ \text{and} \ 
v : = 
\begin{pmatrix}
    0,& \dots &, 0, & 1
\end{pmatrix}'.
\end{align}
We recall that $\bm{1}'$ represents a row vector of dimension $K$ with all elements equal to 1 and $I_K$ denotes the identity matrix of size $K$.
Consider an integer $N\geq K^2+K$ fixed in the following and let $\theta \in \Theta$. We then define the diagonal matrix $V_{\theta}$ by
\begin{align*}
   V_{\theta} := \operatorname{diag}(\tilde{\theta}_{11}^2, \ldots, \tilde{\theta}_{KK}^2)
\end{align*}
and introduce the functions $\psi_k(\theta)$ for each index $1 \leq k \leq N$, defined by
\begin{align*}
   \psi_k(\theta) := c_{k,0}(\theta)= \bm{1}' (A_{\theta} P'_{\theta})^k (I_K - A_{\theta}^2 P'_{\theta})^{-1} V_{\theta}\pi_{\theta}'.
\end{align*}

We then define the function $\Psi^N$ which maps each element $\theta$ from $\mathring{\Theta}$  to a vector in $\mathbb{R}^{N}$ as follows
\begin{equation}
\Psi^N :
    \left\lbrace
        \begin{aligned}
            & \mathring{\Theta} & \rightarrow & \ \mathbb{R}^{N} \\
            & \theta & \mapsto & (\psi_1(\theta), \dots, \psi_{N}(\theta))'.
        \end{aligned}
    \right.
\end{equation}
  The function $\Psi^{N}$ defined in this way is differentiable for all $\theta \in \mathring{\Theta}$ because each $\psi_k$ involves products and compositions of differentiable functions. Using matrix differentiation formulas (see \citet[section 2]{petersen2008matrix}), it is possible to obtain the Jacobian matrix of \(\Psi^{N}\), denoted by \(J_{\Psi^N}(\theta)\in\mathbb{R}^{N\times (K^2+K)}\), in explicit form at any point
 \(\theta\in \mathring{\Theta}\).  For $1\leq i\leq K$ and $1\leq j\leq K-1$, the entries of \(J_{\Psi^N}(\theta)\) are given by
  
\begin{equation}
\label{eq:9966}
\left\lbrace
\begin{array}{ll}
    \frac{\partial \psi_{k}(\theta)}{\partial \theta_{ii}} = & \mathbf{1}'\left(\displaystyle\sum_{r = 0}^{k-1}(A_\theta{P_\theta}')^r \{\partial A_\theta/\partial\theta_{ii}\}{P_\theta}'(A_\theta{P_\theta}')^{k-1-r}(I_K-{A^2_\theta}{P_\theta}')^{-1} \right. \\ 
    & \left. + \ (A_\theta{P_\theta}')^k(I_K-{A^2_\theta}{P_\theta}')^{-1}\left(\displaystyle\sum_{r=0}^1{A_\theta}^r \{\partial A_\theta/\partial\theta_{ii}\}{A_\theta}^{1-r}\right){P_\theta}'(I_K-{A^2_\theta}{P_\theta}')^{-1}V_\theta\pi'_\theta \right) \\[10pt]

    \frac{\partial \psi_k(\theta)}{\partial \overline{\theta}_{ij}} = & \mathbf{1}'\Bigg(\Bigg(\displaystyle\sum_{r = 0}^{k-1}(A_\theta{P_\theta}')^r A_\theta \{\partial P_\theta/\partial\overline{\theta}_{ij}\}'(A_\theta{P_\theta}')^{k-1-r}(I_K-{A^2_\theta}{P_\theta}')^{-1} \\ 
    &  + \ (A_\theta{P_\theta}')^k(I_K-{A^2_\theta}{P_\theta}')^{-1}({A^2_\theta}\{\partial P_\theta/\partial\overline{\theta}_{ij}\}')(I_K-{A^2_\theta}{P_\theta}')^{-1}\Bigg)(V_\theta\pi_\theta') \\ 
    &  + \ (A_\theta{P_\theta}')^k(I_K-{A^2_\theta}{P_\theta}')^{-1}V_\theta(-B_\theta^{-1}\{\partial B_\theta/\partial\overline{\theta}_{ij}\} B_\theta^{-1})v \Bigg) \\[10pt]
  
    \frac{\partial \psi_k(\theta)}{\partial \tilde{\theta}_{ii}} = & \mathbf{1}'\left(A_{\theta}P_{\theta}'\right)^k\left(I_K-A^2_\theta P_{\theta}'\right)^{-1}\{\partial V_\theta/\partial\tilde{\theta}_{ii}\}\pi_{\theta}',
\end{array}
\right.
\end{equation}
where $ {\partial A_{\theta}}/{\partial \theta_{ii}}$, ${\partial B_{\theta}}/{\partial \overline{\theta}_{ij}}$, ${\partial P_{\theta}}/{\partial \overline{\theta}_{ij}}$ and 
${\partial V_\theta}/{\partial \tilde{\theta}_{ii}}$ are then explicit, which allows us to represent the coefficients of the matrix \( J_{\Psi^N}(\theta) \) in a closed form. An alternative and straightforward method to compute this matrix is by using a symbolic computation software such as the \texttt{SymPy} library in Python. We will adopt this second approach for our upcoming simulations.

In order to estimate the parameter \(\theta_0\) we thus have at our disposal the observations $(X_1,\dots,X_n)$. We employ the Newton-Raphson method which is widely-used for finding roots of real-valued functions and even for vector-valued functions. To implement this, we introduce our estimation function \(\mathcal{F}^{N,n}\) defined as
\begin{align}
\label{eq:800}
    \mathcal{F}^{N,n}(\theta):=  J_{\Psi^{N}}(\theta)' F^{N,n}(\theta)\in \mathbb{R}^{K^2+K}\quad \text{for all}\quad\theta \in \Theta ,
\end{align}   
with
\[
   F^{N,n}(\theta): = (\hat{c}_{1,0} - c_{1,0}(\theta), \ldots, \hat{c}_{N,0} - c_{N,0}(\theta))'\in\mathbb{R}^N,
\]
where the random scalar \(\hat{c}_{k,0}\) is an estimator of the theoretical moment \(c_{k,0}(\theta_0)\) and is defined as
\begin{align}
\label{eq:3.1}
\hat{c}_{k,0} := (n-k)^{-1}\sum_{t = 1}^{n-k} X_{t+k}X_t, \ \forall \ 1 \leq k \leq N<n.
\end{align}
Note that  \(\hat{c}_{k,0}\) converges a.s. to $c_{k,0}(\theta_0)$ as $n\to\infty$ for all \(1 \leq k \leq N\) by the ergodic theorem and the fact that the process \((X_t)_{t\in\mathbb{Z}}\) is stationary. 
Let \(\hat{\theta}_n\) be the estimator of \(\theta_0\) obtained by the Newton method through the estimation function \(\mathcal{F}^{N,n}\).
 Formally, for large $n$, we define the random variable \( \hat{\theta}_n \) as the solution to:  
\begin{align}
\label{eq:5}
    \mathcal{F}^{N,n}(\hat{\theta}_n) = 0  \quad \text{and} \quad  \det(\nabla \mathcal{F}^{N,n}(\hat{\theta}_n)) \neq 0 \ \text{a.s} .
\end{align}
The existence of this solution and the consistency are proved in the following Theorem \ref{thm:2}.

We denote by \( J_{\mathcal{F}^{N,n}} \) the Jacobian matrix of the random function \( \mathcal{F}^{N,n}\) defined in Equation \eqref{eq:800}. The construction of the estimator \(\hat{\theta}_n\) via the Newton-Raphson method subject to the constraints imposed by the parameters of model  (\ref{eq:1}) is described as follows:
\begin{center}
\begin{algorithm}[H]
\label{algo:algo_1}
  \SetAlgoLined
  \KwData{An initial parameter $\theta^0 \in \mathbb{R}^{K^2 + K}$, tolerance $\varepsilon > 0$}
  \KwResult{A value of $\theta\in\Theta$ such that $\|\mathcal{F}^{N,n}(\theta)\| < \varepsilon$}

  \textbf{Initialization :} $\theta^0 \in \mathbb{R}^{K^2 + K}$\;
  \For{$k \geq 0$}{
    Compute $J_{\mathcal{F}^{N,n}}(\theta^{(k)})$\;
    Update $\theta^{(k+1)}$ using the equation
    \[
    J_{\mathcal{F}^{N,n}}(\theta^{(k)})(\theta^{(k+1)} - \theta^{(k)}) = -  \mathcal{F}^{N,n}(\theta^{(k)})
    \]
    \If{$\|\mathcal{F}^{N,n}(\theta^{(k+1)})\| < \varepsilon$}{
      \textbf{break} \; 
    }
  }
  \caption{Newton's method for finding a root of $\mathcal{F}^{N,n}$ with specified tolerance.}
\end{algorithm}
\end{center}
\vspace{0.3cm}

 At each iteration $k \in \mathbb{N}$, it is necessary to compute $J_{\mathcal{F}^{N,n}}(\theta^{(k)})$ and to solve a linear system. However, even if we start from an initial point $\theta_0$ where $J_{\mathcal{F}^{N,n}}(\theta_0)$ is invertible, there is no guarantee that $J_{\mathcal{F}^{N,n}}$ remains invertible for $\theta^{(1)}$ and subsequent iterations. Consequently, solving the system $J_{\mathcal{F}^{N,n}}(\theta^{k})\left(\theta^{(k+1)}-\theta^{(k)}\right) = - \mathcal{F}^{N,n}(\theta^{(k)})$ can quickly become very costly in terms of time, not to mention the need to project $\theta^{(k)}$ onto the parameter space $\Theta$.
To optimize computation time, one approach is to replace $\mathcal{F}^{N,n}(\theta^{(k)})$ with a linear approximation, through a $(K^2+K)\times(K^2+K)$  matrix $B^{(k)}$ close to $J_{\mathcal{F}^{N,n}}(\theta^{(k)})$ and easily invertible at each iteration. Hence, we aim to construct a matrix $B^{(k)}$ such that when $\theta^{(k)}$ and $\theta^{(k-1)}$ are known, it satisfies the following condition
\begin{align}
\label{eq:999}
    B^{(k)}(\theta^{(k)}-\theta^{(k-1)}) = \mathcal{F}^{N,n}(\theta^{(k)}) - \mathcal{F}^{N,n}(\theta^{(k-1)}).
\end{align}
One way to choose $B^{(k)}$ is to use Broyden's method, as detailed in \citet[page 312]{gomes1992column}. This involves selecting $B^{k}$ that meets condition (\ref{eq:999}) and such that for every vector $\zeta \in \mathbb{R}^{K^2+K}$ orthogonal to $\delta^{(k)} := \theta^{(k)}-\theta^{(k-1)}$, we have $B^{(k)}\zeta = B^{(k-1)}\zeta$.
Consequently, the Broyden algorithm for estimating \(\hat{\theta}_n\) is formulated as follows
\begin{center}
\begin{algorithm}[H]
  \SetAlgoLined
  \KwData{Initial parameters $\theta^{(0)}, \theta^{(1)} \in \mathbb{R}^{K^2 + K}$, initial matrix $B^{0} \in \mathbb{R}^{(K^2+K)\times(K^2+K)}$}
  \KwResult{A value of $\theta\in\Theta$ such that $\|\mathcal{F}^{N,n}(\theta)\| < \varepsilon$}

  \textbf{Initialization :} $\theta^{(0)}, \theta^{(1)} \in \mathbb{R}^{K^2 + K}$, $B^{0} \in \mathbb{R}^{(K^2+K)\times(K^2+K)}$\;
  \For{$k \geq 1$}{
    Compute $\delta^{(k)} = \theta^{(k)} - \theta^{(k-1)}$\;
    Update $B^{(k)}$ using the equation :
    \[
    B^{(k)} = B^{(k-1)} + \frac{\mathcal{F}^{N,n}(\theta^{(k)}) - \mathcal{F}^{N,n}(\theta^{(k-1)}) - B^{(k-1)}\delta^{(k)}}{\|\delta^{(k)}\|^2} \delta^{(k)^{'}}
    \]
    Update $\theta^{(k)}$ using the equation 
    \[
    B^{(k)}(\theta^{(k)} - \theta^{(k-1)}) =  \mathcal{F}^{N,n}(\theta^{(k)}) - \mathcal{F}^{N,n}(\theta^{(k-1)})
    \]
    \If{ $\|\delta^{(k)}\| < \varepsilon$}{
      \textbf{break} \; 
    }
  }
  \caption{Broyden's algorithm for finding a root of $\mathcal{F}^{N,n}$ with specified tolerance.}
\end{algorithm}
\end{center}

\begin{r1}
One of the primary advantages of the Broyden algorithm is that it eliminates the need to recalculate the Jacobian matrix \(J_{\mathcal{F}^{N,n}}\) at every iteration, a process that proves to be extremely costly in the context of our problem. However, this method introduces a significant drawback, the loss of quadratic convergence. Nevertheless, it is important to emphasize that, although quadratic convergence is lost, this does not significantly affect the accuracy of \(\hat{\theta}_n\). Moreover, the constructed matrix \(B^{(k)}\) can be considered as an approximation of \(J_{\mathcal{F}^{N,n}}(\theta^{(k)})\). This observation will be useful in Section \ref{sect6}, where this algorithm will be put into use.
\end{r1}

\begin{r1}
In the following, the choice of the parameter \(N\) will also play a crucial role.
Indeed, $N$ is chosen large enough to be able to estimate the parameter $\theta_0$.
By denoting \(\mathbf{r}_N(\theta)\) as the rank of \(J_{\Psi^N}(\theta)\) for all $\theta\in \Theta$, we observe that the sequence \((\mathbf{r}_N(\theta))_{N \in \mathbb{N}}\) is increasing. Since this sequence takes values in \(\mathbb{N}\), it converges to a limit which we denote by \(\mathbf{r}(\theta)\) and is stationary from a certain rank onward.
\end{r1}

\section{Asymptotic properties}\label{sect4}
\subsection{Consistency and asymptotic normality of the moments estimator}
The asymptotic properties of the estimator \(\hat{\theta}_n\) obtained via the Newton algorithm \ref{algo:algo_1} are stated in the following two theorems.
\begin{t1}
\label{thm:2}
 Let us assume that the limiting rank of \(J_{\Psi^N}(\theta_0)\) as $N\to \infty$ satisfies $\mathbf{r}(\theta_0) = K^2+K$, and let $N \in\mathbb{N}$ such that $\mathbf{r}_N(\theta_0)=\mathbf{r}(\theta_0)$. There exists a neighborhood $\mathcal{V}_{\theta_0}$ of $\theta_0$ in $\Theta$ and a unique sequence $(\hat{\theta}_n)_{n\in\mathbb{N}}$ taking values in $\mathcal{V}_{\theta_0}$ such that
\begin{align*}
    0 = J_{\Psi^{N}}(\hat{\theta}_n)'F^{N,n}(\hat{\theta}_n) = \mathcal{F}^{N,n}(\hat{\theta}_n).
\end{align*}
Furthermore, we have
\begin{align*}
    \hat{\theta}_n \xrightarrow[n\rightarrow\infty]{\text{a.s}} \theta_0.
\end{align*}
\end{t1}
The proof of this theorem is given in Section \ref{proof_thm2}.

The following theorem establishes the asymptotic normality of \(\hat{\theta}_n\). 
\begin{t1}\label{thm:3}
Assuming that conditions \((\mathbf{A_1})\), \((\mathbf{A_2})\), \((\mathbf{A_3})\), \((\mathbf{A_4})\), \((\mathbf{A_5})\) \text{and} $(\mathbf{A_6})$ are satisfied, that the $\mathbf{r}(\theta_0) = K^2+K$. Let $N \in\mathbb{N}$ such that $\mathbf{r}_N(\theta_0)=\mathbf{r}(\theta_0)$ and let $(\hat{\theta}_n)_{n\in\mathbb{N}}$ be a sequence of moments estimator defined in Equation (\ref{eq:5}), of which existence is justified in Theorem \ref{thm:2}. We have 
\begin{align*}
\sqrt{n}(\hat{\theta}_n - \theta_0) \xrightarrow[n\rightarrow\infty]{\mathcal{D}} \mathcal{N}(0,\Omega:=M^{-1}J'IJ'M^{-1}),
\end{align*}
where the matrices $M$, $I$ and $J$ are defined as follows 
\begin{align}
\label{eq:000}
    J & := J_{\Psi^N}(\theta_0) = \left(\nabla\psi_1(\theta_0),\dots,\nabla\psi_{N}(\theta_0)\right)', \notag\\
    I & := I^N(\theta_0) = \sum_{k=-\infty}^{\infty}\mathrm{Cov}\left(Y_t(\theta_0), Y_{t-k}(\theta_0)\right),\notag\\
    M & : = M^N(\theta_0) = J_{\Psi^N}(\theta_0)'J_{\Psi^N}(\theta_0)  
\end{align}
with
\begin{align}
\label{eq:001}
    Y_t:= Y_t(\theta_0)=X_t(X_{t+1},\dots,X_{t+N})'.
\end{align}
\end{t1}
The proof of this theorem is given in Section \ref{proof_thm3}.

\begin{r1}
It is essential to highlight that the hypothesis \(\mathbf{r}(\theta_0) = K^2 + K \), which implies the invertibility of \( J_{\Psi^N}(\theta_0)'J_{\Psi^N}(\theta_0) \), remains crucial for establishing the asymptotic properties of our moment estimator. Although the choice of \( N \) is necessary to ensure that \( \mathbf{r}(\theta_0) = K^2 + K \), the selection of the parameter \( \theta_0 \) is equally significant. Indeed, for certain choices of \( \theta_0 \), we may have \(\mathbf{r}(\theta_0) < K^2 + K \) regardless of the choice of \( N \).
An example is provided in the following Section \ref{exemple}.
\end{r1}
\subsection{An example where \(\mathbf{r}(\theta_0) < K^2 + K \)}\label{exemple}
%
Assume here that the  matrix $A_{\theta_0}$ and the transition matrix $P_{\theta_0}$ have the particular following forms:
\begin{align*}
    A_{\theta_0} : =a_{0}I_K  \quad  \text{and} \quad P_{\theta_0} : = K^{-1}\mathbf{1}\mathbf{1}',
\end{align*}
where $a_0$  is a non zero scalar and such that $|a_0|<1$, so that the stability condition $(\mathbf{A_6})$ is satisfied. Since \(P_{\theta_0}\) is symmetric, we have 
\begin{align*}
    \left(A_{\theta_0}P_{\theta_0}'\right)^\ell  = a_{0}^\ell P_{\theta_0} \quad \forall \ell \in \mathbb{N}.
\end{align*}
Moreover, one easily computes that $J^{ii} := \partial A_\theta / \partial \theta_{ii} \big|_{\theta = \theta_0}$ is the diagonal $K\times K$ matrix of which entries are $0$ save for the $i$-th diagonal entry which is equal to $1$ for $i=1,\dots,K$, from which we can observe that
\begin{align*}
    J^{ii}P_{\theta_0}' = (P_{\theta_0}J^{ii})' = 
    \begin{pmatrix}
        0 & \cdots & P_{\theta_{0_{1,i}}} & \cdots & 0 \\
        \vdots & \vdots & \vdots & \vdots & \vdots \\
        \cdots & \cdots & P_{\theta_{0_{i,i}}} & \cdots & 0 \\
        \vdots & \vdots & \vdots & \vdots & \vdots \\
        0 & \cdots & P_{\theta_{0_{K,i}}} & \cdots & 0
    \end{pmatrix}' =  
    \begin{pmatrix}
        0 & 0 & 0 & 0 & 0 \\
        \vdots & \vdots & \vdots & \vdots & \vdots \\
        P_{\theta_{0_{i,1}}} & \cdots & P_{\theta_{0_{i,i}}} & \cdots & P_{\theta_{0_{i,K}}} \\
        \vdots & \vdots & \vdots & \vdots & \vdots \\
        0 & 0 & 0 & 0 & 0
    \end{pmatrix}.
\end{align*}
Since the vector \(\left(1/K, \ldots, 1/K\right)' \in \mathbb{R}^K\) is the unique invariant distribution for this chain, it follows by direct computation that
\begin{align*}
    P_{\theta_0}'J^{ii}P_{\theta_0}' = K^{-1}P_{\theta_0},\quad\text{for}\quad i=1,\ldots,K. 
\end{align*}
Therefore, the term $\displaystyle\sum_{r = 0}^{k-1}(A_\theta{P_\theta}')^r \{\partial A_\theta/\partial\theta_{ii}\}{P_\theta}'(A_\theta{P_\theta}')^{k-1-r}$ involved in the derivative ${\partial \psi_{k}(\theta)}/{\partial \theta_{ii}}$ in Equation (\ref{eq:9966}) can be simplified as
\begin{align*}
    \sum_{r = 0}^{k-1}\left(A_{\theta_0}P_{\theta_0}'\right)^r J^{ii} P_{\theta_0}' \left(A_{\theta_0} P_{\theta_0}'\right)^{k-1-r} &= a_{0}^{k-1} \sum_{r=0}^{k-1} (P_{\theta_0}')^r J^{ii} P_{\theta_0}' (P'_{\theta_0})^{k-1-r} \\
    &= a_{0}^{k-1} \left( J^{ii} P_{\theta_0}' + (k-2) P_{\theta_0}' J^{ii} P_{\theta_0}' \right) \\
    &= a_{0}^{k-1} \left( J^{ii} P_{\theta_0}' - P_{\theta_0}' J^{ii} P_{\theta_0}' \right) + (k-1) a_{0}^{k-1} P_{\theta_0}' J^{ii} P_{\theta_0}'.
\end{align*}
The terms \( J^{ii} P_{\theta_0}' + P_{\theta_0}' J^{ii} P_{\theta_0}' \) and \( P_{\theta_0}' J^{ii} P_{\theta_0}' \) are independent of \( k \geq 1 \), from which we deduce that \(\partial \psi_k(\theta_0)/\partial \theta_{ii}\) in Equation (\ref{eq:9966}) can be rewritten in the form
\begin{align}
\label{eq:0101}
  \frac{\partial \psi_k(\theta_0)}{\partial \theta_{ii}} = a_0^{k-1}\mathcal{Q}_1^{ii} + ka_0^{k-1}\mathcal{Q}_2^{ii} + a_0^{k}\mathcal{Q}_3^{ii},
\end{align}
where for $i=1,\dots,K$ the constants \((\mathcal{Q}_s^{ii})_{1 \leq s \leq 3}\) are independent of \(k \geq 1\) and are defined by
\begin{align*}
    \mathcal{Q}_1^{ii}: = & \mathbf{1}' \Bigg\{ \left( J^{ii} P_{\theta_0}' - P_{\theta_0}' J^{ii} P_{\theta_0}' \right) A_{\theta_0}^2 P_{\theta_0}' \left( I_K - A_{\theta_0}^{2} P_{\theta_0}' \right)^{-1} \Bigg\} V_{\theta_0} \pi_{\theta_0}',\\ 
    \mathcal{Q}_2^{ii} : = & \mathbf{1}' \Bigg\{ P_{\theta_0}' J^{ii} P_{\theta_0}' A_{\theta_0}^2 P_{\theta_0}' \left( I_K - A_{\theta_0}^2 P_{\theta_0}' \right)^{-1} \Bigg\} V_{\theta_0} \pi_{\theta_0}', \\
  \mathcal{Q}_3^{ii} : = & 
   \mathbf{1}' \Bigg\{ 2 a_{0} P_{\theta_0}' J^{ii} P_{\theta_0}' \left( I_K - A_{\theta_0}^2 P_{\theta_0}' \right)^{-1} + 2 a_{0} P_{\theta_0}' A_{\theta_0}^2 P_{\theta_0}' \left( I_K - A_{\theta_0}^2 P_{\theta_0}' \right)^{-1} J^{ii} P_{\theta_0}' \left( I_K - A_{\theta_0}^2 P_{\theta_0}' \right)^{-1} \Bigg\} V_{\theta_0} \pi_{\theta_0}'.
\end{align*}
 Following a similar line of reasoning to (\ref{eq:0101}), for $i,j=1,\dots,K$ there exist matrix coefficients $(\mathcal{R}_s^{ij})_{1\leq s \leq 3}$ and $(\mathcal{W}_s^{ii})_{1\leq s \leq 3}$ independent of $k\geq 1$ such that
\begin{align}
\label{eq:0102}
    \frac{\partial \psi_k(\theta_0)}{\partial \overline{\theta}_{ij}} &= a_0^{k-1} \mathcal{R}_1^{ij} + k a_0^{k-1} \mathcal{R}_2^{ij} + a_0^{k} \mathcal{R}_3^{ij} & \text{and} & \quad \frac{\partial \psi_k(\theta)}{\partial \tilde{\theta}_{ii}} = a_0^{k-1} \mathcal{W}_1^{ii} + k a_0^{k-1} \mathcal{W}_2^{ii} + a_0^{k} \mathcal{W}_3^{ii}.
\end{align}
In view of Equations (\ref{eq:0101}) and (\ref{eq:0102}), let us demonstrate that $J_{\Psi^N}(\theta_0)$ has a rank of at most $2$ for all $N\ge 3$, it suffices to identify constants \((\lambda_s)_{1 \leq s \leq 3}\), not all zero, such that the following system holds 
\begin{equation}
\label{eq:0103}
\left\lbrace
\begin{aligned}
	\lambda_1 a_0 + \lambda_2 a_0^2 + \lambda_3 a_0^3 = & \ 0 \\
        \lambda_1 + 2\lambda_2 a_0 + 3\lambda_3 a_0^2 = & \ 0 \\
        \lambda_1 + \lambda_2 a_0 + \lambda_3 a_0^2 = & \ 0.
\end{aligned}
 \right.
\end{equation}
By closely examining the system described by Equation (\ref{eq:0103}), we observe that this is equivalent to showing that \(a_0\) is a double root of the polynomial $P_3$ defined by \(P_3:x \mapsto \lambda_1 x + \lambda_2 x^2 + \lambda_3 x^3\). Thanks to this point of view, one can check that such a non-zero solution to the system (\ref{eq:0103}) may for example be given by 
\begin{align*}
    (\lambda_1, \lambda_2, \lambda_3) = (1, -2a_0, a_0^2).
\end{align*}
This proves that the matrix \(J_{\Psi^N}(\theta_0)\) has a rank of at most 2. Therefore, the matrix \( J_{\Psi^N}(\theta_0)'J_{\Psi^N}(\theta_0) \) is non-invertible for this particular choice of the parameter \(\theta_0\), since the first three rows of the matrix are linearly dependent. We may go even further by showing that, in the particular case when
\begin{equation}\label{Cond_J_Psi}
    a_0\in \left(\frac{1}{\sqrt{2}},1  \right)
\end{equation}
then for any \(N\ge 3\), \(J_{\Psi^N}(\theta_0)\) has a rank \textit{exactly} equal to 2 for this choice of the parameter \(\theta_0\). Indeed, for all \(m \geq 3\), let us show that \(\partial \psi_{m}(\theta_0) / \partial \theta_{ii}\) (respectively \(\partial \psi_{m}(\theta_0) / \partial \overline{\theta}_{ij}\), \(\partial \psi_{m}(\theta_0) / \partial \tilde{\theta}_{ii}\)) is a linear combination of \(\partial \psi_1(\theta_0) / \partial \theta_{ii}\) (respectively \(\partial \psi_1(\theta_0) / \partial \overline{\theta}_{ij}\), \(\partial \psi_1(\theta_0) / \partial \tilde{\theta}_{ii}\)) and \(\partial \psi_2(\theta_0) / \partial \theta_{ii}\) (respectively \(\partial \psi_2(\theta_0) / \partial \overline{\theta}_{ij}\), \(\partial \psi_2(\theta_0) / \partial \tilde{\theta}_{ii}\)) under Condition \eqref{Cond_J_Psi}.
\\
For this and  similarly to Equation (\ref{eq:0103}), it suffices to show that there exist constants \((\lambda_i)_{1\leq i\leq 3}\) such that $\lambda_3\neq 0$ and
\begin{equation}
\left\lbrace
\begin{aligned}
	\lambda_1 a_0 + \lambda_2 a_0^2 + \lambda_3 a_0^m = & \ 0 \\
        \lambda_1 + 2\lambda_2 a_0 + n\lambda_3 a_0^{m-1} = & \ 0 \\
        \lambda_1 + \lambda_2 a_0 + \lambda_3 a_0^{m-1} = & \ 0.
\end{aligned}
 \right.
\end{equation}
As in Equation (\ref{eq:0103}), this amounts to finding \(\lambda_1, \lambda_2, \lambda_3\) such that \(a_0\) is a double root of the polynomial of degree \(m\), defined as \(P_m:x \mapsto \lambda_1 x + \lambda_2 x^2 + \lambda_3 x^m\). Therefore, we seek \((\lambda_i)_{1 \leq i \leq 3}\)
and $(\nu_j)_{0\leq j\leq m-2}$, such that
\begin{align}
\label{eq:0104}
    (x - a_0)^2 \left( \sum_{i=0}^{m-2} \nu_i x^i \right) = \lambda_1 x + \lambda_2 x^2 + \lambda_3 x^m,\;\forall x\in\mathbb{R}.
\end{align}
Expanding \eqref{eq:0104} and identifying the coefficients of the polynomials yields
\begin{eqnarray}
a_0^2\nu_0 & = & 0, \label{L0}\\
- 2 a_0\nu_0 + a_0^2\nu_1 & = & \lambda_1, \label{L1}\\
-2 a_0 \nu_1 + a_0^2 \nu_2 & = & \lambda_2, \label{L2}\\
\nu_{k-2}- 2 a_0 \nu_{k-1}+ a_0^2 \nu_k & = & 0,\quad k=3,\ldots,m-2,\label{Lk}\\
\nu_{m-3}- 2 a_0 \nu_{m-2}& =& 0, \label{Lm-1}\\
\nu_{m-2}& = & \lambda_3.\label{Lm}
\end{eqnarray}
\eqref{L0} and \eqref{L1} respectively imply $\nu_0=0$ and $\lambda_1=a_0^2\nu_1$.  Now, $\nu_k$, $k=3,\ldots,m-2$, satisfies the second order recurrence relation \eqref{Lk} of which general expression can be verified to be
\begin{equation}
    \nu_k= C_1 a_0^k+ C_2k a_0^k,\quad k=1,\ldots,m-2,
\end{equation}
for some constants $C_1$ and $C_2$ that verify $2a_0\lambda_3=C_1 a_0^{m-3}+ C_2(m-3) a_0^{m-3}$ and $\lambda_3=C_1 a_0^{m-2}+ C_2(m-2) a_0^{m-2}$ thanks to \eqref{Lm-1} and \eqref{Lm},  from which one easily checks that, setting $C_1=1$ and $C_2=\frac{1/2-a_0^2}{(m-2)a_0^2-(m-3)/2}$ yields, after a bit of computation, the expressions of the coefficients
\begin{multline*}
\lambda_1=a_0^3\frac{(m-3)(a_0^2-1/2)+1/2}{(m-2) a_0^2-(m-3)/2},\quad \lambda_2= a_0^2\left[ -2+a_0^2+2(a_0^2-1) \frac{1/2-a_0^2}{(m-2)a_0^2-(m-3)/2}   \right],\\ 
    \lambda_3 = a_0^{m-2}\frac{1/2}{(m-2)a_0^2-(m-3)/2}.
\end{multline*}
Note that it is not difficult to check that $(m-2)a_0^2-(m-3)/2$ in the expression above is indeed different from $0$ for all $m\ge 3$ when Condition \eqref{Cond_J_Psi} is satisfied, so that $\lambda_3$ is well defined and different from $0$.


%

\subsection{Expression of the matrix $I$ when $(\eta_t)_{t\in\mathbb{Z}}$ is assumed i.i.d.}\label{I_ind}
\noindent The aim of this subsection is to show that the covariance matrix \(I=I^N(\theta_0)\) given in Theorem~\ref{thm:3} has an explicit, albeit not simple, expression in the particular case when the noise sequence is i.i.d. This will be important in comparing the performance of the constant estimator \(\hat{\theta}_n\) defined in Theorem \ref{thm:2} in the upcoming numerical Section \ref{sect6}, as opposed to the case where the noise $(\eta_t)_{t\in\mathbb{Z}}$ is non correlated but exhibits a dependence structure.
\\
Let \(m_1, m_2 \in \{1, \dots,N\}\). Starting from the expression for \(I\) given in \eqref{eq:000} and by stationarity of the process \((X_t)_{t\in \mathbb{Z}}\), we have
\begin{align}
\label{Im1m2}
  I(m_1, m_2) = \sum_{k=-\infty}^{\infty} \Cov\left(X_0X_{m_1}, X_{-k}X_{-k+m_2}\right).  
\end{align}
In view of Equation (\ref{eq:4}), for any \(t\) in \(\mathbb{Z}\), we may then write
\begin{align}
\label{eq:377}
    X_t = \sum_{i=0}^{\infty} d_i^t \eta_{t-i},
\end{align}
where $d_i^t := \left(\prod_{j=0}^{i-1}a(\Delta_{t-j})\right)f(\Delta_{t-i})=\left(\prod_{j=-t}^{i-t-1}a(\Delta_{-j})\right)f(\Delta_{t-i})$. 
\\
By substituting \(X_0\), \(X_{m_1}\), \(X_{-k}\), and \(X_{-k+m_2}\) into \(\Cov(X_0X_{m_1}, X_{-k}X_{-k+m_2})\), it follows from Assumption \((\mathbf{A_3})\) that
 \begin{align}
\label{eq:888}
   \Cov\left(X_0X_{m_1},X_{-k}X_{-k+m_2}\right) = & \displaystyle\sum_{i_1,\ldots,i_4 = 0}^{\infty} \mathbb{E}\left(d_{i_1}^0d_{i_2}^{-k}d_{i_3}^{m_1}d_{i_4}^{m_2-k}\right)\mathbb{E}\left(\eta_{-i_1}\eta_{-k-i_2}\eta_{m_1-i_3}\eta_{m_2-k-i_4}\right) \notag \\
   & \hspace{0.7cm} - \left(\displaystyle\sum_{i_1,i_3 = 0}^{\infty}\mathbb{E}\left(d^0_{i_1} d^{m_1}_{i_3}\right)\mathbb{E}\left(\eta_{-i_1}\eta_{m_1-i_3}\right)\right) \times \notag \\ & \hspace{0.7cm} \left(\displaystyle\sum_{i_2,i_4=0}^{\infty}\mathbb{E}\left(d^{-k}_{i_2}d^{-k+m_2}_{i_4}\right)\mathbb{E}\left(\eta_{-k-i_2}\eta_{m_2-k-i_4}\right)\right).
\end{align}
 Furthermore, since \((\eta_t)_{t\in \mathbb{Z}}\) is assumed to be i.i.d., it is possible to distinguish the cases where the different moments mentioned above namely $\mathbb{E}\left(\eta_{-i_1}\eta_{-k-i_2}\eta_{m_1-i_3}\eta_{m_2-k-i_4}\right)$, $\mathbb{E}\left(\eta_{-i_1}\eta_{m_1-i_3}\right)$ and $\mathbb{E}\left(\eta_{-k-i_2}\eta_{m_2-k-i_4}\right)$ are not zero. 
\\
More precisely, we can easily observe that,
\begin{multline*}
    \mathbb{E}\left(\eta_{-i_1} \eta_{-k-i_2} \eta_{m_1-i_3} \eta_{m_2-k-i_4}\right) 
    =  \ \mathbb{E}(\eta_0^4) \, \mathds{1}_{ [(i_1, i_2, i_3, i_4, k) \in \mathcal{A}_1]} 
     \\+ \left(\mathbb{E}(\eta_0^2)\right)^2 \Big( 
    \mathds{1}_{ [(i_1, i_2, i_3, i_4, k) \in \mathcal{A}_2]} + \mathds{1}_{ [( i_1, i_2, i_3, i_4, k) \in \mathcal{A}_3]} + \mathds{1}_{  [(i_1, i_2, i_3, i_4, k) \in \mathcal{A}_4]}
    \Big)
\end{multline*}
and
\begin{align*}
    \mathbb{E}\left(\eta_{-i_1}\eta_{m_1-i_3}\right) \mathbb{E}\left(\eta_{-k-i_2}\eta_{m_2-k-i_4}\right) 
    = & \ \left(\mathbb{E}(\eta_0^2)\right)^2 \mathds{1}_{  [(i_1, i_3) \in \mathcal{A}_5]} \mathds{1}_{  [( i_2,i_4, k) \in \mathcal{A}_6]},
\end{align*}
where the sets $\mathcal{A}_i$, $i=1,\ldots,6$, are defined as follows:
\begin{align*}
    \mathcal{A}_1 & := \{(i_1, i_2, i_3, i_4, k) \in \mathbb{N}^4 \times \mathbb{Z} : 
    -i_1 = -k-i_2 = m_1-i_3 = m_2-k-i_4\}, \\
    \mathcal{A}_2 & := \{(i_1, i_2, i_3, i_4, k) \in \mathbb{N}^4 \times \mathbb{Z} : 
    -i_1 = -k-i_2, \, m_1-i_3 = m_2-k-i_4, \, i_1 \neq i_3 - m_1\}, \\
    \mathcal{A}_3 & := \{(i_1, i_2, i_3, i_4, k) \in \mathbb{N}^4 \times \mathbb{Z} : 
    -i_1 = m_1-i_3, \, -k-i_2 = m_2-k-i_4, \, i_1 \neq k + i_2\}, \\
    \mathcal{A}_4 & := \{(i_1, i_2, i_3, i_4, k) \in \mathbb{N}^4 \times \mathbb{Z} : 
    -i_1 = m_2-k-i_4, \, -k-i_2 = m_1-i_3, \, i_1 \neq k + i_2\}, \\
    \mathcal{A}_5 & := \{(i_1, i_3) \in \mathbb{N}^2 : -i_1 = m_1-i_3\}, \\
    \mathcal{A}_6 & := \{(i_2, i_4, k) \in \mathbb{N}^2 \times \mathbb{Z} : -k-i_2 = m_2-k-i_4\}.
\end{align*}
 Thus, we obtain
\begin{align}
\label{eq:777}
& \displaystyle\sum_{i_1,\ldots,i_4 = 0}^{\infty} \mathbb{E}\left(d_{i_1}^0d_{i_2}^{-k}d_{i_3}^{m_1}d_{i_4}^{m_2-k}\right)\mathbb{E}\left(\eta_{-i_1}\eta_{-k+i_2}\eta_{m_1-i_3}\eta_{m_2-k-i_4}\right) \notag \\ & = \ \mathbb{E}(\eta_0^4)\displaystyle\sum_{i_1=0}^{\infty}\mathbb{E}\left(d^0_{i_1}d^{-k}_{-k+i_1}d^{m_1}_{i_1+m_1}d^{m_2-k}_{m_2-k+i_1}\right)  + \ \sum_{\substack{i_1, i_3=0 \\ i_1 \neq i_3-m_1}}^{\infty} \mathbb{E}\left(d^0_{i_1}d^{-k}_{i_1-k}d^{m_1}_{i_3}d^{m_2-k}_{m_2-m_1-k+i_3}\right) \notag \\ &  + \ \displaystyle\sum_{\substack{i_1,i_2=0 \\ i_1\neq k+i_2}}^{\infty}\mathbb{E}\left(d^0_{i_1}d^{-k}_{i_2}d^{m_1}_{m_1+i_1}d^{m_2-k}_{m_2+i_2}\right)   + \ \displaystyle\sum_{\substack{i_1,i_2=0 \\ i_1\neq k+i_2}}^{\infty}\mathbb{E}\left(d^0_{i_1}d^{-k}_{i_2}d^{m_1}_{m_1+k+i_2}d^{m_2-k}_{m_2-k+i_1}\right)
\end{align}
and 
\begin{align}
\label{eq:666}
  &\left(\displaystyle\sum_{i_1,i_3 = 0}^{\infty}\mathbb{E}\left(d^0_{i_1} d^{m_1}_{i_3}\right)\mathbb{E}\left(\eta_{-i_1}\eta_{m_1-i_3}\right)\right)\left(\displaystyle\sum_{i_2,i_4=0}^{\infty}\mathbb{E}\left(d^{-k}_{i_2}d^{-k+m_2}_{i_4}\right)\mathbb{E}\left(\eta_{-k-i_2}\eta_{m_2-k-i_4}\right)\right) \notag \\  & = \left(\displaystyle\sum_{i_1=0}^{\infty}\mathbb{E}\left(d^0_{i_1}d^{m_1}_{m_1+i_1}\right)\right)\left(\displaystyle\sum_{i_2=0}^{\infty}\mathbb{E}\left(d^{-k}_{i_2}d^{-k+m_2}_{m_2+i_2}\right)\right).
\end{align}
 Denoting \( c_{-t}^{i-t-1}f_{t-i} := d_i^t = \left(\prod_{j=-t}^{i-t-1}a(\Delta_{-j})\right)f(\Delta_{t-i}) \), 
and combining Equations (\ref{eq:888}), (\ref{eq:777}), and (\ref{eq:666}), 
the expression for \( I(m_1, m_2) \) simplifies as follows
\begin{align}
\label{eq:555}
  I(m_1,m_2) = & \displaystyle\sum_{k=-\infty}^{\infty} \Bigg[\mathbb{E}(\eta_0^4)\displaystyle\sum_{i_1=0}^{\infty}\mathbb{E}\left(c^{i_1-1}_0c^{i_1-1}_{k}c^{i_1-1}_{-m_1}c^{i_1-1}_{-m_2+k}f_{-i_1}^4\right) \notag\\ & +\displaystyle\sum_{\substack{i_1,i_2=0 \notag\\i_1\neq i_2-m_1}}^{\infty}\mathbb{E}\left(c^{i_1-1}_{0}c^{i_1-1}_{k}f_{-i_1}^2c^{i_2-m_1-1}_{-m_1}c^{i_2-m_1-1}_{-m_2+k}f_{m_1-i_2}^2\right) \notag \\ & +\ \displaystyle\sum_{\substack{i_1,i_2=0 \notag \\ i_1\neq k+i_2}}^{\infty}\mathbb{E}\left(c^{i_1-1}_{0}c^{i_1-1}_{-m_1}f^{2}_{-i_1}c^{i_2+k-1}_{k}c^{i_2+k-1}_{-m_2+k}f^{2}_{-i_2-k}\right) \\ & + \displaystyle\sum_{\substack{i_1,i_2=0\\ i_1\neq k+i_2}}^{\infty}\mathbb{E}\left(c^{i_1-1}_{0}c^{i_1-1}_{-m_2+k}f^{2}_{-i_1}c^{i_2+k-1}_{k}c^{i_2+k-1}_{-m_1}f^{2}_{-i_2-k}\right) \notag \\ & - \ \left(\displaystyle\sum_{i_1=0}^{\infty}\mathbb{E}\left(c^{i_1-1}_{0}c^{i_1-1}_{-m_1}f^{2}_{-i_1}\right)\right)\left(\displaystyle\sum_{i_1=0}^{\infty}\mathbb{E}\left(c^{i_1+k-1}_{k}c^{i_1+k-1}_{-m_2+k}f^{2}_{-i_1-k}\right)\right)\Bigg].
\end{align}
To conclude, we will express \( I(m_1, m_2) \) in terms of the parameters of the model (\ref{eq:1}). To do this, we will consider two different cases to express the various terms in (\ref{eq:555}) as functions of the parameters of the model (\ref{eq:1}) using Lemma \ref{lem:1}. The two following cases explain how to obtain closed form expressions for the generic quantities respectively of the form $\mathbb{E}\left(\prod_{i=1}^{s} c^{\iota_i}_{\kappa_i}f_{-\iota_i-1}\right)$ and $ \mathbb{E}\left(\prod_{i=1}^{2} c^{\iota_i}_{\kappa_i}f_{-\iota_i-1}\displaystyle\prod_{i=3}^{4}c^{\iota_i}_{\kappa_i}f_{-\iota_i-1}\right)$, $s=2,4$, $\iota_i\in \mathbb{N}$, $i\in\{1,2,3,4\}$, that appear in \eqref{eq:555}.

\subparagraph{$\diamond$ Case 1: Expression for $\mathbb{E}\left(\prod_{i=1}^{s} c^{\iota_i}_{\kappa_i}f_{-\iota_i-1}\right)$.}\ \\

\begin{figure}[H]
    \centering
    \includegraphics[width=0.8\textwidth]{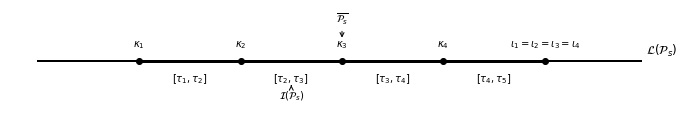}
    \caption{An illustration example of the sets \( \overline{\mathcal{P}_s} \), \( \mathcal{L}(\mathcal{P}_s) \), and \( \mathcal{I}(\mathcal{P}_s) \) for \( s=4 \). This shows how each element of \( \overline{\mathcal{P}_s} \) is positioned on the number line \( \mathcal{L}(\mathcal{P}_s) \) and the intervals between consecutive points represent the elements of \( \mathcal{I}(\mathcal{P}_s) \), demonstrating their sequential relationships and distribution.}
    \label{fig:1}
\end{figure}

 Let \( s \in \{2,4\} \) be fixed. Define the set 
\[
\mathcal{P}_s := \{(\kappa_i, \iota_i)_{i=1,\dots,s} \in (\mathbb{Z}^2)^s \mid  \iota_i = \iota_{i+1},  i =1,\dots,s-1\}.
\]
We then define \( \overline{\mathcal{P}_s} \) as follows
\[
\overline{\mathcal{P}_s} := \{\kappa_i, \iota_i, i = 1, \ldots, s \ \text{and} \  (\kappa_i, \iota_i)_{i=1,\cdots,s} \in \mathcal{P}_s\}
\]
representing all distinct individual elements extracted from each pair in \( \mathcal{P}_s \). This construction ensures that \( \overline{\mathcal{P}_s} \) contains only unique values from both components of the pairs.
\\
 Let \( \mathcal{L}(\mathcal{P}_s) \) be the set of sorted elements of  \( \overline{\mathcal{P}_s} \) on the real line.
%
  We then define the set of intervals \( \mathcal{I}(\mathcal{P}_s) \) as follows 
\[
\mathcal{I}(\mathcal{P}_s) := \{[\tau_i, \tau_{i+1}], i=1,\dots,s, \tau_i \leq \tau_{i+1}\ \text{and}\ \tau_i, \tau_{i+1} \in \mathcal{L}(\mathcal{P}_s) \cap \overline{\mathcal{P}}_s\}.
\]
which represents a set of intervals consisting of which endpoints are the consecutive elements of \( \mathcal{L}(\mathcal{P}_s) \).
An illustrative example for \( s=4 \) of the sets \( \overline{\mathcal{P}_s} \), \( \mathcal{L}(\mathcal{P}_s) \), and \( \mathcal{I}(\mathcal{P}_s) \) is presented below (see Fig.~\ref{fig:1}). Each point represents an element of \( \overline{\mathcal{P}_s} \) placed on the number line \( \mathcal{L}(\mathcal{P}_s) \), and the intervals between consecutive points illustrate the elements of \( \mathcal{I}(\mathcal{P}_s) \).

Thus, in view of this modeling and utilizing Lemma \ref{lem:1}, we have

\begin{align}
\label{eq:444}
\mathbb{E}\left(\prod_{i=1}^{s} c^{\iota_i}_{\kappa_i}f_{-\iota_i-1}\right) = & \ \bm{1}^\prime\left(\prod_{i=1}^{s} Q_{a^{\varphi([\tau_i,\tau_{i+1}])}}^{\tau_{i+1}-\tau_i}\right) \bm{\pi}_{f^s} \notag\\
= & \ \bm{1}^\prime\left(\prod_{\zeta \in \mathcal{I}(\mathcal{P}_s)} Q_{a^{\varphi(\zeta)}}^{\length\zeta}\right) \bm{\pi}_{f^s}
\end{align}
where  $Q_{a^{\nu_1}} = A_{\theta_0}^{\nu_1} P_{\theta_0}'$, ${Q}_{{f}^{s}} = \text{diag}\left({f}^{s}(1),\ldots,{f}^{s}(K)\right) P_{\theta_0}'$  and the column vector $\bm\pi_{f^{s}} = \left(f^{s}(1)\pi(1),\ldots,f^{s}(K)\pi(K)\right)'$ for some $\nu_1>0$ and where
\begin{align*}
    \varphi([\tau_i, \tau_{i+1}])  := \sum_{j=1}^{s} \mathds{1}_{[\kappa_j, \iota_j] \supset [\tau_i, \tau_{i+1}]}, \ 1 \leq i < s, \ \varphi(\zeta)  := \sum_{j=1}^{s} \mathds{1}_{[\kappa_j, \iota_j] \supset \zeta}, \ \zeta \in \mathcal{I}(\mathcal{P}_s) \ \text{and} 
    \length\zeta
\end{align*} \text{denotes the length of an element $\zeta\in$ } $\mathcal{I}(\mathcal{P}_s)$.
\subparagraph{$\diamond$ Case 2: Expression for $\mathbb{E}\left(\prod_{i=1}^{2} c^{\iota_i}_{\kappa_i}f_{-\iota_i-1}\displaystyle\prod_{i=3}^{4}c^{\iota_i}_{\kappa_i}f_{-\iota_i-1}\right)$.} \ \\
 Let \(\mathcal{H}\) be defined as \begin{align*}\mathcal{H} := \{(\kappa_i,\iota_i)_{i=1,\dots,4} \in (\mathbb{Z}^2)^4 \mid\  \iota_1 = \iota_2, \ \iota_3 = \iota_4 \text{ and } \iota_2 \neq \iota_3\}.\end{align*} As previously, we define \(\overline{\mathcal{H}}\) as \begin{align*} \overline{\mathcal{H}} := \{\kappa_i,\iota_i, i=1,\dots,4 \ \text{and} \ (\kappa_i, \iota_i)_{i=1,\dots,4} \in \mathcal{H}\},\end{align*} representing all distinct individual elements extracted from each pair in \(\mathcal{H}\). 
\\
Let \( \mathcal{L}(\mathcal{H}) \) be the set of sorted elements of  \( \overline{\mathcal{H}} \) on the real line. Next, we form 
\[
\mathcal{J}(\mathcal{H}) := \{[\tau_{i}, \tau_{i+1}] \mid i = 1, \cdots, 4, \ \tau_{i} \leq \tau_{i+1}, \ \tau_{i}, \tau_{i+1} \in \mathcal{L}(\mathcal{H}) \cap \overline{\mathcal{H}}\},
\]
a set of intervals consisting of consecutive points on this line, where each interval is formed between consecutive entries in the representation \(\mathcal{L}(\mathcal{H})\). With this consideration and in view of Lemma~\ref{lem:1}, it also follows that

\begin{align}
\label{eq:4444}
\mathbb{E}\left(\prod_{i=1}^{2} c^{\iota_i}_{\kappa_i}f_{-\iota_i-1}\displaystyle\prod_{i=3}^{4}c^{\iota_i}_{\kappa_i}f_{-\iota_i-1}\right) = &  \ 
\bm{1}^\prime\left(\prod_{\zeta \in \mathcal{J}(\mathcal{H})\setminus{\{\zeta_{\star}\}}} Q_{a^{\varphi(\zeta)}}^{\length\zeta}\right)Q_{f^2} Q_{a^{\varphi(\zeta_{\star})}}^{\length{\zeta_{\star}}^{}}\bm{\pi}_{f^2}
\end{align}
where $\zeta_{\star}$ represents the interval with the highest index of $\mathcal{J}(\mathcal{H})$, $\varphi(\zeta) := \sum_{j=1}^{4} \mathds{1}_{[\kappa_j, \iota_j] \supset \zeta}$, $\zeta\in\mathcal{J}(\mathcal{H})$ and $\length \zeta$ denotes the length of an element $\zeta\in\mathcal{J}(\mathcal{H})$.
\\
Finally, the computations carried out in the above Cases 1 and 2 yield the existence of a family of sets $\left(\mathcal{I}_{i}\right)_{1\leq i \leq 3}$, $\left(\mathcal{J}_j\right)_{1\leq j \leq 3}$ such that $I(m_1,m_2)$ in \eqref{Im1m2} can be expressed, in view of Equations (\ref{eq:444}) and (\ref{eq:4444}), as
\begin{multline*}
   I(m_1,m_2) = \displaystyle\sum_{k=-\infty}^{\infty}\left[\mathbb{E}(\eta_0^4)\displaystyle\sum_{i_1=0}^{\infty} \bm{1}^\prime\left(\prod_{\zeta \in \mathcal{I}_1} Q_{a^{\varphi(\zeta)}}^{\length\zeta}\right) \bm{\pi}_{f^4} \right. \\ + \displaystyle\sum_{\substack{i_1,i_2=0 \notag \\ i_1\neq i_2-m_1}}^{\infty}\bm{1}^\prime\left(\prod_{\zeta \in \mathcal{J}_1\setminus{\{\zeta_{\star}\}}} Q_{a^{\varphi(\zeta)}}^{\length\zeta^{}}Q_{f^2} Q_{a^{\varphi(\zeta_{\star})}}^{\length{\zeta_{\star}}}\right) \bm{\pi}_{f^2}  \\  + \displaystyle\sum_{\substack{i_1,i_2=0 \\ i_1\neq k+i_2}}^{\infty}\bm{1}^\prime\left(\prod_{\zeta \in \mathcal{J}_2\setminus{\{\zeta_{\star}\}}} Q_{a^{\varphi(\zeta)}}^{\length\zeta}Q_{f^2} Q_{a^{\varphi(\zeta_{\star})}}^{\length{\zeta_{\star}}^{}}\right) \bm{\pi}_{f^2} \ + \displaystyle\sum_{\substack{i_1,i_2=0\\ i_1\neq k+i_2}}^{\infty}\bm{1}^\prime\left(\prod_{\zeta \in \mathcal{J}_3\setminus{\{\zeta_{\star}\}}} Q_{a^{\varphi(\zeta)}}^{\length\zeta^{}}Q_{f^2} Q_{a^{\varphi(\zeta_{\star})}}^{\length{\zeta_{\star}}^{}}\right) \bm{\pi}_{f^2} \\  \ -\left(\  \displaystyle\sum_{i_1=0}^{\infty} \bm{1}^\prime\left(\prod_{\zeta \in \mathcal{I}_2}Q_{a^{\varphi(\zeta)}}^{\length\zeta^{}}\right) \bm{\pi}_{f^2}\right)  \ \times \left.\left(\  \displaystyle\sum_{i_1=0}^{\infty} \bm{1}^\prime\left(\prod_{\zeta \in \mathcal{I}_3}Q_{a^{\varphi(\zeta)}}^{\length\zeta^{}}\right) \bm{\pi}_{f^2}\right)\right]
\end{multline*}
where the sets are more precisely defined as
\begin{align*}
\mathcal{I}_1 & := \mathcal{I}(\{(0,i_1-1),(k,i_1-1),(-m_1,i_1-1),(-m_2+k,i_1-1)\}), \\
\mathcal{J}_1 & := \mathcal{J}(\{(0,i_1-1),(k,i_1-1),(-m_1,i_2-m_1-1),(-m_2+k,i_2-m_1-1)\}), \\
\mathcal{J}_2 & := \mathcal{J}(\{(0,i_1-1),(k,i_2+k-1),(-m_1,i_1-1),(-m_2+k,i_2-1)\}), \\
\mathcal{J}_3 & := \mathcal{J}(\{(0,i_1-1),(k,i_2+k-1),(-m_1,i_2+k-1),(-m_2+k,i_1-1)\}), \\
\mathcal{I}_2 & := \mathcal{I}(\{(0,i_1-1),(-m_1,i_1-1)\}), \\
\mathcal{I}_3 & := \mathcal{I}(\{(k,i_1+k-1),(-m_2+k,i_1+k-1)\}).
\end{align*}
Notice that, under Assumption $(\mathbf{A_2})$, the respective spectral radii of the matrices $Q_{a^{\varphi(\zeta)}}$, for $\zeta \in \bigcup_{i=1}^{3}\mathcal{I}_i\cup\mathcal{J}_i$, are strictly less than $1$.
\\
Note also that, in practice the infinite sums involved in $I(m_1,m_2)$ are truncated.
\qed
\section{Estimation of the asymptotic covariance matrix}\label{sect5}
This section aims to propose a consistent estimator for the variance-covariance matrix $\Omega$ obtained in Theorem~\ref{thm:3}.  It is about proposing a consistent estimator of the matrix $I$ as well as the matrix $J$. For the matrix $J$, a simple estimator in the context of our study is given by
\begin{align}
\label{eq:5.1}
    \hat{J}_n: = J_{\mathcal{F}^{N,n}}(\hat{\theta}_n),
\end{align}
where $\hat{\theta}_n$ represents an estimator of $\theta_0$. However, estimating the matrix $I$ turns out to be more complex than estimating the matrix $J$. Various approaches can be considered for estimating \(I\): a non-parametric kernel estimation (see \cite{andrews1991heteroskedasticity} and \cite{newey1987hypothesis} for general references) as well as a spectral density-based estimation (see \cite{berk1974consistent} and \cite{den1996practitioner} for general references). In this paper, we focus on an estimator based on spectral density by interpreting $(2\pi)^{-1}I$ as the spectral density of the stationary process $(\mathcal{Y}_t:= Y_t(\theta_0)-\mathbb{E}[Y_t(\theta_0)])_{t\in\mathbb{Z}}$ evaluated at frequency zero (see \citet[ p. 459]{brockwell1991time}).
 A similar approach to estimate the matrix $I$ can be found in \citet[Theorem 3.10, p. 10]{boubacar2020estimation}. This technique involves writing the matrix \(I\) as:
\begin{align*}
    I = \bm{\varphi}(1)^{-1}\Sigma_u\bm{\varphi}'(1)^{-1},
\end{align*}
when $(\mathcal{Y}_t)_{t\in\mathbb{Z}}$ exhibits an AR$(\infty)$ structure 
\begin{align}
\label{eq:20}
    \bm{\varphi}(L)\mathcal{Y}_t: = \mathcal{Y}_t -\sum_{i=1}^{\infty}\varphi_i\mathcal{Y}_{t-i} = u_t,\quad\text{with }\, \bm{\varphi}(L)= I_N -\sum_{i=1}^{\infty}\varphi_i L^i
\end{align}
where $(u_t)_{t\in\mathbb{Z}} \in \mathbb{R}^{N}$ is a weak white noise with variance-covariance matrix $\Sigma_u$, $L$ stands for the back-shift operator, and $I_N$ is the identity operator. 
Even though the sequence \((X_t)_{t = 1, \dots, n}\) is observable, \(\mathcal{Y}_t\) is not observable because \(\mathbb{E}[Y_t(\theta_0)]\) is unknown. An estimator \(\hat{\mathcal{Y}}_t\) of \(\mathcal{Y}_t\) is thus obtained by replacing \(\mathbb{E}[Y_t(\theta_0)]\) by its empirical estimator \((\hat{c}_{1,0}, \dots, \hat{c}_{N,0})'\) in the expression of \(\mathcal{Y}_t\), so that
$$\hat{\mathcal{Y}}_t=Y_t-(\hat{c}_{1,0}, \dots, \hat{c}_{N,0})',\quad t=1,\ldots,n .$$
We also define $\hat{\varphi}_{r,1},\dots,\hat{\varphi}_{r,r}$ as the coefficients of the regression of $\hat{\mathcal{Y}}_t$ on $\{\hat{\mathcal{Y}}_{t-1},\dots,\hat{\mathcal{Y}}_{t-r}\}$, $\hat{u}_{r,t}$ as the residual from this regression and $\hat{\Sigma}_{\hat{u}_r}$ as the covariance matrix of the residuals $\hat{u}_{r,1},\dots,\hat{u}_{r,n}$. Formally, $\hat{\mathcal{Y}}_t$ obeys to the equation 
\begin{align}
\label{eq:21}
   \hat{\bm{\varphi}}(L)\hat{\mathcal{Y}}_t : = \hat{\mathcal{Y}}_t - \sum_{i=1}^{r}\hat{\varphi}_{r,i}\hat{\mathcal{Y}}_{t-i} = \hat{u}_{r,t}.
\end{align}
The asymptotic study of the estimator of $I$ using the spectral density method is given in the following theorem.
\begin{t1}\label{thm:4}
Let the conditions of Theorem~\ref{thm:3} be satisfied. Additionally, we assume that $\mathbb{E}|\eta_t|^{8+4\nu}<\infty$ for some $\nu>0$ and the process $(\mathcal{Y}_t)_{t\in\mathbb{Z}}$ has an $AR(\infty)$ representation as specified in Equation $(\ref{eq:20})$. Moreover, suppose that $\|\varphi_i\| = \mathrm{o}(i^{-2})$ as $i\rightarrow\infty$, the roots of $\det(\bm{\varphi}(z)) = 0,\;z\in\mathbb{C}$, are outside the unit disk and the matrix $\Sigma_u$ is non-singular. Under these conditions, the spectral estimator of the matrix $I$ holds:
\begin{align*}
    \hat{I}^{SP} := \hat{\bm{\varphi}}(1)^{-1}\hat{\Sigma}_{\hat{u}_r}{\hat{\bm{\varphi}}'(1)^{-1}}
\end{align*}
converges in probability to $I= {\bm{\varphi}}(1)^{-1}{\Sigma}_{{u}}{{\bm{\varphi}}'(1)^{-1}}$ when $r = r(n)\rightarrow\infty$ and $r = \mathrm{o}(n^{1/3})$ as $n\to\infty$. 
\end{t1}
The proof of this theorem is given in Section \ref{proof_thm4}.

Consequently a weakly consistent estimator of $\Omega$ is
\begin{align*}
   \hat{\Omega} := \hat{M}_n^{-1}\hat{J}_{n}'\hat{I}^{SP}\hat{J}_{n}'\hat{M}_n^{-1},
\end{align*} where $\hat{M}_n := \hat{J}_n'\hat{J}_n$ with $\hat{J}_n$ defined in Equation (\ref{eq:5.1}).
\\ Let $\hat{Q}_{\hat{a}^{\nu_1}} = A_{\hat{\theta}_n}^{\nu_1} P_{\hat{\theta}_n}'$, $\hat{Q}_{\hat{f}^{s}} = \text{diag}\left(\hat{f}^{s}(1),\ldots,\hat{f}^{s}(K)\right) P_{\hat{\theta}_n}'$ and the column vector $\hat{\bm\pi}_{\hat{f}^{s}} = \left(\hat{f}^{s}(1)\hat{\pi}(1),\ldots,\hat{f}^{s}(K)\hat{\pi}(K)\right)'$ for some $\nu_1>0$ and $s\in\{2,4\}$.
In the standard strong ARHMC case i.e when the noise \((\eta_t)_{t \in \mathbb{Z}}\) is independent (particularly when $\eta_t\stackrel{\mathcal{D}}{=}\mathcal{N}(0,1)$), in view of Section \ref{I_ind}, we have \(\hat{\Omega}_S := \hat{M}_n^{-1}\hat{J}_{n}'\hat{I}_{S}\hat{J}_{n}'\hat{M}_n^{-1}\)
where \(\hat{I}_S\) is a consistent estimator of  the matrix \(I_S\) defined for a fixed integers $r_, r_2>0$ as:

\begin{multline*}
   \hat{I}_S(m_1,m_2) = \displaystyle\sum_{k=-r_1}^{r_1}\left[3\displaystyle\sum_{i_1=0}^{r_2} \bm{1}^\prime\left(\prod_{\zeta \in \mathcal{I}_1} \hat{Q}_{\hat{a}^{\varphi(\zeta)}}^{\length\zeta}\right) \hat{\bm{\pi}}_{\hat{f}^4} \right. \\ + \displaystyle\sum_{\substack{i_1,i_2=0 \notag \\ i_1\neq i_2-m_1}}^{r_2}\bm{1}^\prime\left(\prod_{\zeta \in \mathcal{J}_1\setminus{\{\zeta_{\star}\}}} \hat{Q}_{\hat{a}^{\varphi(\zeta)}}^{\length\zeta^{}}\hat{Q}_{\hat{f}^2} \hat{Q}_{\hat{a}^{\varphi(\zeta_{\star})}}^{\length{\zeta_{\star}}}\right) \hat{\bm{\pi}}_{\hat{f}^2}  \\  + \displaystyle\sum_{\substack{i_1,i_2=0 \\ i_1\neq k+i_2}}^{r_2}\bm{1}^\prime\left(\prod_{\zeta \in \mathcal{J}_2\setminus{\{\zeta_{\star}\}}} \hat{Q}_{\hat{a}^{\varphi(\zeta)}}^{\length\zeta}\hat{Q}_{\hat{f}^2} \hat{Q}_{\hat{a}^{\varphi(\zeta_{\star})}}^{\length{\zeta_{\star}}^{}}\right) \hat{\bm{\pi}}_{\hat{f}^2} \ + \displaystyle\sum_{\substack{i_1,i_2=0\\ i_1\neq k+i_2}}^{r_2}\bm{1}^\prime\left(\prod_{\zeta \in \mathcal{J}_3\setminus{\{\zeta_{\star}\}}} \hat{Q}_{\hat{a}^{\varphi(\zeta)}}^{\length\zeta^{}}\hat{Q}_{\hat{f}^2} \hat{Q}_{\hat{a}^{\varphi(\zeta_{\star})}}^{\length{\zeta_{\star}}^{}}\right) \hat{\bm{\pi}}_{\hat{f}^2} \\  \ -\left(\  \displaystyle\sum_{i_1=0}^{r_2} \bm{1}^\prime\left(\prod_{\zeta \in \mathcal{I}_2}\hat{Q}_{\hat{a}^{\varphi(\zeta)}}^{\length\zeta^{}}\right) \hat{\bm{\pi}}_{\hat{f}^2}\right)  \ \times \left.\left(\  \displaystyle\sum_{i_1=0}^{r_2} \bm{1}^\prime\left(\prod_{\zeta \in \mathcal{I}_3} \hat{Q}_{\hat{a}^{\varphi(\zeta)}}^{\length\zeta^{}}\right) \hat{\bm{\pi}}_{\hat{f}^2}\right)\right]
\end{multline*}
with 
for \(m_1, m_2 \in \{1, \dots,N\}\) and where the sets $\mathcal{I}_1, \mathcal{I}_2, \mathcal{I}_3, \mathcal{J}_1,
\mathcal{J}_2$ and $\mathcal{J}_3$ are  defined in Section \ref{I_ind}
\section{Numerical illustrations }\label{sect6}
In this section, we investigate the finite sample properties of the asymptotic results that we introduced in this work. For that sake we use Monte Carlo experiments. The numerical
illustrations of this section are made with the \texttt{Python} software. 
\\
We examine a specific case of the model presented in Equation (\ref{eq:1}) by choosing the number of regimes $K=2$, which significantly reduces the number of parameters of the model (\ref{eq:1}). In this two-regime configuration, the total number of parameters $K^2+K$ is then $6$. Fig. \ref{fig: simulation_X} below illustrates the evolution of the process $(X_t)_{t\in\mathbb{Z}}$ under the influence of different types of noises $(\eta_t)_{t\in\mathbb{Z}}$, particularly in the cases of strong and weak white noises. To compare our results, we used the same initial parameter values $\theta_0$ as the one used by \cite{xie2008general} who studied a generalized version of the model (\ref{eq:1}) under the strong noise assumption. The linear innovation $(\eta_t)_{t\in\mathbb{Z}}$ is a function of a process $(u_t)_{t\in\mathbb{Z}}$, simulated according to a standard normal distribution $(u_t\stackrel{\mathcal{D}}{=}\mathcal{N}(0,1))$.
Table \ref{tab:noise_cases} summarizes the different noise cases used in the study, detailing their mathematical expressions and providing brief descriptions of their characteristics.

\begin{table}[H]
\centering
\renewcommand{\arraystretch}{1.5} 
\setlength{\tabcolsep}{5pt} 
\begin{tabular}{llp{7cm}}
\toprule
\toprule
\textbf{Noise type} & \textbf{Expression for $\eta_t$} & \textbf{Description} \\
\midrule
\text{Strong} \label{strong}   & $\eta_t = u_t$ & Basic noise model, considering independent and identically distributed random variables. \\
\text{Weak 1} \label{weak1}   & $\eta_t = u_t u_{t-1}$ & Weak noise with a dependence on the previous observation $u_{t-1}$. \\
\text{Weak 2} \label{weak2}   & $\eta_t = u_t^2 u_{t-1}$ & Quadratic dependence on the current value $u_t$ and linear dependence on the previous value $u_{t-1}$. \\
\text{Weak 3} \label{weak3}   & $\eta_t = u_t \left( |u_{t-1}| + 1 \right)^{-1}$ & Another weak noise case with inverse scaling by $u_{t-1}$ to reduce the impact of previous values. \\
\text{GARCH} \label{garch}    & $\eta_t = h_t^{1/2} u_t, h_t:=\omega_0 + a_0 \eta_{t-1}^2 + \beta_0 h_{t-1}$ & GARCH model incorporating volatility dynamics with conditional heteroscedasticity. \\
\bottomrule
\bottomrule
\end{tabular}
\caption{Different cases of noise experimented in the study.}
\label{tab:noise_cases}
\end{table}

Also note that the noises defined by Weak 1, Weak 2 and Weak 3 are direct generalizations of the weak white noises defined  by \citet[Example 2.1 and 2.2]{romano1996inference}. Consequently, it is straightforward to verify that they meet the criteria for weak white noises. Contrary to Weak 1, Weak 3 and GARCH, the Weak 2 noise is not a martingale difference sequence for which the limit theory is more classical.


We conducted simulations by generating \( R = 1,000 \) independent trajectories, for each of two series of length \( n = \{300;2,000\} \), based on the model described in Equation (\ref{eq:1}). The simulations are carried out to highlight the different types of noise defined in Table \ref{tab:noise_cases} in order to illustrate a range of scenarios. For each experiment, \( R \) independent realizations were generated and 
we estimated the coefficient vector $\theta_0 := \left(\alpha_{11}, \alpha_{22}, \beta_{11}, \beta_{21}, \gamma_{11}, \gamma_{22}\right)' = \left(-0.4, 0.3, 0.3, 0.2, 1.0, 0.5\right)'$. The parameter space $\Theta$ associated is chosen to satisfy the assumptions of Theorem \ref{thm:3}. The simulation procedure was as follows: starting with $\theta_0$, we simulated $R$ trajectories based on a noise type given in Table \ref{tab:noise_cases}. For each simulated trajectory, we use the estimation function $\mathcal{F}^{N,n}(\cdot)$ to generate an estimate of $\theta_0$.

 Tables \ref{tab:strong} through \ref{tab:weak_4} presented below summarize the statistical characteristics of the simulations conducted using model (\ref{eq:1}) with the various noises types defined in the Table \ref{tab:noise_cases}, thereby providing an overview of the distribution of the  estimator \(\hat{\theta}_n\). More precisely, for each element of \(\hat{\theta}_n\) they detail: the mean, representing the average value observed throughout the simulations; the standard deviation (\emph{Std}), the minimum (\emph{Min}) and maximum (\emph{Max}) values, highlighting the dataset's range. Additionally, the tables include the first (\(Q_1\)) and third (\(Q_3\)) quartiles, offering a deeper insight into the data distribution by showing the values below which a certain percentage of the data falls.
As expected, Tables \ref{tab:strong} through \ref{tab:weak_4} show that the bias and the RMSE decrease when the size of the sample increases.

Fig. \ref{fig: density} and Fig. \ref{fig: qqplot} compare the distribution of the moments estimator in the strong and weak noises cases. The distributions of \(\hat{\beta}_{11} \) and \(\hat{\beta}_{22}\) are similar in three cases (Strong ARHMC,  Weak 2 ARHMC and GARCH ARHMC) and they are more accurate than in Weak 1 ARHMC case. Whereas the moments estimator of \(\hat{\alpha}_{11}\), \(\hat{\alpha}_{22}\),  \(\hat{\gamma}_{11}\) and \(\hat{\gamma}_{22}\) are more accurate in the strong case  than in the Weak 1, Weak 2 and GARCH cases. 
This is in accordance with the results of \cite{romano1996inference} who showed that,
with similar noises, the asymptotic covariance of the sample autocorrelations can be greater (for Weak 1 or Weak 2 noises) or less (for Weak 3 noise) than 1 as well (1 is the asymptotic covariance for a strong  noise).

Fig. \ref{fig: omega} below compares the standard sandwich estimator \(\hat{\Omega}_S = \hat{M}_n^{-1}\hat{J}_{n}'\hat{I}_{S}\hat{J}_{n}'\hat{M}_n^{-1}\) 
 and our general estimator \(\hat{\Omega} = \hat{M}_n^{-1}\hat{J}_{n}'\hat{I}^{SP}\hat{J}_{n}'\hat{M}_n^{-1}\) introduced in Section~\ref{sect5}. For the calculation of \(\hat{I}^{SP}\), we used the \texttt{statsmodels} function from the \texttt{Python} package \textbf{VAR}. The order \(r\) of the AR model is automatically selected by AIC (Akaike Information Criterion). 

In the strong case we know that the two estimators are consistent. As shown in the two top panels of Fig. \ref{fig: omega}, the standard sandwich estimator \(\hat{\Omega}_S\) is more precise than \(\hat{\Omega}\) in the strong case, where it exhibits less bias and better accuracy. Whereas when examining the weak cases (Weak 1 and Weak 2), \(\hat{\Omega}_S\) performs poorly. In contrast, the sandwich estimator \(\hat{\Omega}\) proves to be much more robust across in all scenarios, although it may be slightly less precise in the strong case than \(\hat{\Omega}_S\), it remains consistent and performs well in both weak cases (see the middle and bottom subfigures of Fig. \ref{fig: omega}). More precisely, it is clear that in the weak cases  $n(\hat{\theta}_n - \theta_0)^2$ is better estimated by \(\text{diag}(\hat{\Omega})\) (see the box-plots
(a)-(f) of the right center-bottom and the right-bottom panel of Fig. \ref{fig: omega}) than by \(\text{diag}(\hat{\Omega}_S)\) (see the box-plots (a)-(f) of the left center-bottom and the left-bottom panel of Fig. \ref{fig: omega}). The failure of the standard estimator of $\Omega$ in the weak ARHMC setting may have important consequences in terms of hypothesis testing for instance.

\begin{figure}[H]
    \centering
    \includegraphics[width=0.95\textwidth]{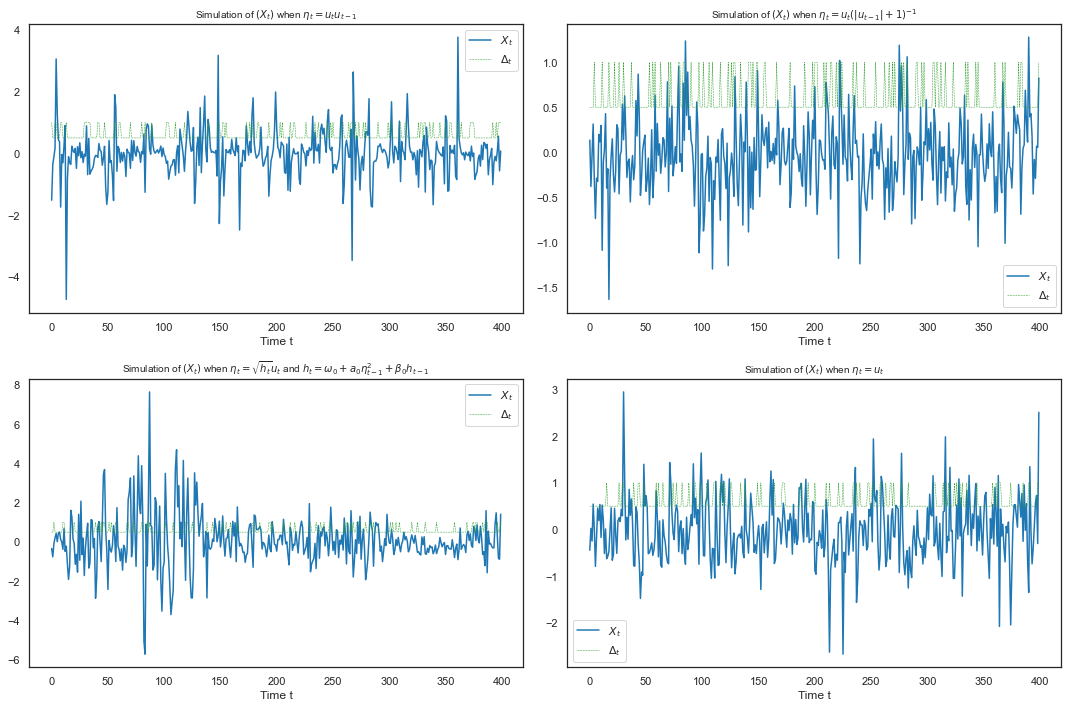}
    \caption{Simulation of length 400 of model (\ref{eq:1}) with $\theta_0: = \left(\alpha_{11},\alpha_{22},\beta_{11},\beta_{21},\gamma_{11},\gamma_{22}\right)' = \left(-0.4,0.3,0.3,0.2,1.0,0.5\right)'$ and $(\omega_0, a_0, \beta_0) = (0.2,0.1,0.5)$.}
    \label{fig: simulation_X}
\end{figure}

\begin{table}[H]
\centering
\setlength{\tabcolsep}{10pt}
\begin{tabular}{lccccccc}
\toprule
&  & $\alpha_{11}$ & $\alpha_{22}$ & $\beta_{11}$ & $\beta_{21}$ & $\gamma_{11}$ & $\gamma_{22}$  \\
\cmidrule(lr){3-8}
& $\theta_0$ & -0.4 & 0.3 & 0.3 & 0.2 & 1.0 & 0.5 \\
\cmidrule(lr){3-8}
& Min & -1.4109 & -0.28563 & 0.01091 & 0.01936 & -0.53100 & -0.08625  \\ 
& $Q_1$ & -0.48516 & 0.25538 & 0.25817 & 0.15796 & 0.95890 & 0.45766 \\
& Mean & -0.39240 & 0.30055 & 0.30457 & 0.19217 & 0.98579 & 0.48383 \\ 
& Rmse & 0.22969 & 0.09666 & 0.100314 & 0.05927 & 0.10078 & 0.06969 \\ 
\multicolumn{1}{c}{$n = 300$} & Bias & 0.00759 & 0.00055 & 0.00457 & -0.00782 & -0.01420 & -0.01616 \\
& $Q_2$ & -0.39768 & 0.30294 & 0.29809 & 0.19146 & 0.99388 & 0.49004 \\ 
& Std & 0.22968 & 0.09670 & 0.10026 & 0.05878 & 0.09982 & 0.06782 \\ 
& $Q_3$ & -0.30108 & 0.35261 & 0.33874 & 0.21860 & 1.02265 & 0.51514 \\ 
& Max& 1.14141 & 0.60950 & 0.81411 & 0.49516 & 1.33618 & 1.10513 \\ 
\cmidrule(lr){3-8}
& Min & -1.48591 & -0.00907 & 0.04130 & 0.00701 & 0.24168 & 0.12564 \\
& $Q_1$ & -0.43255 & 0.28103 & 0.28366 & 0.18570 & 0.98196 & 0.48412  \\
& Mean & -0.39784 & 0.29772 & 0.29956 & 0.20298 & 0.99339 & 0.49566 \\ 
& Rmse & 0.15457 & 0.05665 & 0.06168 & 0.04542 & 0.06052 & 0.04564 \\ 
\multicolumn{1}{c}{$n = 2,000$} & Bias & 0.00215 & -0.00227 & -0.00043 & 0.00298 & -0.00660 & -0.00433 \\
& $Q_2$ & -0.40017 & 0.29947 & 0.29912 & 0.19997 & 0.99916 & 0.49750 \\ 
& Std & 0.15463 & 0.05663 & 0.06171 & 0.04535 & 0.06018 & 0.04546 \\ 
& $Q_3$ & -0.36289 & 0.31741 & 0.31488 & 0.21367 & 1.01299 & 0.51153 \\ 
& Max& 0.69552 & 0.61279 & 0.73094 & 0.63367 & 1.39112 & 0.69572 \\
\bottomrule
\end{tabular}
\caption{Summary statistics of parameters for model (\ref{eq:1}) with noise \(\eta_t = u_t\), based on 1000 replications for sequence sizes of 300 and 2000, respectively.}
\label{tab:strong}
\end{table}

\begin{table}[H]
\centering
\setlength{\tabcolsep}{10pt}
\begin{tabular}{lccccccc}
\toprule
&  & $\alpha_{11}$ & $\alpha_{22}$ & $\beta_{11}$ & $\beta_{21}$ & $\gamma_{11}$ & $\gamma_{22}$  \\
\cmidrule(lr){3-8}
& $\theta_0$ & -0.4 & 0.3 & 0.3 & 0.2 & 1.0 & 0.5 \\
\cmidrule(lr){3-8}
& Min & -1.62424 & -0.28668 & 0.00115 & 0.00456 & 0.26374 & -0.22960\\ 
& $Q_1$ & -0.52372 & 0.20381 & 0.22390 & 0.12777 & 0.92073 & 0.42161  \\
& Mean & -0.38442 & 0.27706 & 0.28836 & 0.18181 & 0.97581 & 0.45942 \\ 
& Rmse & 0.29932 & 0.13104 & 0.11636 & 0.08371 & 0.11454 & 0.09859 \\ 
\multicolumn{1}{c}{$n = 300$} & Bias & 0.01558 & -0.02293 & -0.01163 & -0.01818 & -0.02418 & -0.04057 \\
& $Q_2$ & -0.39040 & 0.28515 & 0.28478 & 0.18034 & 0.98602 & 0.47409 \\ 
& Std & 0.29907 & 0.12908 & 0.11584 & 0.08176 & 0.11201 & 0.08990 \\ 
& $Q_3$ & -0.23993 & 0.35569 & 0.34196 & 0.22443 & 1.03396 & 0.51061 \\ 
& Max & 1.07378 & 0.70402 & 0.89723 & 0.85197 & 1.60616 & 0.81016 \\ 
\cmidrule(lr){3-8}
& Min & -1.20361 & -0.30287 & 0.00536 & 0.02943 & 0.55062 & 0.06970 \\
& $Q_1$ & -0.45436 & 0.26146 & 0.26839 & 0.17190 & 0.97583 & 0.47427  \\
& Mean & -0.38610 & 0.29112 & 0.29021 & 0.19766 & 0.99631 & 0.49304 \\ 
& Rmse & 0.15457 & 0.05665 & 0.06168 & 0.04542 & 0.06052 & 0.04564 \\ 
\multicolumn{1}{c}{$n = 2,000$} & Bias & 0.00215 & -0.00227 & -0.00043 & 0.00298 & -0.00660 & -0.00433 \\
& $Q_2$ & -0.40012 & 0.2968 & 0.29639 & 0.19558 & 0.9980 & 0.49692 \\ 
& Std & 0.19234 & 0.07453 & 0.07750 & 0.05539 & 0.06850 & 0.05259 \\ 
& $Q_3$ & -0.33254 & 0.32352 & 0.31871 & 0.21692 & 1.01961 & 0.51592 \\ 
& Max & 1.12381 & 0.59356 & 0.87746 & 0.83714 & 1.42285 & 0.66952 \\
\bottomrule
\end{tabular}
\caption{Summary statistics of parameters for model (\ref{eq:1}) with noise \(\eta_t = u_tu_{t-1}\), based on 1000 replications for sequence sizes of 300 and 2000, respectively.}
\label{tab:weak_1}
\end{table}

\begin{table}[H]
\centering
\setlength{\tabcolsep}{10pt}
\begin{tabular}{lccccccc}
\toprule
&  & $\alpha_{11}$ & $\alpha_{22}$ & $\beta_{11}$ & $\beta_{21}$ & $\gamma_{11}$ & $\gamma_{22}$  \\
\cmidrule(lr){3-8}
& $\theta_0$ & -0.4 & 0.3 & 0.3 & 0.2 & 1.0 & 0.5 \\
\cmidrule(lr){3-8}
& Min & -1.9678 & -0.43683 & 0.00045 & 0.00039 & -0.06681 & -0.22733 \\ 
& $Q_1$ & -0.54833 & 0.14011 & 0.15916 & 0.05114 & 0.722264 & 0.24915  \\
& Mean & -0.27079 & 0.20452 & 0.27219 & 0.13215 & 0.85198 & 0.29883 \\ 
& Rmse & 0.42777 & 0.15546 & 0.16409 & 0.14106 & 0.24889 & 0.22494 \\ 
\multicolumn{1}{c}{$n = 300$} & Bias & 0.12920 & -0.09547 & -0.02780 &  -0.06785 & -0.14801 & -0.20116 \\
& $Q_2$ & -0.25576 & 0.20369 & 0.25490 & 0.09102 & 0.87387 & 0.31681 \\ 
& Std & 0.40800 & 0.12275 & 0.16180 & 0.12373 & 0.20019 & 0.10070 \\ 
& $Q_3$ & -0.03088 & 0.27088 & 0.35315 & 0.16551 & 0.98917 & 0.36757 \\ 
& Max & 1.76732 & 0.68227 & 0.89650 & 0.91458 & 2.05482 & 0.61907 \\ 
\cmidrule(lr){3-8}
& Min & -2.45929 & -0.29875 & 0.00417 & 0.00185 & -0.12604 & -0.33428 \\
& $Q_1$ & -0.52847 & 0.15590 & 0.15208 & 0.05051 & 0.75660 & 0.25282  \\
& Mean & -0.30056 & 0.20639 & 0.25876 & 0.12658 & 0.87164 & 0.30351 \\ 
& Rmse & 0.40494 & 0.13894 & 0.15527 & 0.14347 & 0.24629 & 0.22152 \\ 
\multicolumn{1}{c}{$n = 2,000$} & Bias & 0.09943 & -0.09360 & -0.04123 & -0.07341 & -0.12835 & -0.19648 \\
& $Q_2$ & -0.26905 & 0.20141 & 0.24084 & 0.09055 & 0.88985 & 0.32210 \\ 
& Std & 0.39274 & 0.10273 & 0.14977 & 0.12332 & 0.21031 & 0.10236 \\ 
& $Q_3$ & -0.06193 & 0.25227 & 0.34025 & 0.15130 & 1.00659 & 0.37440 \\ 
& Max & 1.14336 & 0.80772 & 0.84607 & 0.89439 & 1.65371 & 0.60730 \\
\bottomrule
\end{tabular}
\caption{Summary statistics of parameters for model (\ref{eq:1}) with noise $\eta_t  = u_t\left(|u_{t-1}|+1\right)^{-1}$ , based on 1000 replications for sequence sizes of 300 and 2000, respectively.}
\label{tab:weak_2}
\end{table}

\begin{table}[H]
\centering
\setlength{\tabcolsep}{10pt}
\begin{tabular}{lccccccc}
\toprule
&  & $\alpha_{11}$ & $\alpha_{22}$ & $\beta_{11}$ & $\beta_{21}$ & $\gamma_{11}$ & $\gamma_{22}$  \\
\cmidrule(lr){3-8}
& $\theta_0$ & -0.4 & 0.3 & 0.3 & 0.2 & 1.0 & 0.5 \\
\cmidrule(lr){3-8}
& Min & -1.51651 & -0.52734 & 0.02252 & 0.01694 & -0.67843 & -0.18852 \\ 
& $Q_1$ & -0.62838 & 0.26199 & 0.22981 & 0.19064 & 0.95115 & 0.47577  \\
& Mean & -0.39454 & 0.36688 & 0.33641 & 0.27990 & 1.03722 & 0.57338 \\ 
& Rmse & 0.37315 & 0.22583 & 0.15232 & 0.14947 & 0.17451 & 0.17459 \\ 
\multicolumn{1}{c}{$n = 300$} & Bias & 0.00545 & 0.06688 & 0.03641 & 0.07990 & 0.03722 &  0.07338 \\
& $Q_2$ &  	-0.42387 & 0.38325 & 0.32349 & 0.26263 & 1.03171 & 0.56237 \\ 
& Std & 0.37330 & 0.21581 & 0.14798 & 0.12638 & 0.17058 & 0.15850 \\ 
& $Q_3$ & -0.22073 & 0.50270 & 0.42388 & 0.35306 & 1.12341 & 0.66141 \\ 
& Max & 1.31006 & 1.38173 & 0.94588 & 0.90364 & 1.94964 & 1.27796 \\ 
\cmidrule(lr){3-8}
& Min & -1.96784 & -0.43683 & 0.00045 & 0.0003 & -0.06681 & -0.22733\\
& $Q_1$ & -0.67414 & 0.25572 & 0.28388 & 0.27025 & 1.02601 & 0.58111  \\
& Mean & -0.33839 & 0.40620 & 0.38929 & 0.36814 & 1.14900 & 0.73178\\
& Rmse & 0.45897 & 0.29333 & 0.18361 & 0.22587 & 0.27650 & 0.32975 \\ 
\multicolumn{1}{c}{$n = 2,000$} & Bias & 0.06160 & 0.10620 & 0.08929 & 0.16814 & 0.14900 & 0.23178 \\
& $Q_2$ & -0.44907 & 0.46629 & 0.37886 & 0.35567 & 1.14965 & 0.73647\\ 
& Std & 0.45504 & 0.27357 & 0.16052 & 0.15089 & 0.23304 & 0.23467 \\ 
& $Q_3$  & -0.05628 & 0.58892 & 0.48749 & 0.44952 & 1.28062 & 0.88501\\ 
& Max & 1.34049 & 1.14302 & 0.91743 & 0.98627 & 2.03401 & 1.54148 \\
\bottomrule
\end{tabular}
\caption{Summary statistics of parameters for model (\ref{eq:1}) with noise $\eta_t = u_t^2u_{t-1}$, based on 1000 replications for sequence sizes of 300 and 2000, respectively.}
\label{tab:weak_3}
\end{table}

\begin{table}[H]
\centering
\setlength{\tabcolsep}{10pt}
\begin{tabular}{lccccccc}
\toprule
&  & $\alpha_{11}$ & $\alpha_{22}$ & $\beta_{11}$ & $\beta_{21}$ & $\gamma_{11}$ & $\gamma_{22}$  \\
\cmidrule(lr){3-8}
& $\theta_0$ & -0.4 & 0.3 & 0.3 & 0.2 & 1.0 & 0.5 \\
\cmidrule(lr){3-8}
& Min & -2.4364 & -0.49696 & 0.00052 & 0.00149 & -0.15889 & -0.38022 \\ 
& $Q_1$ & -0.55356 & 0.16637 & 0.17976 & 0.06301 & 0.83500 & 0.32146  \\
& Mean & -0.35396 & 0.24289 & 0.27857 & 0.13906 & 0.92532 & \\ 
& Rmse & 0.36926 & 0.14894 & 0.14436 & 0.12194 & 0.18385 & 0.17384 \\ 
\multicolumn{1}{c}{$n = 300$} & Bias & 0.04603 & -0.05710 & -0.0214 & -0.06093 &  -0.07467 & -0.13079\\
& $Q_2$ & -0.35941 & 0.24907 & 0.27072 & 0.11641 & 0.93218 & 0.38476 \\ 
& Std & 0.36656 & 0.13762 & 0.14283 & 0.10567 & 0.16809 &  0.11456\\ 
& $Q_3$ & -0.16511 & 0.32332 & 0.36559 & 0.18516 & 1.02073 & 0.44270\\ 
& Max & 1.27712 & 0.92749 & 0.78052 & 0.88233 & 2.17583 & 0.67173 \\ 
\cmidrule(lr){3-8}
& Min & -2.00079 & -0.16932 & 0.00187 & 0.00112 & 0.23911 & -0.12263 \\
& $Q_1$ & -0.51990 & 0.18937 & 0.19211 & 0.09294 & 0.87228 & 0.36941  \\
& Mean & -0.37027 & 0.24453 & 0.27217 & 0.14891 & 0.94797 & 0.41751 \\ 
& Rmse & 0.31198 & 0.10950 & 0.12014 & 0.09890 & 0.14060 & 0.12109 \\ 
\multicolumn{1}{c}{$n = 2,000$} & Bias & 0.02972 & -0.05546 & -0.02782 & -0.05108 & -0.05203 & -0.08248 \\
& $Q_2$ & -0.37526 & 0.24178 & 0.26365 & 0.13646 & 0.95252 & 0.42640 \\ 
& Std & 0.31071 & 0.09446 & 0.11694 & 0.08474 & 0.13068 & 0.08869 \\ 
& $Q_3$ & -0.21407 & 0.29864 & 0.33859 & 0.19526 & 1.03054 & 0.47740 \\ 
& Max & 1.77742 & 0.69406 & 0.88300 & 0.63743 & 1.51641 & 0.75889 \\
\bottomrule
\end{tabular}
\caption{Summary statistics of parameters for model (\ref{eq:1}) with noise $\eta_t = (\omega_0 + a_0 \eta_{t-1}^2 + \beta_0 h_{t-1})^{-1/2}$, based on 1000 replications for sequence sizes of 300 and 2000, respectively.}
\label{tab:weak_4}
\end{table}



%
%


\begin{figure}[H]
    \centering
    \includegraphics[width=1\textwidth]{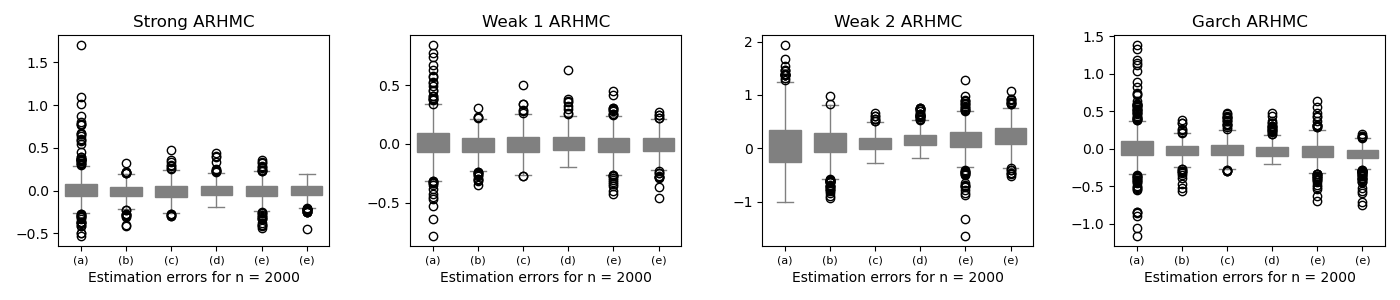}
    \includegraphics[width=1\textwidth]{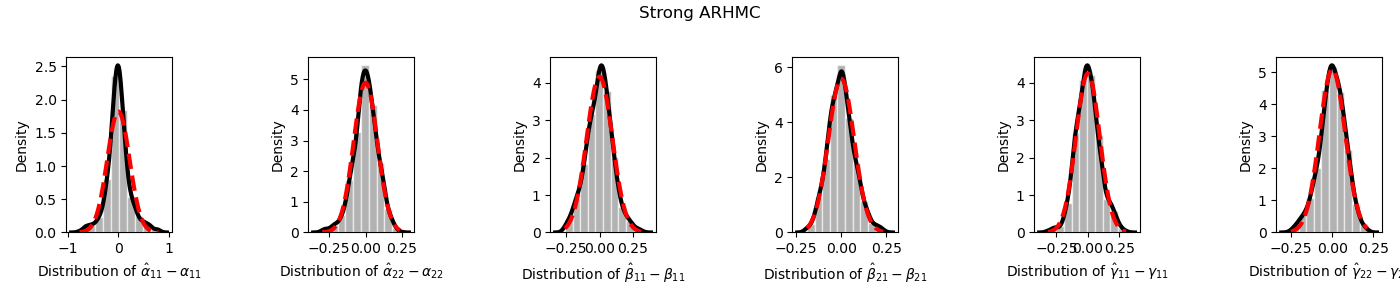}
    \includegraphics[width=1\textwidth]{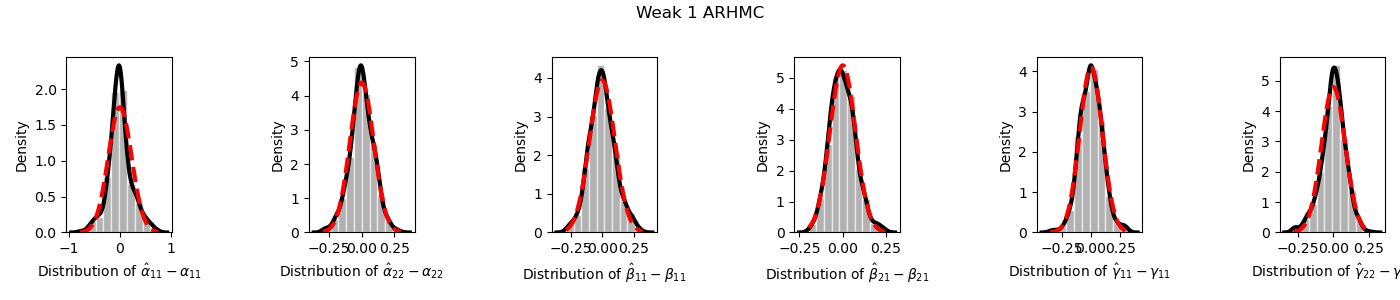}
    \includegraphics[width=1\textwidth]{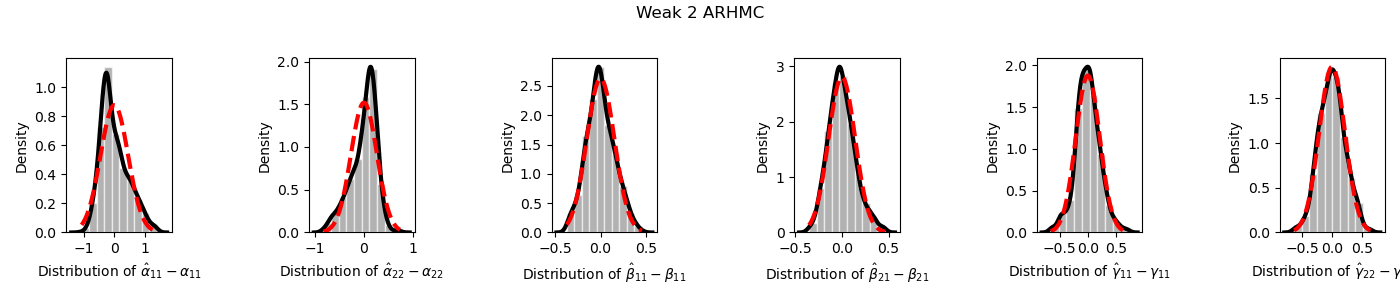}
    \includegraphics[width=1\textwidth]{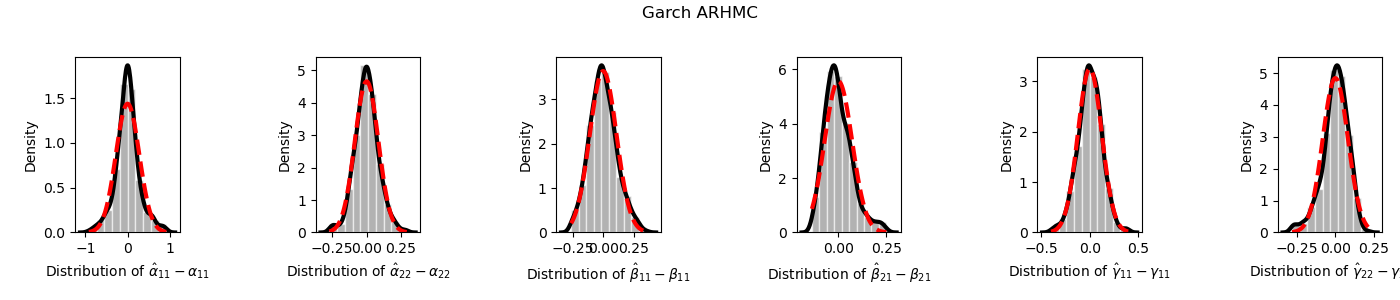}
    \caption{Boxplots and distribution of errors \(\hat{\theta}_n(i) - \theta_0(i)\) for \(i = 1, \ldots, 6\), where the noise \(\eta_t\) is defined as \(u_t\), \(u_t u_{t-1}\), \(u_t^2 u_{t-1}\), and \((\omega_0 + a_0 \eta_{t-1}^2 + \beta_0 h_{t-1} \, u_t)^{1/2}\), respectively. The kernel density estimate is displayed in full line and the centered Gaussian density with the same variance is plotted in dotted line.}
    \label{fig: density}
\end{figure}

\begin{figure}[H]
    \centering
     \includegraphics[width=1\textwidth]{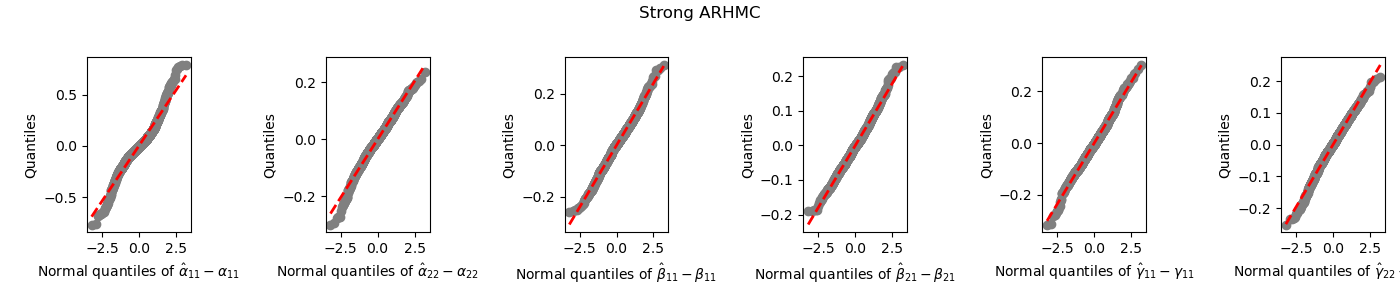}
     \includegraphics[width=1\textwidth]{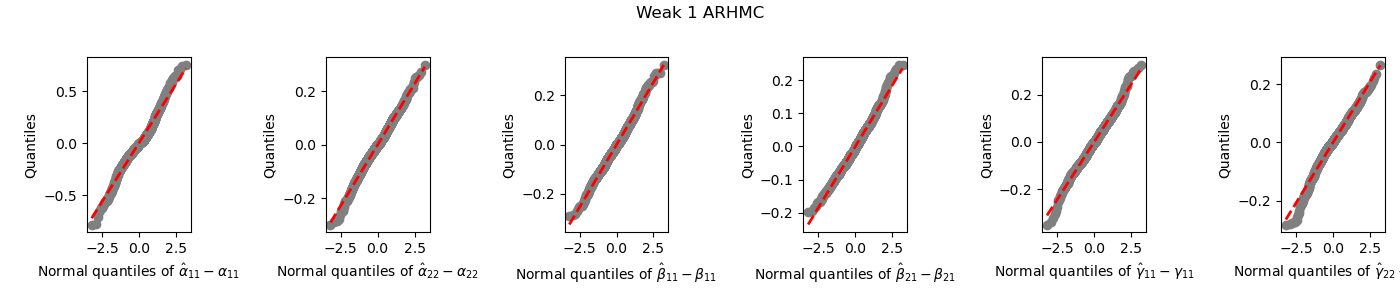}
      \includegraphics[width=1\textwidth]{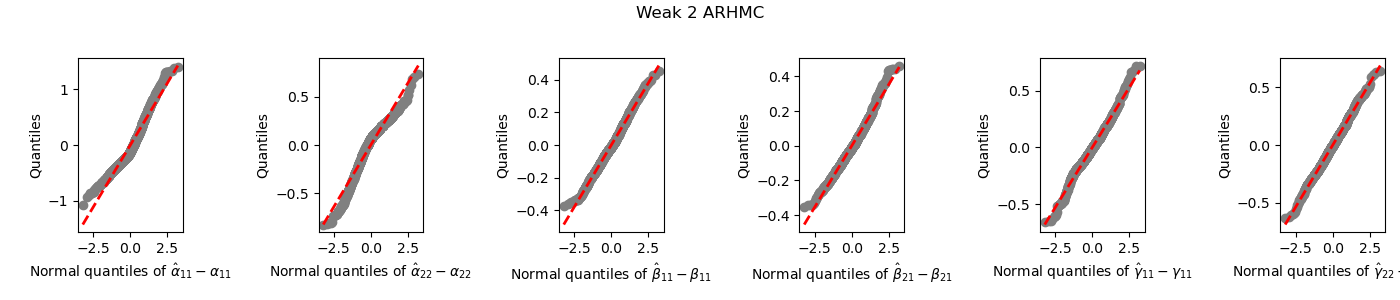}
    \includegraphics[width=1\textwidth]{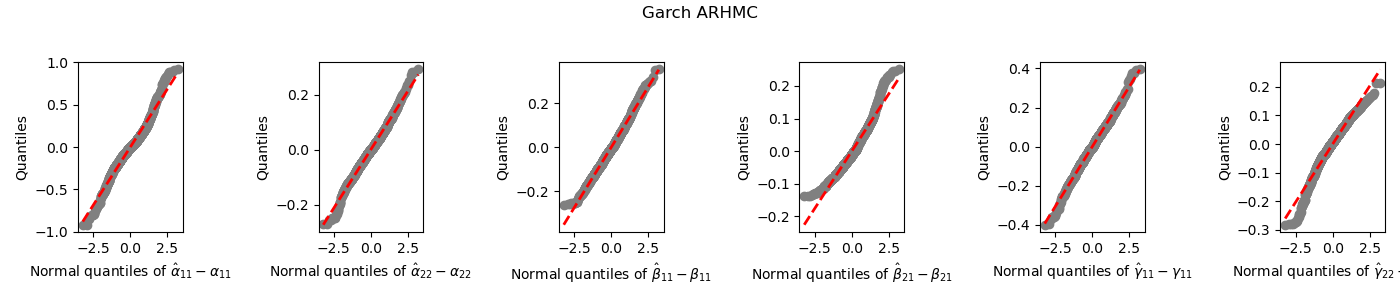}
    \caption{QQ-plots of errors \(\hat{\theta}_n(i) - \theta_0(i)\) for \(i = 1, \ldots, 6\), where the noise \(\eta_t\) is defined as \(u_t\), \(u_t u_{t-1}\), \(u_t^2 u_{t-1}\), and \((\omega_0 + a_0 \eta_{t-1}^2 + \beta_0 h_{t-1} \, u_t)^{1/2}\), respectively.}
    \label{fig: qqplot}
\end{figure}


\begin{figure}[H]
\begin{adjustwidth}{-1.8cm}{-1cm}
    \centering
     \includegraphics[width=1.2\textwidth]{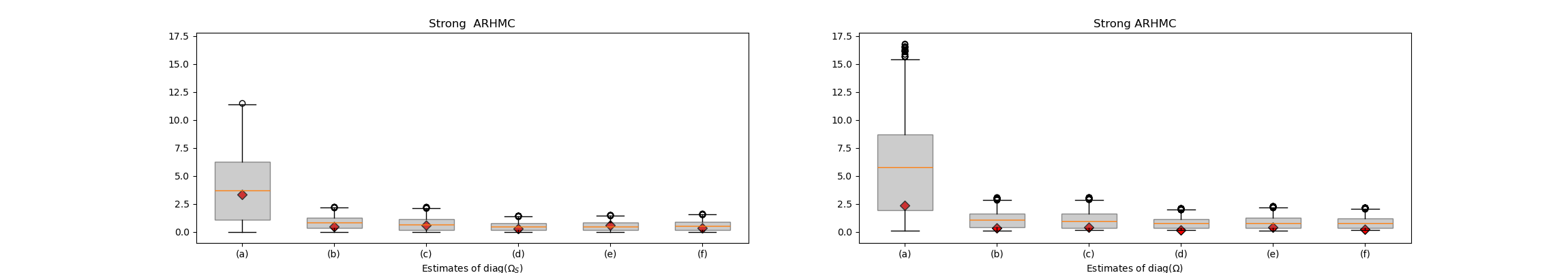}
      \includegraphics[width=1.2\textwidth]{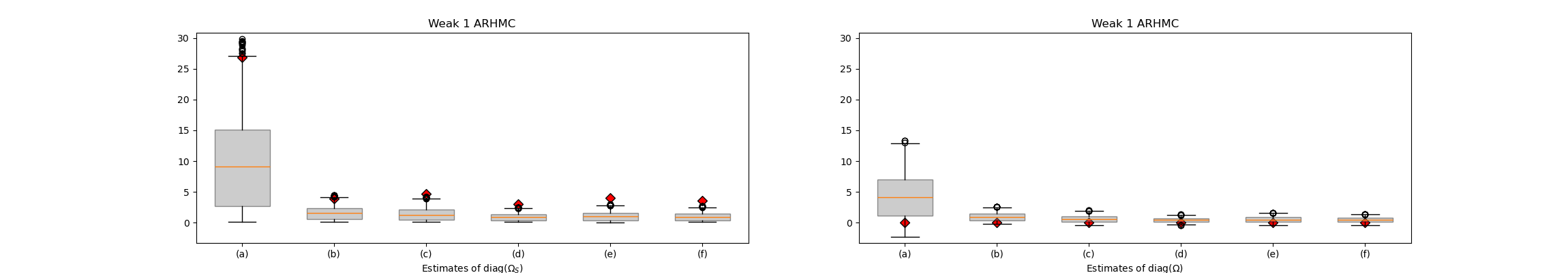}
      \includegraphics[width=1.2\textwidth]{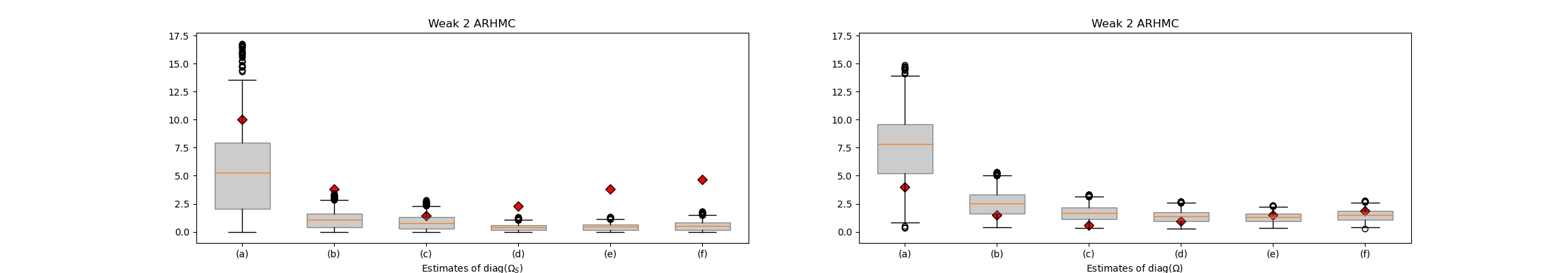}
    \end{adjustwidth}
\end{figure}
\begin{figure}[H]
\centering
    \caption{Comparison of standard and modified estimates of the asymptotic covariance matrix $\Omega$ of the moment estimator on the simulated model presented in Fig.
    \ref{fig: density} and Fig. \ref{fig: qqplot}
    \label{fig: omega}. Strong ARHMC corresponds to the model (\ref{eq:1}) with the noise \textit{Strong}, Weak 1 ARHMC corresponds to the model (\ref{eq:1}) with the noise \textit{Weak 1}, and Weak 2 ARHMC corresponds to the model (\ref{eq:1}) with the noise \textit{Weak 2}. The red diamond symbols represent the mean over $R=1,000$ replications of standardized squared errors: $n\{\hat{\alpha}_{11}+0.4\}^2$ for (a), $n\{\hat{\alpha}_{22}-0.3\}^2$ for (b), $n\{\hat{\beta}_{11}-0.3\}^2$ for (c), $n\{\hat{\beta}_{22}-0.3\}^2$ for (d), $n\{\hat{\gamma}_{11}-1.0\}^2$ for (e), and $n\{\hat{\gamma}_{22}-0.5\}^2$ for (f).}
\end{figure}



\section{Application to real data}
\begin{figure}[H]
 \centering
\includegraphics[width=1.0\textwidth]{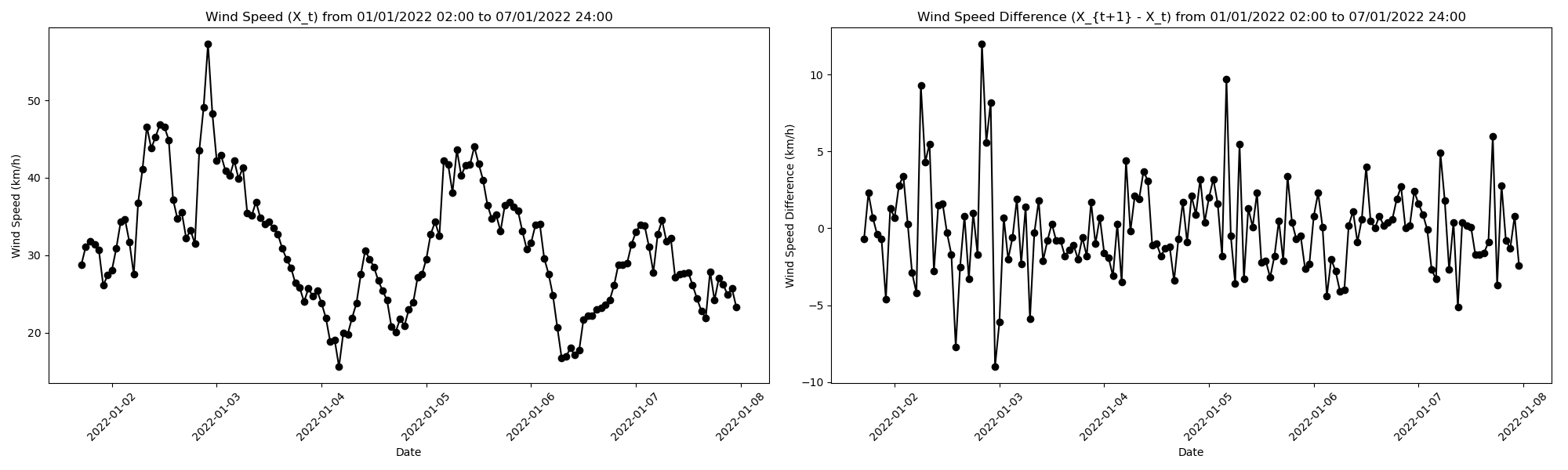}
\caption{    The figure presents two graphs representing wind data from January 1st to January 7, 2022. 
    The left graph shows the hourly wind speed at 100 meters above ground level \{\(H_t\)\}, 
    where higher values indicate stronger winds. 
    The right graph illustrates the difference in wind speed between consecutive hours \{\(H_{t+1} - H_t\)\}, 
    displaying a succession of tumultuous and calm periods. Large variations in this graph indicate periods of rapid wind speed changes (tumultuous), 
    while small variations represent stable, calm conditions.
}
\label{fig:wind-speed}
\end{figure}

In this section, we consider the hourly meteorological data from the Los Angeles region from January 1st to January 31, 2022, denoted by \( (H_t)_{t = 1, \dots, 744} \). The data were obtained from the website \textit{Open-Meteo.com}. In this dataset, we are specifically interested in the variable \textit{wind-speed-100m}, which represents the wind speed at 100 meters above the ground, measured in km/h. Fig. \ref{fig:wind-speed} plots the hourly wind speed at 100 meters and the differenced hourly wind speed between consecutive hours \{\(H_{t+1} - H_t\)\} from January 1 to January 7, 2022.
An initial analysis of the data shows a high variability in wind speed, with periods of strong winds towards the end of the month and calmer winds at the beginning. The maximum wind speed observed was 70.1 km/h, recorded on January 29, 2022, at 19:00, while the minimum wind speed of 3.8 km/h was observed on January 10, 2022, at 04:00. We consider \((X_t)_{t=1,\dots,743}\), the mean corrected series of the differenced series: $H_{t+1} - H_t$ to investigate the temporal variations in wind speed.

To illustrate the process of identifying the number of regimes \( K \) in model (\ref{eq:1}), we estimated the model parameters for \( K = 2 \), \( K = 3 \) and \( K = 4 \). The results of these estimations are presented in Tables \ref{tab:tab6}, \ref{tab:tab7} and \ref{tab:tab8}, respectively. The objective of these simulations is to determine the most appropriate number of regimes for accurately fitting the data.

 By analyzing the results of Tables \ref{tab:tab6}, \ref{tab:tab7} and \ref{tab:tab8}, we observe that  \( K=3 \) is the most relevant number of regimes to fit the data \( (X_t)_{t=1, \dots, 743} \). Indeed, for \( K=2 \), \( K=3 \) and \( K=4 \), we respectively obtain \((\hat{\beta}_{11}, \hat{\beta}_{22}) = (0.816, 0.321)\), \((\hat{\beta}_{11}, \hat{\beta}_{22}, \hat{\beta}_{33}) = (0.528, 0.427, 0.369)\), and \((\hat{\beta}_{11}, \hat{\beta}_{22}, \hat{\beta}_{33}, \hat{\beta}_{44}) = (0.503, 0.307, 0.324, 0.092)\). These results show that the fourth regime is less persistent than the first three, suggesting that a four-regime model includes an unnecessary regime. Noting also that the fact that, in certain regimes, the autoregressive parameters are in absolute value greater than $1$ (i.e explosive) such as the parameters $\hat{\alpha}_{22} = 1.534$ for $K=2$ and $\hat{\alpha}_{33} = 1.013$ for $K=3$ does not contradict the stationarity of the process $\{X_t\}_{t=1,\dots,744}$. Indeed, the stationarity condition given in Equation (\ref{eq:p1234}) is satisfied since, for $K=2$, $\hat{\pi}_1\log|\hat{\alpha}_{11}| + \hat{\pi}_2\log|\hat{\alpha}_{22}| = 0.787 \times \log|-0.342| + 0.212 \times \log|1.534| \approx -0.752 < 0$ and for $K=3$, $\hat{\pi}_1\log|\hat{\alpha}_{11}| + \hat{\pi}_2\log|\hat{\alpha}_{22}| + \hat{\pi}_3\log|\hat{\alpha}_{33}| = -0.275 \times \log|-0.893| + 0.325 \times \log|-0.287| + 0.400 \times \log|1.013| \approx -0.4312 < 0$. In other words, the presence of explosive regimes does not preclude the strict stationarity of the process unlike standard autoregressive processes. The process remains globally stationary as noted by \citet[page 343]{francq2001stationarity}.

To evaluate the significance of the autoregressive and those of the volatility parameters, their $p-$values and their \textit{standard errors} 
are presented in Table \ref{tab:tab10}. In view of Table \ref{tab:tab10}, it seems that only the autoregressive coefficient of the second regime $\hat{\alpha}_{22}$ is statistically insignificant at the 5\% significance level. The other coefficients, namely \(\hat{\alpha}_{11}\), \(\hat{\alpha}_{33}\), \(\hat{\gamma}_{11}\), \(\hat{\gamma}_{22}\), and \(\hat{\gamma}_{33}\) are all significant at  5\% level. 
Then, in a second step, the reduced ARHMC$(1)$ model was estimated with constraints on the autoregressive parameters with non-significant $p-$values in Table \ref{tab:tab10}, namely, the  coefficient of the regime 2 (\(\alpha_{22}\) is setting to be zero). The moments estimator of the final model are presented in Table \ref{tab:tab13}.

 As above-mentioned, the process $(X_t)_{t=1,\dots,743}$ is globally stationary, although we observe an explosive regime 1 with autoregressive parameter $\hat{\alpha}_{11}$. At the 5\% significance level, the autoregressive parameters $\hat{\alpha}_{11}$ and $\hat{\alpha}_{33}$, as well as the volatility parameters $\hat\gamma_{11}$, $\hat\gamma_{22}$ and $\hat\gamma_{33}$ are all significant  (see Table \ref{tab:tab13}). The fact that in regime 2, the  autoregressive parameter $\alpha_{22}$ is not significant can be explained by the nature of the process in this regime: $(X_t)_{t=1,\dots,743}$ behaves like multiplicative noise without any dependence on its past values. However, there is still some variability in this regime, as the absolute value of $\hat\gamma_{22}$ is 1.275. In contrast, the regimes 1 and 3 model calm and tumultuous periods of the process, respectively. In these regimes, the process is explained by its historical values, unlike in regime 2. It is also noted that the estimated coefficient $\hat\beta_{22}=0.188$ indicates that the process remains in regime 2 for a short period, which can be seen as a brief phase where the wind blows relatively strongly, before transitioning either to regime 1 with a probability $\hat\beta_{21}=0.641$ or to regime 3 with a probability $\hat\beta_{23}= 0.170$.

\begin{table}[H]
\centering
\setlength{\tabcolsep}{8pt}
\begin{tabular}{lccccccccc}
\toprule
\toprule
& \textbf{Parameter:} & $\alpha_{11}$ & $\alpha_{22}$ & $\beta_{11}$ & $\beta_{21}$ & $\beta_{12}$ & $\beta_{22}$ & $\gamma_{11}$ & $\gamma_{22}$  \\
\cmidrule(lr){3-10}
& \textbf{Estimate:} & -0.342 & 1.534 & 0.816 & 0.678 & 0.183 & 0.321 & 1.505 & 0.151 \\
\bottomrule
\bottomrule
\end{tabular}
\caption{Estimation of the parameters of model (\ref{eq:1}) based on  $(X_t)_{t=1,\dots, 743}$ with $K=2$.}
\label{tab:tab6}
\end{table}


\begin{table}[H]
\centering
\setlength{\tabcolsep}{7pt} 
\renewcommand{\arraystretch}{1.2} 
\begin{tabular}{lccccccccccccc}
\toprule
\toprule
& \textbf{Parameter:} & $\alpha_{11}$ & $\alpha_{22}$ & $\alpha_{33}$ & $\beta_{11}$ & $\beta_{21}$ & $\beta_{31}$ & $\beta_{12}$ & $\beta_{22}$ \\
\cmidrule(lr){3-10}
& \textbf{Estimate:} & -0.893 & -0.287 & 1.013 & 0.528 & 0.111 & 0.181 & 0.140 & 0.427 \\
\cmidrule(lr){3-10}
& \textbf{Parameter:} & $\beta_{32}$ & $\beta_{13}$ & $\beta_{23}$ & $\beta_{33}$ & $\gamma_{11}$ & $\gamma_{22}$ & $\gamma_{33}$ \\
\cmidrule(lr){3-9}
& \textbf{Estimate:} & 0.449 & 0.331 & 0.460 & 0.369 & 0.569 & -1.911 & -1.636 \\
\bottomrule
\bottomrule
\end{tabular}
\caption{Estimation of the parameters of model (\ref{eq:1}) based on $(X_t)_{t=1,\dots, 743}$ with $K=3$.}
\label{tab:tab7}
\end{table}

\begin{table}[H]
\centering
\setlength{\tabcolsep}{3pt} 
\renewcommand{\arraystretch}{1.2} 
\begin{tabular}{lccccccccccccc}
\toprule
\toprule
& \textbf{Parameter:} & $\alpha_{11}$ & $\alpha_{22}$ & $\alpha_{33}$ & $\alpha_{44}$ & $\beta_{11}$ & $\beta_{21}$ & $\beta_{31}$ & $\beta_{41}$ & $\beta_{12}$ & $\beta_{22}$ &$\beta_{32}$ & $\beta_{42}$ \\
\cmidrule(lr){3-14}
& \textbf{Estimate:} & 0.775 & -0.115 & 0.639 & 0.768 & 0.503 & 0.232 & 0.326 & 0.180 & 0.090 & 0.307 & 0.307 & 0.570\\
\cmidrule(lr){3-14}
& \textbf{Parameter:} & $\beta_{13}$ & $\beta_{23}$ & $\beta_{33}$ & $\beta_{43}$ & $\beta_{14}$ & $\beta_{24}$&$\beta_{34}$ & $\beta_{44}$& $\gamma_{11}$ & $\gamma_{22}$ & $\gamma_{33}$ & $\gamma_{44}$ \\
\cmidrule(lr){3-14}
& \textbf{Estimate:} & 0.133 & 0.336 & 0.324 & 0.157 & 0.180 & 0.570 &  0.157 & 0.092&-0.150 & 0.630 & 0.820 & -1.464  \\
\bottomrule
\bottomrule
\end{tabular}
\caption{Estimation of the parameters of model (\ref{eq:1}) based on  $(X_t)_{t=1,\dots, 743}$ with $K=4$.}
\label{tab:tab8}
\end{table}

\begin{table}[H]
\centering
\setlength{\tabcolsep}{10pt} 
\renewcommand{\arraystretch}{1.2} 
\begin{tabular}{lcccc}
\toprule
\toprule
& \textbf{Regimes} & $K=2$ & $K=3$ & $K=4$ \\
\midrule
\textbf{Estimated stationary distributions} & $\hat{\pi}_1$ & 0.787 & 0.275 & 0.277\\
& $\hat{\pi}_2$ & 0.212 & 0.325 & 0.306 \\
& $\hat{\pi}_3$ & --     & 0.400 & 0.295 \\
& $\hat{\pi}_4$ & --     & --    & 0.122 \\
\bottomrule
\bottomrule
\end{tabular}
\caption{Stationary distributions for regimes \( K = 2 \), \( K = 3 \) and \( K = 4 \).}
\label{tab:stationary_distributions}
\end{table}

\begin{table}[H]
\centering
\setlength{\tabcolsep}{8pt}
\begin{tabular}{lcccccccccc}
\toprule
\toprule
& \textbf{Estimate:} &$\hat{\alpha}_{11}$& $\hat{\alpha}_{22}$ & $\hat{\alpha}_{33}$ &  $\hat{\gamma}_{11}$& $\hat{\gamma}_{22}$ & $\hat{\gamma}_{33}$\\
\cmidrule(lr){3-9}
& \textbf{Standard error:} & 0.234 & 0.331 & 0.361 & 0.263 &  0.255 &  0.276 &   \\
\cmidrule(lr){3-9}
& \textbf{ p-value:} & \textbf{0.0001}  & 0.3851 & \textbf{0.0051} & \textbf{0.0307} & \textbf{0.0000} & \textbf{0.0000} &    \\
\bottomrule
\bottomrule
\end{tabular}
\caption{Summary of p-values and standard errors for the estimated autoregressive parameters \(\hat{\alpha}_{11}\), \(\hat{\alpha}_{22}\), \(\hat{\alpha}_{33}\) and the volatility parameters \(\hat{\gamma}_{11}\), \(\hat{\gamma}_{22}\), \(\hat{\gamma}_{33}\) given in Table \ref{tab:tab7}. The $p-$values less than 5\% are in bold.}
\label{tab:tab10}
\end{table}


\begin{table}[H]
\centering
\setlength{\tabcolsep}{7pt} 
\renewcommand{\arraystretch}{1.2} 
\begin{tabular}{lcccccccccccc}
\toprule
\toprule
& \textbf{Parameter:} & $\alpha_{11}$ & $\alpha_{33}$ & $\beta_{11}$ & $\beta_{21}$ & $\beta_{31}$ & $\beta_{12}$ & $\beta_{22}$ \\
\cmidrule(lr){3-9}
& \textbf{Estimate:} & -1.337 & 0.800 & 0.474 & 0.641 & 0.066 & 0.438 & 0.188 \\
\cmidrule(lr){3-9}
& \textbf{Standard error:}  & 0.2341 & 0.361 & -- & -- & -- & -- & -- \\
\cmidrule(lr){3-9}
& \textbf{p-value:} & 0.000 & 0.026 &  -- & -- & -- & -- & -- \\
\toprule
\toprule
& \textbf{Parameter:} & $\beta_{32}$ & $\beta_{13}$ & $\beta_{23}$ & $\beta_{33}$ & $\gamma_{11}$ & $\gamma_{22}$ & $\gamma_{33}$ \\
\cmidrule(lr){3-9}
& \textbf{Estimate:} & 0.263 & 0.086 & 0.170 & 0.670 & 0.637 & -1.275 & -0.988 \\
\cmidrule(lr){3-9}
& \textbf{Standard error:} & -- & -- & -- & -- & 0.263 & 0.255 & 0.276 \\
\cmidrule(lr){3-9}
& \textbf{p-value:} & -- & -- &  -- & -- & 0.015 & 0.004 & 0.000 \\
\bottomrule
\bottomrule
\end{tabular}
\caption{Re-estimation of the parameters of model (\ref{eq:1}) based on $(X_t)_{t=1,\dots, 743}$ with $K=3$.}
\label{tab:tab13}
\end{table}

\section{Conclusion}
In this paper, we studied a first-order autoregressive model where the parameters depend on a hidden Markov chain under the assumption that the noise is uncorrelated but not necessarily independent. 
First, we estimate the model parameters using the method of moments and establish its asymptotic properties which presented a significant challenge, given the limited literature on applying this method to such models, in contrast to the more extensive work on maximum likelihood or least squares approaches. We demonstrated the consistency and asymptotic normality of the moments estimator, followed by the estimation of the asymptotic variance-covariance matrix under certain mixing conditions applied to both the noise and the hidden Markov chain. 

In comparison to the work of \cite{boubacar2020estimation}, where the Markov chain is observed, we obtain similar sandwich expression for the asymptotic variance-covariance matrix. Furthermore, we conducted numerical simulations that perfectly illustrated our theoretical results, particularly the asymptotic normality established in Theorem \ref{thm:3}. We also applied our findings to real-world data following an approach similar to that of \cite{francq1997white} to determine the optimal number of regimes, represented by the parameter \(K\), in order to efficiently fit the data. Looking forward, our next challenge will be to validate this type of model by developing a portmanteau statistical test.

\section{Proofs}\label{proofs}

\subsection{Proof of Theorem \ref{thm:1}}\label{proof_thm1}
In order to proof Theorem \ref{thm:1}, one needs  the following lemma.
\begin{l1}[\cite{francq2004estimation}]\label{lem:1}
Let $(\Delta_t)_{t\in \mathbb{Z}}$ be an irreducible, aperiodic, and stationary Markov chain with state space in $\{1,\cdots,K\}$, transition probability $(p_{ij})_{i,j = 1,\ldots,K}$, and stationary distribution $\pi = \{\pi(1),\cdots,\pi(K)\}.$ Then, for any $k\geq 1$ and functions $f_1,\cdots,f_{k+1}$ defined on $\{1,\cdots,K\}$, we have
\begin{equation}
\label{eq:3}
    \mathbb{E}(f_1(\Delta_{t-1})\cdots f_k(\Delta_{t-k}) f_{k+1}(\Delta_{t-k-1})) = \bm{1}'\left(\prod_{\ell=1}^{k}Q_{f_\ell}\right)\bm{\pi}_{f_{k+1}}
\end{equation}
where for $1\leq \ell\leq k,$ $f_\ell(i)p_{ji}$ is the $(i,j)-$th element of a square matrix $Q_{f_\ell}$ of order $K$ and $\bm{\pi}_{f_{k+1}}: = \{f_{k+1}(1)\pi(1),\dots,f_{k+1}(K)\pi(K) \}'$.
\qed
\end{l1}



 
 Using Assumption $(\mathbf{A_3})$, Equation (\ref{eq:4}) and the fact that $\mathbb{E}(\eta_t) = 0$ for all $t$, we have for all $k\geq 0,$
\begin{align*}
c_{k,0}: = & \ \mathbb{E}(X_kX_0) = \mathbb{E}\left(\left[\sum_{n=0}^{\infty}\left(\prod_{j=0}^{n-1}a(\Delta_{k-j})\right)f(\Delta_{k-n})\eta_{k-n}\right]\left[\sum_{m=0}^{\infty}\left(\prod_{j=0}^{m-1}a(\Delta_{-j})\right)f(\Delta_{-m})\eta_{-m}\right]\right)\\
= & \sum_{n=0}^\infty\left(\sum_{i = 0}^n\mathbb{E}\left(\prod_{j=0}^{i-1}a(\Delta_{k-j})f(\Delta_{k-i})\eta_{k-i}\prod_{j = 0}^{n-i-1}a(\Delta_{-j})f(\Delta_{-n+i})\eta_{-n+i}\right)\right) \\
= & \sum_{n = 0}^{\infty}\left(\sum_{i=0}^{n}\mathds{1}_{\{i = \frac{k+n}{2}\}}\mathbb{E}\left(\prod_{j=0}^{i-1}a(\Delta_{k-j})f(\Delta_{k-i})\prod_{j = 0}^{n-i-1}a(\Delta_{-j})f(\Delta_{-n+i})\right)\right)\\
= & \sum_{n=0}^{\infty}\left(\mathbb{E}\left(\prod_{j=0}^{\frac{k+n}{2}-1} a(\Delta_{k-j})f(\Delta_{k-\frac{k+n}{2}})\prod_{j=0}^{n-\frac{k+n}{2}-1}a(\Delta_{-j})f(\Delta_{-n+\frac{k+n}{2}})\right)\mathds{1}_{\{k+n\in 2\mathbb{N}\}}\right).
\end{align*}

\noindent To obtain an explicit expression of \(c_{k,0}\), we need to distinguish between the case when \(k\) is even and the case when \(k\) is odd.
\medskip

\subparagraph{$\diamond$ Case 1 : $k$ is even.}\ \\

For even $k$ $(k=2k_1,k_1\in\mathbb{N})$, we have
\begin{align*}
c_{2k_1,0} = & \ \sum_{n_1 = k_1}^{\infty}\left(\mathbb{E}\left(\prod_{j=0}
^{n_1+k_1-1}a(\Delta_{2k_1-j})f(\Delta_{k_1-n_1})\prod_{j=0}^{n_1-k_1-1}a(\Delta_{-j})f(\Delta_{-n_1+k_1})\right)\right).
\end{align*}
Making the change of variable $j: = j-2k_1$, we then obtain
\begin{align*}
c_{2k_1,0} = & \ \sum_{n_1=k_1}^{\infty}\left(\mathbb{E}\left(\prod_{j=-2k_1}^{n_1-k_1-1}a(\Delta_{-j})f(\Delta_{k_1-n_1})\prod_{j=0}^{n_1-k_1-1}a(\Delta_{-j})f(\Delta_{-n_1+k_1})\right)\right)\\
 = & \sum_{n_1=k_1}^{\infty}\left(\mathbb{E}\left(\prod_{j=-2k_1}^{-1}a(\Delta_{-j})\prod_{j=0}^{n_1-k_1-1}a^2(\Delta_{-j})f^2(\Delta_{-n_1+k_1})\right)\right).
\end{align*}
In view of Lemma~\ref{lem:1} we get for even $k$ that:
\begin{align*}
c_{2k_1,0} = & \ \sum_{n_1=k_1}^{\infty}\bm{1}'\left(\prod_{l=1}^{2k_1}Q_a\prod_{l=1}^{n_1-k_1}Q_{a^2}\right)\pi_{f^2}=\bm{1}'\sum_{n_1=k_1}^{\infty}\left(Q_a^{2k_1}Q_{a^2}^{n_1-k_1}\right)\bm\pi_{f^2}.
\end{align*}

\subparagraph{$\diamond$ Case 2 : $k$ is odd.}\ \\

For odd $k$ $(k = 2k_1+1,k_1\in\mathbb{N})$, we have
\begin{align*}
c_{2k_1+1,0} = & \ \sum_{n_1 = k_1}^{\infty}\left(\mathbb{E}\left(\prod_{j=0}
 ^{n_1+k_1}a(\Delta_{2k_1+1-j})f(\Delta_{k_1-n_1})\prod_{j=0}^{n_1-k_1-1}a(\Delta_{-j})f(\Delta_{-n_1+k_1})\right)\right)\\
  = & \ \sum_{n_1=k_1}^{\infty}\left(\mathbb{E}\left(\prod_{j=-2k_1-1}^{n_1-k_1-1}a(\Delta_{-j})f(\Delta_{k_1-n_1})\prod_{j=0}^{n_1-k_1-1}a(\Delta_{-j})f(\Delta_{-n_1+k_1})\right)\right)\\
  = & \ \sum_{n_1=k_1}^{\infty}\left(\mathbb{E}\left(\prod_{j=-2k_1-1}^{-1}a(\Delta_{-j})\prod_{j=0}^{n_1-k_1-1}a^2(\Delta_{-j})f^2(\Delta_{-n_1+k_1})\right)\right).
\end{align*}
Using Lemma~\ref{lem:1}, we obtain for odd $k$ that:
\begin{align*}
c_{2k_1+1,0} =  & \ \sum_{n_1 = k_1}^{\infty}\left(\bm{1}'\left(\prod_{l=1}^{2k_1+1}Q_a\prod_{l=1}^{n_1-k_1}Q_{a^2}\right)\bm\pi_{f^2}\right) =\bm{1}'\sum_{n_1=k_1}^{\infty}\left(Q_a^{2k_1+1}Q_{a^2}^{n_1-k_1}\right)\bm\pi_{f^2}.
\end{align*}
For any $k\in\mathbb{N}$ we draw the conclusion  that:
\begin{equation}
\label{eq:9}
c_{k,0} =  \ \bm{1}'Q_a^{k}(I_K-Q_{a^2})^{-1}\bm\pi_{f^2}.
\end{equation}
The proof of Theorem \ref{thm:1} is then complete.
\qed 

\subsection{Proof of Theorem \ref{thm:2}}\label{proof_thm2}
The proof of this theorem is a direct consequence of the implicit function theorem.

 Consider the following differentiable function
\begin{equation*}
\Phi :
    \left\lbrace
        \begin{aligned}
            & \mathring\Theta \times \mathbb{R}^{N} & \rightarrow & \ \quad \mathbb{R}^{K^2+K} \\
            & \quad (\theta, z) & \mapsto & \quad J_{\Psi^{N}}(\theta)'\left(z-\Psi^N(\theta)\right).
        \end{aligned}
    \right.
\end{equation*}
Differentiation  with respect to $\theta$ yields 
\begin{align*}
   \frac{\partial}{\partial\theta}\Phi(\theta,z)  = \frac{\partial}{\partial\theta}J_{\Psi^N}(\theta)'\left(z-\Psi^N(\theta)\right) - J_{\Psi^{N}}(\theta)'J_{\Psi^{N}}(\theta),
\end{align*}
which implies that
\begin{align*}
   \frac{\partial}{\partial\theta}\Phi(\theta_0,\mathbf{c}_{0}^N)  = - J_{\Psi^{N}}(\theta_0)'J_{\Psi^{N}}(\theta_0),
\end{align*}
where
 $\mathbf{c}_{0}^N: = \left(\mathbb{E}(X_0X_k)\right)_{k=1\dots N}$. Since by hypothesis: $\bm{r}(\theta_0) = K^2+K$, it follows that the matrix $J'_{\Psi^{N}}(\theta_0)J_{\Psi^{N}}(\theta_0)$ is invertible. Thus by the implicit function theorem, there exists a neighborhood $\mathcal{V}_{\theta_0}$ of $\theta_0$ in $\Theta$, a neighborhood $\mathcal{V}_{\mathbf{c}_0^N}$ of $\mathbf{c}_{0}^N$ in $\mathbb{R}^{N}$ and a continuous function $\vartheta$ from $\mathcal{V}_{\mathbf{c}_0^N}$ to $\mathcal{V}_{\theta_0}$ such that
\begin{align}
\label{eq:0001}
  J_{\Psi^N}(\theta(z))'\left(z-\Psi^N(\theta)\right) = 0 \Leftrightarrow \theta = \vartheta(z) \quad \forall (\theta,z)\in\mathcal{V}_{\theta_0}\times\mathcal{V}_{\mathbf{c}_0^N} .
\end{align}
Next, we set $\hat{\theta}_n := \vartheta(\hat{c}_0^{N})$ where $\hat{c}_0^{N} : = (\hat{c}_{k,0})_{k=1,\dots,N}$. By the ergodic theorem, $\hat{c}_0^{N}$ converges a.s. to $c_0^{N}$, as $n\to\infty$. Thus by the continuity of the function $\vartheta(\cdot)$, we obtain
\begin{align*}
    \hat{\theta}_n = \vartheta(\hat{c}_0^{N}) \xrightarrow[n\rightarrow\infty]{a.s.} \vartheta(c_0^N) = \theta_0,
\end{align*}
which proves the consistency of the estimator $\hat{\theta}_n$. Furthermore, the uniqueness of the solutions to (\ref{eq:0001}) in $\mathcal{V}_{\theta_0} \times \mathcal{V}_{c_0^N}$ implies the uniqueness of the existence of $\hat{\theta}_n$.

\qed
\subsection{Proof of Theorem \ref{thm:3}}\label{proof_thm3}
In this subsection, we shall give the proof of the asymptotic distribution of \(\hat{\theta}_n\) based on the following series of Lemmas. Lemma~\ref{lem:5} provides Davydov's inequality, a crucial result for analyzing strongly mixing processes. Lemma~\ref{lem:6} establishes the conditions under which the process \((X_t)_{t \in \mathbb{Z}}\) defined in (\ref{eq:1}) is not only strictly stationary but also admits moments of sufficiently high order necessary for the proof of asymptotic normality. Lemmas~\ref{lem:111}, ~\ref{lem:7} and ~\ref{lem:8} respectively confirm the finiteness of the moments of the process \((X_t)_{t \in \mathbb{Z}}\), the existence of the asymptotic variance matrix and the asymptotic distribution of the random vector \(\sqrt{n}F^{n}(\theta_0)\).

Let us suppose that the conditions of Theorem~\ref{thm:3} are satisfied. Since the functions $\psi_k$, $1\leq k\leq N$ are smooth functions over all $\theta$ in $\Theta$, it follows that 
\begin{align*}
    \displaystyle\sup_{n}\displaystyle\sup_{\theta\in\Theta}\left|\frac{\partial^2 \mathcal{F}^{N,n}(\theta)}{\partial\theta(i)\partial\theta(j)}\right|<\infty,\quad \text{ for }\quad 1\leq i,j\leq K^2+K.
\end{align*}
In view of Theorem~\ref{thm:2}, we have almost surely that $\hat{\theta}_n$ converges to $\theta_0\in\mathring\Theta$. By employing a standard Taylor expansion around $\theta_0$ and noting that $\mathcal{F}^{N,n}(\hat{\theta}_n) = 0$, we obtain 
\begin{align}
\label{eq:00012}
  0 = \sqrt{n}\mathcal{F}^{N,n}(\theta_0) + \sqrt{n}\nabla\mathcal{F}^{N,n}(\theta_n^\star)\left(\hat{\theta}_n-\theta_0\right),
\end{align}
where the parameter ${\theta_n^\star}$ lies on the segment in $\mathbb{R}^{K^2+K}$ with endpoint $\hat{\theta}_n$ and $\theta_0$. Proceeding with another Taylor expansion, we also have
\begin{align}
\label{eq:299}
    \|\nabla\mathcal{F}^{N,n}(\theta_n^\star) - \nabla\mathcal{F}^{N,n}(\theta_0)\| &= \|\nabla^2\mathcal{F}^{N,n}(\theta)\left(\theta_n^\star - \theta_0\right)\| \notag\\
    &\leq \sup_{n}\sup_{\theta \in \Theta}\|\nabla^2\mathcal{F}^{N,n}(\theta)\|\left\|\theta_n^\star - \theta_0\right\| \xrightarrow[n\rightarrow\infty]{a.s} 0.
\end{align}
Using the ergodic theorem once again, we easily observe that
\begin{align}
\label{eq:normality32}
    \nabla\mathcal{F}^{N,n}(\theta_0) &= \left.\frac{\partial}{\partial \theta}J_{\Psi^N}(\theta)\right|_{\theta=\theta_0} F^{N,n}(\theta_0) + J_{\Psi^N}(\theta_0)'J_{\Psi^N}(\theta_0) \notag\\
    &\xrightarrow[n\rightarrow\infty]{a.s} J_{\Psi^N}(\theta_0)'J_{\Psi^N}(\theta_0) := M^N(\theta_0).
\end{align}
Along with Equations (\ref{eq:00012}), (\ref{eq:299}) and (\ref{eq:normality32}), we obtain
\begin{align}
\label{eq:2999}
\notag
 0 & = \sqrt{n}\mathcal{F}^{N,n}(\theta_0)+(\nabla\mathcal{F}^{N,n}(\theta_n^\star)-\nabla\mathcal{F}^{N,n}(\theta_0)+\nabla\mathcal{F}^{N,n}(\theta_0))\sqrt{n}(\hat{\theta}_n - \theta_0) \notag\\&  =  \sqrt{n}\mathcal{F}^{N,n}(\theta_0)+\nabla\mathcal{F}^{N,n}(\theta_0)\sqrt{n}(\hat{\theta}_n -\theta_0)+\mathrm{o}_\mathbb{P}(1) \notag\\ & = \sqrt{n}\mathcal{F}^{N,n}(\theta_0)+(\nabla\mathcal{F}^{N,n}(\theta_0)-M^N(\theta_0)+M^N(\theta_0))\sqrt{n}(\hat{\theta}_n - \theta_0)+\mathrm{o}_\mathbb{P}(1) \notag\\ & = \sqrt{n}\mathcal{F}^{N,n}(\theta_0)+M^N(\theta_0)\sqrt{n}(\hat{\theta}_n - \theta_0)+\mathrm{o}_\mathbb{P}(1).
\end{align}

\noindent Consequently, under the hypothesis $\mathbf{r}(\theta_0) = K^2 + K$, so that the matrix $M^N(\theta_0)$ is invertible, it follows from Equation (\ref{eq:2999}) that:
\begin{align*}
    \sqrt{n}\left(\hat{\theta}_n - \theta_0\right) = \left(M^N(\theta_0)\right)^{-1} \sqrt{n}\mathcal{F}^{N,n}(\theta_0)
\end{align*}
The proof of Theorem \ref{thm:3} then directly follows from Lemma \ref{lem:7} and Lemma \ref{lem:8} below and by using Slutsky's Theorem, we obtain that $\sqrt{n}\left(\hat{\theta}_n - \theta_0\right)$ has a limiting normal distribution with mean 0 and covariance matrix $M^{-1}J'I JM^{-1}$.
\qed

\begin{l1}[\text{\cite{davydov1968convergence}}]\label{lem:5}
Let \(X\) and \(Y\) be two random variables, and let \(\sigma(X)\) and \(\sigma(Y)\) be the \(\sigma\)-fields generated by \(X\) and \(Y\) respectively. Consider three strictly positive numbers \(p\), \(q\), and \(r\) such that \(p^{-1} + q^{-1} + r^{-1} = 1\). Then,
\[
|\Cov(X,Y)| \leq C_0 \|X\|_p \|Y\|_q \left\{\alpha\left(\sigma(X), \sigma(Y)\right)\right\}^\frac{1}{r},
\]
where \(\|.\|_p\) denotes the \(L^p\)-norm, \(C_0\) is a universal constant and 
$\alpha\left(\sigma(X), \sigma(Y)\right)$ denotes the strong mixing coefficient between the $\sigma$-fields $\sigma(X)$ and $\sigma(Y)$ generated by $X$ and $Y$.
\end{l1}
\begin{l1}\label{lem:6}
Under  Assumptions $(\mathbf{A_4})$ and $(\mathbf{A_6})$, for \(\nu>0\) we have 
\begin{align}
    \mathbb{E}\left(\left[\left(\prod_{j = 0}^{i-1}a(\Delta_{t-j})\right)f(\Delta_{t-i})\right]^{8+4\nu}\right) \leq C\rho^i,
\end{align}
where \(C\) is a positive universal constant and \(\rho\) is a constant in \((0,1)\).
\end{l1}
\subsubsection*{Proof.}
In view of Lemma \ref{lem:1} we obtain 

\[
\mathbb{E}\left(\left(\prod_{j = 0}^{i-1}a^{8+4\nu}(\Delta_{t-j})\right)f^{8+4\nu}(\Delta_{t-i})\right) = \bm{1}'\left(\prod_{l = 1}^{i}Q_{a^{8+4\nu}}\right)\bm\pi_{f^{8+4\nu}},
\]
where $Q_{a^{8+4\nu}} = A_{\theta_0}^{8+4\nu}P_{\theta_0}'$ and the column vector $\bm\pi_{f^{8+4\nu}} = \left(f^{8+4\nu}(1)\pi(1),\ldots,f^{8+4\nu}(K)\pi(K)\right)'$. 

 Let \(\epsilon > 0\). Then on \(\mathbb{R}^{K}\) there exists an induced norm \(\|\cdot\|_{\epsilon}\) such that \(\|A_{\theta_0}^{8+4\nu}P_{\theta_0}'\|_{\epsilon}\leq \varrho + \epsilon\), where \(\varrho\) is the spectral radius of the matrix \(A_{\theta_0}^{8+4\nu}P_{\theta_0}'\). Using this induced norm, we have 
\[
\left|\bm{1}'\left(\prod_{l = 1}^{i}Q_{a^{8+4\nu}}\right)\bm\pi_{f^{8+4\nu}} \right| \leq \|\bm{1}'\|_{\epsilon}\|A_{\theta_0}^{8+4\nu}P_{\theta_0}'\|_{\epsilon}^{i}\|\pi_{f^{8+4\nu}}\|_{\epsilon} \leq (\varrho+\epsilon)^{i}\|\bm\pi_{f^{8+4\nu}}\|_{\epsilon}.
\]
 Finally, under Assumptions \((\mathbf{A_2})\) and \((\mathbf{A_4})\), taking \(\epsilon\) small enough, we can find a positive constant \(C\) and \(\rho \in (0,1)\) such that 
\[
\mathbb{E}\left(\prod_{j = 0}^{i-1}a^{8+4\nu}(\Delta_{t-j})f^{8+4\nu}(\Delta_{t-i})\right) \leq C\rho^{i}.
\]

 The conclusion follows from these arguments.
\qed

\begin{not1}
\label{not:1}
 In the remainder of this paper, to simplify the notation, we define for all \((t, i) \in \mathbb{Z} \times \mathbb{N}\).
\begin{align}
\label{eq:not}
    c(\theta_0, \Delta_t, \ldots, \Delta_{t-i}) := \left(\prod_{j = 0}^{i-1} a(\Delta_{t-j})\right) f(\Delta_{t-i}).
\end{align}
\end{not1}

\begin{l1}
\label{lem:111}
Under the assumptions  $(\mathbf{A_3})$ and  $(\mathbf{A_4})$, for $\nu>0$ we have 
\begin{align}
 \left\|X_0\right\|_{4+2\nu} < \infty.   
\end{align}
\end{l1}

\subsubsection*{Proof.}
Using the independence between the chain $(\Delta_t)_{t\in\mathbb{Z}}$ and the noise $(\eta_t)_{t\in\mathbb{Z}}$, there exists a positive constant $C$ and $\rho \in (0,1)$ such that
\begin{align*}
    \left\|X_0\right\|_{4+2\nu} = \left\|\sum_{i=0}^{\infty}c(\theta_0,\Delta_0, \ldots, \Delta_{-i})\eta_i\right\|_{4+2\nu} \leq \sum_{i=0}^{\infty}\left\|c(\theta_0,\Delta_0, \ldots, \Delta_{-i})\right\|_{4+2\nu}\left\|\eta_i\right\|_{4+2\nu}\leq \left(\sum_{i=0}^{\infty} C\rho^i\right)\left\|\eta_0\right\|_{4+2\nu}
\end{align*}

\noindent Hence, the proof is concluded based on Assumption $(\mathbf{A_4})$.
\qed

\begin{l1}\label{lem:7}
Let Assumptions of the Theorem~\ref{thm:3} be satisfied. Then the following convergence holds
$$I^N(\theta_0) : = \displaystyle\lim_{n\rightarrow \infty}\Var\Big(\sqrt{n}F^{N,n}(\theta_0)\Big)\in\mathbb{R}^{N\times N}.$$
\end{l1}
\subsubsection*{Proof.}
By the stationarity of the process $(X_t)_{t\in\mathbb{Z}}$, we have
\[
I_n^N(\theta_0) := \text{Var}\left(\sqrt{n}F^{N,n}(\theta_0)\right) = n \ \text{Cov}\left(\frac{1}{n}\sum_{t = 1}^{n}Y_t, \frac{1}{n}\sum_{s = 1}^{n}Y_s\right) = \frac{1}{n}\sum_{k = -n+1}^{n-1}\left(n-|k|\right)\text{Cov}(Y_t,Y_{t-k}).
\]
Next, we introduce the notation 
\[
c_k(l,r) := \text{Cov}(Y_t(l),Y_{t-k}(r)) = \Cov(X_tX_{t+l},X_{t-k}X_{t-k+r}),
\]
where \(Y_t(l)\) denotes the \(l\)-th element of \(Y_t\) defined in (\ref{eq:001}) with $l,r = 1,\dots,N$.

 In view of Equation \eqref{eq:4}, the coefficients \((c_k(l,r))_{k \in \mathbb{Z}}\) can be rewritten as:
\begin{align*}
    c_k(l,r) = \displaystyle\sum_{i_1 = 0}^{\infty}\displaystyle\sum_{i_2 = 0}^{\infty}\displaystyle\sum_{i_3 = 0}^{\infty}\displaystyle\sum_{i_4=0}^\infty\zeta_{i_1,i_2,i_3,i_4,k}(l,r)
\end{align*}
where the terms \(\zeta_{i_1,i_2,i_3,i_4,k}(l,r)\) are defined as lengthy covariances and expectations involving the coefficients \(c(\theta_0,\Delta_t, \ldots, \Delta_{t-i})\) and the noise terms \(\eta_{t-i}\). More formally by using $(\mathbf{A_3})$, we obtain
\begin{align*}
    \zeta_{i_1,i_2,i_3,i_4,k}(l,r) = &\ \text{Cov}\Big(c(\theta_0,\Delta_t,\cdots,\Delta_{t-i_1})\eta_{t-i_1}) c(\theta_0,\Delta_{t+l},\cdots,\Delta_{t+l-i_2})\eta_{t+l-i_2}),\\& \hspace{1cm} c(\theta_0,\Delta_{t-k},\cdots,\Delta_{t-k-i_3})\eta_{t-k-i_3}c(\theta_0,\Delta_{t-k+r},\cdots,\Delta_{t-k+r-i_4})\eta_{t-k+r-i_4}\Big) \\ = &
    \ \mathbb{E}\Big[c(\theta_0,\Delta_t,\cdots,\Delta_{t-i_1})c(\theta_0,\Delta_{t+l},\cdots,\Delta_{t+l-i_2})\\& \hspace{0.5cm} c(\theta_0,\Delta_{t-k},\cdots,\Delta_{t-k-i_3})c(\theta_0,\Delta_{t-k+r},\cdots,\Delta_{t-k+r-i_4})\Big] \\ & \times \mathbb{E}\left(\eta_{t-i_1}\eta_{t+l-i_2}\eta_{t-k-i_3}\eta_{t-k+r-i_4}\right)-\mathbb{E}\Big[c(\theta_0,\Delta_{t},\cdots,\Delta_{t-i_1})c(\theta_0,\Delta_{t+l},\cdots,\Delta_{t+l-i_2})\Big]\\ &  \times \mathbb{E}\Big[c(\theta_0,\Delta_{t-k},\cdots,\Delta_{t-k-i_3})c(\theta_0,\Delta_{t-k+r},\cdots,\Delta_{t-k+r-i_4})\Big]\\ & \times\mathbb{E}\left(\eta_{t-k-i_3}\eta_{t-k+r-i_4}\right) \times\mathbb{E}\left(\eta_{t-i_1}\eta_{t+l-i_2}\right)\\ 
 = &\mathbb{E}\Big[c(\theta_0,\Delta_t,\cdots,\Delta_{t-i_1})c(\theta_0,\Delta_{t+l},\cdots,\Delta_{t+l-i_2})\\& \hspace{0.5cm} c(\theta_0,\Delta_{t-k},\cdots,\Delta_{t-k-i_3})c(\theta_0,\Delta_{t-k+r},\cdots,\Delta_{t-k+r-i_4})\Big]\\ &  \times \text{Cov}(\eta_{t-i_1}\eta_{t+l-i_2},\eta_{t-k-i_3}\eta_{t-k+r-i_4}) \\ & +  \text{Cov}(c(\theta_0,\Delta_t,\cdots,\Delta_{t-i_1}) c(\theta_0,\Delta_{t+l},\cdots,\Delta_{t+l-i_2}),\\& \hspace{1.2cm} c(\theta_0,\Delta_{t-k},\cdots,\Delta_{t-k-i_3})c(\theta_0,\Delta_{t-k+r},\cdots,\Delta_{t-k+r-i_4}))\\ & \times \mathbb{E}\left(\eta_{t-i_1}\eta_{t+l-i_2}\right)\mathbb{E}\left(\eta_{t-k-i_3}\eta_{t-k+r-i_4}\right).
\end{align*}
Then applying Lemma~\ref{lem:6} and using the Cauchy-Schwarz inequality, we have
\begin{align*}
    &\Big|\mathbb{E}(c(\theta_0,\Delta_t,\cdots,\Delta_{t-i_1})c(\theta_0,\Delta_{t+l},\cdots,\Delta_{t+l-i_2})\\ & \hspace{0.5cm} c(\theta_0,\Delta_{t-k},\cdots,\Delta_{t-k-i_3})c(\theta_0,\Delta_{t-k+r},\cdots,\Delta_{t-k+r-i_4}))\Big| \notag \\
    & \leq \left(\mathbb{E}\left[c(\theta_0,\Delta_t,\cdots,\Delta_{t-i_1})^4\right]\mathbb{E}\left[c(\theta_0,\Delta_{t+l},\cdots,\Delta_{t+l-i_2})^4\right]\right)^{\frac{1}{4}} \\ 
    & \times \left(\mathbb{E}\left[c(\theta_0,\Delta_{t-k},\cdots,\Delta_{t-k-i_3})^4\right]\mathbb{E}\left[c(\theta_0,\Delta_{t-k+r},\cdots,\Delta_{t-k+r-i_4})^4\right]\right)^{\frac{1}{4}} \notag \\
    & \leq C\rho^{i_1+i_2+i_3+i_4}
\end{align*}
where $C$ is a positive constant and $\rho\in (0,1)$.
\subparagraph*{$\diamond$ Assuming \(k \geq 0\).}\ \\

 Applying the dominated convergence theorem, we can control \(|c_k(l,r)|\) as follows
\begin{align*}
    \left|c_k(l,r)\right| = & \left|\text{Cov}(Y_t(l),Y_{t-k}(r)\right|\\ = &\left|\displaystyle\sum_{i_1 = 0}^{\infty}\displaystyle\sum_{i_2 = 0}^{\infty}\displaystyle\sum_{i_3 = 0}^{\infty}\displaystyle\sum_{i_4 = 0}^{\infty}\zeta_{i_1,i_2,i_3,i_4,k}(l,r)\right|\\ \leq & \displaystyle\sum_{i_1 = 0}^{\infty}\displaystyle\sum_{i_2 = 0}^{\infty}\displaystyle\sum_{i_3 = 0}^{\infty}\displaystyle\sum_{i_4 = 0}^{\infty}\left|\zeta_{i_1,i_2,i_3,i_4,k}(l,r)\right| \\ \leq & \  u_1 + u_2 + u _3 + u_4 + u_5 + v_1 + v_2 + v_3
\end{align*}
where 
\begin{align*}
   u_1 & := \sum_{i_1 > \frac{k}{2}} \sum_{i_2 = 0}^{\infty} \sum_{i_3 = 0}^{\infty} \sum_{i_4 = 0}^{\infty} C\rho^{i_1+i_2+i_3+i_4} \left|\text{Cov}(\eta_{t-i_1}\eta_{t+l-i_2}, \eta_{t-k-i_3}\eta_{t-k+r-i_4})\right|, \\
   u_2 & := \sum_{i_1 =0}^{\infty} \sum_{i_2>\frac{k}{2}} \sum_{i_3 = 0}^{\infty} \sum_{i_4 = 0}^{\infty} C\rho^{i_1+i_2+i_3+i_4} \left|\text{Cov}(\eta_{t-i_1}\eta_{t+l-i_2}, \eta_{t-k-i_3}\eta_{t-k+r-i_4})\right|, \\
   u_3 & := \sum_{i_1 = 0}^{\infty} \sum_{i_2 = 0}^{\infty} \sum_{i_3>\frac{k}{2}} \sum_{i_4 = 0}^{\infty} C\rho^{i_1+i_2+i_3+i_4} \left|\text{Cov}(\eta_{t-i_1}\eta_{t+l-i_2}, \eta_{t-k-i_3}\eta_{t-k+r-i_4})\right|,\\
   u_4 & := \sum_{i_1 = 0}^{\infty} \sum_{i_2 = 0}^{\infty} \sum_{i_3 = 0}^{\infty} \sum_{i_4>\frac{k}{2}} C\rho^{i_1+i_2+i_3+i_4} \left|\text{Cov}(\eta_{t-i_1}\eta_{t+l-i_2}, \eta_{t-k-i_3}\eta_{t-k+r-i_4})\right|, \\
    u_5 & := \sum_{i_1\leq\frac{k}{2}} \sum_{i_2\leq\frac{k}{2}} \sum_{i_3 \leq\frac{k}{2}} \sum_{i_4\leq\frac{k}{2}} C\rho^{i_1+i_2+i_3+i_4} \left|\text{Cov}(\eta_{t-i_1}\eta_{t+l-i_2}, \eta_{t-k-i_3}\eta_{t-k+r-i_4})\right|, 
\\
    v_1 & := \sum_{i_1 > \frac{k}{2}} \sum_{i_3 = 0}^{\infty} \left|\text{Cov}\left(c(\theta_0,\Delta_t,\cdots,\Delta_{t-i_1}) c(\theta_0,\Delta_{t-l},\cdots,\Delta_{t-i_1}), \right.\right.\\
    & \hspace{5.5em} \left.\left. c(\theta_0,\Delta_{t-k},\cdots,\Delta_{t-k-i_3})c(\theta_0,\Delta_{t-k+r},\cdots,\Delta_{t-k-i_3})\right)\right|, \\
    v_2 & := \sum_{i_1 = 0}^{\infty} \sum_{i_3 > \frac{k}{2}} \left|\text{Cov}\left(c(\theta_0,\Delta_t,\cdots,\Delta_{t-i_1}) c(\theta_0,\Delta_{t-l},\cdots,\Delta_{t-i_1}), \right.\right.\\
    & \hspace{5.5em} \left.\left. c(\theta_0,\Delta_{t-k},\cdots,\Delta_{t-k-i_3})c(\theta_0,\Delta_{t-k+r},\cdots,\Delta_{t-k-i_3})\right)\right|, \\
    v_3 & := \sum_{i_1 \leq \frac{k}{2}} \sum_{i_3 \leq \frac{k}{2}} \left|\text{Cov}\left(c(\theta_0,\Delta_t,\cdots,\Delta_{t-i_1}) c(\theta_0,\Delta_{t-l},\cdots,\Delta_{t-i_1}), \right.\right.\\
    & \hspace{5.5em} \left.\left. c(\theta_0,\Delta_{t-k},\cdots,\Delta_{t-k-i_3})c(\theta_0,\Delta_{t-k+r},\cdots,\Delta_{t-k-i_3})\right)\right|.
\end{align*}
Moreover, employing the Cauchy-Schwarz inequality again and using Assumption (\(\mathbf{A_4}\)), we obtain
\begin{align*}
\left|\text{Cov}(\eta_{t-i_1}\eta_{t+l-i_2}, \eta_{t-k-i_3}\eta_{t-k+r-i_4})\right| & \leq 
\mathbb{E}|\eta_{t-i_1}\eta_{t+l-i_2}\eta_{t-k-i_3} \eta_{t-k+r-i_4}| + \mathbb{E}|\eta_{t-i_1}\eta_{t+l-i_2}|\mathbb{E}|\eta_{t-k-i_3}\eta_{t-k+r-i_4}|\\ & \leq 
\left(\mathbb{E}|\eta_{t-i_1}|^4\right)^{\frac{1}{4}}\left(\mathbb{E}|\eta_{t+l-i_2}|^4\right)^{\frac{1}{4}}\left(\mathbb{E}|\eta_{t-k-i_3}|^4\right)^{\frac{1}{4}}\left(\mathbb{E}|\eta_{t-k+r-i_4}|^4\right)^{\frac{1}{4}} \\ & + (\mathbb{E}|{\eta_{t-i_1}}|^2)^\frac{1}{2}(\mathbb{E}|{\eta_{t+l-i_2}}|^2)^\frac{1}{2}+\mathbb{E}|{\eta_{t-k-i_3}}|^2)^\frac{1}{2}(\mathbb{E}|{\eta_{t-k+r-i_4}}|^2)^\frac{1}{2}\\ & \leq  \mathbb{E}|\eta_t|^4 + (\mathbb{E}|\eta_{t}|^2)^2 <  \infty.
\end{align*}
Hence, there exists a set of positive constants \((C_i)_{1\leq i \leq 4}\) satisfying
\begin{align*}
    u_i\leq C_i\rho^{\frac{k}{2}},\quad \forall \ 1\leq i \leq 4.
\end{align*}
Under Assumption (\(\mathbf{A_4}\)) and by the Cauchy-Schwarz inequality, we have \(\mathbb{E}|\eta_t\eta_{t^{\prime}}|^{2+\nu} < \infty\) for some \(\nu>0\). Given that $(2+\nu)^{-1} + (2+\nu)^{-1} + \nu(2+\nu)^{-1} = 1,\quad$ \(0\leq i_1,i_2,i_3,i_4\leq k/2\) and \(1\leq l,r\leq N\), 
by applying Lemma~\ref{lem:5}, we can find two positive constants \(C^{\prime}_5,C_5\) such that 
\begin{align*}
   u_5 & := \sum_{i_1\leq\frac{k}{2}} \sum_{i_2\leq\frac{k}{2}} \sum_{i_3 \leq\frac{k}{2}} \sum_{i_4\leq\frac{k}{2}} C\rho^{i_1+i_2+i_3+i_4} \left|\text{Cov}(\eta_{t-i_1}\eta_{t+l-i_2}, \eta_{t-k-i_3}\eta_{t-k+r-i_4})\right| \\ & \leq \sum_{i_1\leq\frac{k}{2}} \sum_{i_2\leq\frac{k}{2}} \sum_{i_3 \leq\frac{k}{2}} \sum_{i_4\leq\frac{k}{2}}C^{\prime}_5 \rho^{i_1+i_2+i_3+i_4}\left(\mathbb{E}|\eta_{t-i_1}\eta_{t+l-i_2}|^{2+\nu}\right)^{\frac{1}{2+\nu}} \left(\mathbb{E}|\eta_{t-k-i_3}\eta_{t-k-l+r-i_4}|^{2+\nu}\right)^{\frac{1}{2+\nu}}\\ & \hspace{1.8cm} \times \Big(\alpha_{\eta}\left(\min\left\{k+i_4-i_1-r,k+i_4-i_2-r+l,k+i_3-i_1,k+i_3-i_2+l\right\}\right)\Big)^{\frac{\nu}{2+\nu}}\\ & \leq C_5\Big(\alpha_{\eta}\Big(\lfloor k/2\rfloor-r\Big)\Big)^{\frac{\nu}{2+\nu}}.
\end{align*}
Thanks to Cauchy-Schwarz inequality, there exists two positive constants $C^{\prime}_1$ and $C^{\prime\prime}_1$ such that 
\begin{align}
\label{eq:14}
v_1 & := \sum_{i_1 > \frac{k}{2}} \sum_{i_3 = 0}^{\infty} \left|\text{Cov}\left(c(\theta_0,\Delta_t,\cdots,\Delta_{t-i_1}) c(\theta_0,\Delta_{t-l},\cdots,\Delta_{t-i_1}), \right.\right. \notag\\
    & \hspace{5.5em} \left.\left. c(\theta_0,\Delta_{t-k},\cdots,\Delta_{t-k-i_3})c(\theta_0,\Delta_{t-k+r},\cdots,\Delta_{t-k-i_3})\right)\right| \notag \\ 
    & \leq \sum_{i_1 > \frac{k}{2}} \sum_{i_3 = 0}^{\infty} \mathbb{E}|c(\theta_0,\Delta_t,\cdots,\Delta_{t-i_1})c(\theta_0,\Delta_{t+l},\cdots,\Delta_{t-k-i_1})\notag\\ & \hspace{1.6cm} c(\theta_0,\Delta_{t-k},\cdots,\Delta_{t-i_3})c(\theta_0,\Delta_{t-k+r},\cdots,\Delta_{t-k-i_3})| \notag \\ &  \hspace{1.5cm} + \mathbb{E}|c(\theta_0,\Delta_t,\cdots,\Delta_{t-i_1})c(\theta_0,\Delta_{t+l},\cdots,\Delta_{t-k-i_1})|\notag\\ &    \hspace{1.5cm} \times \mathbb{E}|c(\theta_0,\Delta_{t-k},\cdots,\Delta_{t-i_3})c(\theta_0,\Delta_{t-k+r},\cdots,\Delta_{t-k-i_3})| \notag \\
    &  \leq \sum_{i_1 > \frac{k}{2}} \sum_{i_3 = 0}^{\infty} C^{\prime\prime}_1\rho^{i_1+i_3} \notag \\ & \leq C^{\prime}_1 \rho^{\frac{k}{2}}.
\end{align}
Similarly, following the reasoning for \( v_1 \) there exists a positive constant \( C^{\prime}_2 \) such that 
\begin{align*}
    v_2 \leq C^{\prime}_2\rho^{\frac{k}{2}}.
\end{align*}
Notice that given the assumptions made on the Markov chain \( (\Delta_t)_{t\in\mathbb{Z}} \), it follows from \citet[Theorem 3.1 or Theorem 3.2]{bradley2005basic} that the process \( (\Delta_t)_{t\in\mathbb{Z}} \) is strongly mixing and satisfies 
$\sum_{h = 0}^{\infty} \alpha_{\Delta}(h)^{\frac{\nu}{2+\nu}} < \infty$, for a certain $\nu>0$ and 
where \( \alpha_\Delta \) is defined as in Equation (\ref{eq:12}).
Hence, using once again Lemma~\ref{lem:5} and the fact that $0\leq i_1,i_3\leq k/2$ and $1\leq l,r\leq N$, we also have
\begin{align*}
   v_3 & := \sum_{i_1 \leq \frac{k}{2}} \sum_{i_3 \leq \frac{k}{2}} \left|\text{Cov}\left(c(\theta_0,\Delta_t,\cdots,\Delta_{t-i_1}) c(\theta_0,\Delta_{t-l},\cdots,\Delta_{t-i_1}), \right.\right.\\
    & \hspace{5.5em} \left.\left. c(\theta_0,\Delta_{t-k},\cdots,\Delta_{t-k-i_3})c(\theta_0,\Delta_{t-k+r},\cdots,\Delta_{t-k-i_3})\right)\right|\\
    & \leq \Big(\mathbb{E}\left|c(\theta_0,\Delta_t,\cdots,\Delta_{t-i_1}) c(\theta_0,\Delta_{t-l},\cdots,\Delta_{t-i_1})\right|^{2+\nu}\Big)^{\frac{1}{2+\nu}}\\ & \times \Big(\mathbb{E}\left|c(\theta_0,\Delta_{t-k},\cdots,\Delta_{t-k-i_3})c(\theta_0,\Delta_{t-k+r},\cdots,\Delta_{t-k-i_3})\right|^{2+\nu}\Big)^{\frac{1}{2+\nu}}\\ & \times \Big(\alpha_{\Delta}\left(\min\{k+i_3-i_1,k+i_3-l-1\}\right)\Big)^{\frac{2}{2+\nu}}\\ & \leq C^{\prime}_3\alpha^{\frac{\nu}{2+\nu}}_{\Delta}\left(\min\Big\{k-l-1,\lfloor k/2\rfloor \Big\}\right),
\end{align*}
where $C^{\prime}_3$ is a strictly positive constant. 

In conclusion, considering all the preceding upper bounds, it follows that for $k\geq 0$ there exists three positive constants $M_1,M_2,M_3$ such that
\begin{align*}
    \displaystyle\sum_{k = 0}^{\infty}|c_k(l,r)|\leq M_1\displaystyle\sum_{k = 0}^{\infty}\rho^{\frac{k}{2}} + M_2\displaystyle\sum_{k = 0}^{\infty}\alpha_{\eta}^{\frac{\nu}{2+\nu}}\left(\lfloor k/2\rfloor-r\right) + M_3\displaystyle\sum_{k = 0}^{\infty}\alpha^{\frac{\nu}{2+\nu}}_{\Delta}\left(\min\Big\{k-l-1,\lfloor k/2\rfloor\Big\}\right)< \infty. 
\end{align*}
\subparagraph*{$\diamond$ The same bounds then clearly holds for \(k \leq 0\):}\ \\
 $$\displaystyle\sum_{k = -\infty}^{0}|c_k(l,r)| < \infty.$$
Therefore we have
\begin{align*}
   \displaystyle\sum_{k = -\infty}^{\infty}|c_k(l,r)| < \infty.
\end{align*}
Then by applying the dominated convergence theorem, we obtain
\begin{align*}
     I_n^N(\theta_0) = \frac{1}{n}\displaystyle\sum_{k = -n+1}^{n-1}\left(n-|k|\right)\text{Cov}(Y_t,Y_{t-k})\xrightarrow[n\rightarrow\infty]\displaystyle\sum_{k = -\infty}^{\infty}|c_k(l,r)|.
\end{align*}
The proof is complete. 
\qed

\begin{l1}\label{lem:8}
Under the assumptions of Theorem~\ref{thm:3}, the random vector $\sqrt{n}F^{N,n}(\theta_0)$ has a limiting normal distribution with mean $0$ and covariance matrix $I$.
\end{l1}
\subsubsection*{Proof.}
Using the definition of $F^{N,n}$ given in Equation (\ref{eq:800}), the fact that $\mathbb{E}_{\theta_0}(\hat{c}_{k,0}) = c_{k,0}$ for all $1 \leq k \leq N$ entails that $\mathbb{E}_{\theta_0}(\sqrt{n} F^{N,n}(\theta_0)) = 0$. In other words, the statistic $\sqrt{n}F^{N,n}(\theta_0)$ is centered. 
\\
Let $1\leq p\leq N.$ Applying the Cauchy's product formula to two convergent series, we have
\begin{align*}
    X_tX_{t+p} = & \left(\displaystyle\sum_{i_1 = 0}^{\infty}\left(\prod_{j_1 = 0}^{i_1-1} a(\Delta_{t-j_1})f(\Delta_{t-i_1})\right)\eta_{t-i_1}\right)\left(\displaystyle\sum_{i_2=0}^{\infty}\left(\displaystyle\prod_{j_2=0}^{i_2-1}a(\Delta_{t+p-j_2})f(\Delta_{t+p-i_2})\right)\eta_{t+p-i_2}\right)\\
     = & \displaystyle\sum_{i = 0}^{\infty}\displaystyle\sum_{k = 0}^{i}\left(\displaystyle\prod_{j = 0}^{k-1} a(\Delta_{t+p-j})f(\Delta_{t+p-k})\right)\left(\displaystyle\prod_{j=0}^{i-k-1}a(\Delta_{t-j})f(\Delta_{t+k-i})\right)\eta_{t+p-k}\eta_{t+k-i}\\
     = & \displaystyle\sum_{i=0}^{\infty}\displaystyle\sum_{k = 0}^{i}c(\theta_0,\Delta_{t+p},\cdots,\Delta_{t+p-k})c(\theta_0,\Delta_{t},\cdots,\Delta_{t-i+k})\eta_{t+p-k}\eta_{t+k-i}.
\end{align*}
Let $s$ be a positive integer. We introduce the following convenient notation 
 \begin{align}
\label{eq:113}
    d_{k}^{t}:= d_{k}^{t}(\theta_0) = & \ c(\theta_0,\Delta_{t},\cdots,\Delta_{t-k}), \ k\geq 0
\end{align} and 
\begin{align*}
    Y_{t,s} := & \displaystyle\sum_{i=0}^{s}\displaystyle\sum_{k=0}^{i}\left(d_{k}^{t+1}d_{i-k}^{t}\eta_{t+1-k}\eta_{t+k-i},\cdots,d_{k}^{t+N}d_{i-k}^{t}\eta_{t+N-k}\eta_{t+k-i}\right)^{\prime},\\
    Z_{t,s} := & \displaystyle\sum_{i=s+1}^{\infty}\displaystyle\sum_{k=0}^{i}\left(d_{k}^{t+1}(\theta_0)d_{i-k}^{t}\eta_{t+1-k}\eta_{t+k-i},\cdots,d_{k}^{t+N}(\theta_0)d_{i-k}^{t}(\theta_0)\eta_{t+N-k}\eta_{t+k-i}\right)^{\prime},
\end{align*}
so that $Y_{t,s} + Z_{t,s} = Y_t$ and where we recall that 
$c(\theta_0,\Delta_{t},\dots,\Delta_{t-k})$ is defined in Notation~\ref{not:1} (See Equation (\ref{eq:not})). We then have
\begin{align*}
    \sqrt{n}F^{N,n}({\theta_0}) = &\frac{1}{\sqrt{n}}\displaystyle\sum_{t=1}^{n}(Y_t(\theta_0)-\mathbb{E}_{\theta_0}(Y_t(\theta_0))
    =  \frac{1}{\sqrt{n}}\displaystyle\sum_{t=1}^{n}\left(Y_{t,s}-\mathbb{E}_{\theta_0}(Y_{t,s})\right) + \frac{1}{\sqrt{n}}\displaystyle\sum_{t=1}^{n}\left(Z_{t,s}-\mathbb{E}_{\theta_0}(Z_{t,s})\right).
\end{align*} 
The process \( (Y_{t,s})_{t\in\mathbb{Z}} \) depends on $\eta_k$ and $\Delta_k$ for $k$ in finite set. Furthermore, as the processes \( (\Delta_t)_{t\in\mathbb{Z}} \) and \( (\eta_t)_{t\in\mathbb{Z}} \) are strongly mixing, according to the assumption \( (\mathbf{A_4}) \), in view of \citet[Theorem 14.1 p. 210]{davidson1994stochastic} it follows that the process \( (Y_{t,s}) \) is strongly mixing.\ In addition, it can be deduced from \citet[Theorem 5.1]{bradley2005basic} that the mixing coefficients \( (\alpha_{Y,s}(h))_{h \in \mathbb{Z}} \) of the process \( (Y_{t,s} )\) satisfy \( \alpha_{Y,s}(h) \leq \alpha_{\Delta,\eta}(\max\{0,h-s\}) \leq \alpha_{\Delta}(\max\{0,h-s+1\}) + \alpha_{\eta}(\max\{0,h-s\})\). Applying the limit central theorem for the process strongly mixing  (see \citet[see Theorem 1.7, p. 367]{ibragimov1962some}), it follows that 
\begin{align*}
    \frac{1}{\sqrt{n}}\displaystyle\sum_{t=1}^{n}(Y_{t,s}-\mathbb{E}_{\theta_0}(Y_{t,s}))
\end{align*}
has a limiting normal $\mathcal{N}(0,I_s)$ distribution with $I_s : = \sum_{h=-\infty}^{\infty}\text{Cov}(Y_{t,s}Y_{t-h,s})\xrightarrow[s\rightarrow\infty]{} I.$ 

To complete the proof, it will suffice to show that $$\mathbb{E}\left(\left(\frac{1}{\sqrt{n}}\sum_{t=1}^{n}(Z_{t,s}-\mathbb{E}_{\theta_0}(Z_{t,s})\right)\left(\frac{1}{\sqrt{n}}\sum_{t=1}^{n}(Z_{t,s}-\mathbb{E}_{\theta_0}(Z_{t,s})\right)^{\prime}\right)$$ converges uniformly to zero as outlined in \citet[Lemma 3]{francq1998estimating} or in \citet[Lemma A.3]{boubacar2020estimation}  and we will conclude thanks to a result given by \citet[Corollary 7.7.1, p. 426]{anderson1971family}.

Due stationarity of $(X_t)_{t\in\mathbb{Z}}$ and using the fact that the process $(Z_{t,s}-{\mathbb{E}}_{\theta_0}(Z_{t,s}))_{t\in\mathbb{Z}}$ is centered, for all $1\leq p \leq N$ it follows that
\begin{align}
\label{eq:1800}
 \Var\left(\frac{1}{\sqrt{n}}\displaystyle\sum_{t=1}^{n}(Z_{t,s}(p)-\mathbb{E}_{\theta_0}(Z_{t,s}(p)))\right)  = & 
   \ \frac{1}{n}\displaystyle\sum_{t=1}^{n}\displaystyle\sum_{r=1}^{n}\Cov\left(Z_{t,s},Z_{r,s}\right) \notag\\  
 =  & \ \frac{1}{n}\displaystyle\sum_{h = -n+1}^{n-1}\left(n-|h|\right)\Cov\left(Z_{t,s}(p),Z_{t-h,s}(p)\right) \notag\\ 
 \leq & \ \displaystyle\sum_{h=-\infty}^{\infty}|c_s^Z(h)|
\end{align}
where  we define
\begin{align*}
   c_s^Z(h) = & \ \text{Cov}(Z_{t,s},Z_{t-h,s}). 
\end{align*}
We prove in what follows that the series in \eqref{eq:1800} tends to zero as $n\to\infty$ by bounding appropriately $c_s^Z(h)$.
\subparagraph{$\diamond$ Case 1 : suppose that $h\geq 0$ and $\lfloor h/2\rfloor>s$.}\ \\

 Then we can write 
\begin{align*}
    Z_{t,s}(p) : = Z^{h-}_{t,s}(p) + Z^{h+}_{t,s}(p), \ p = 1,\dots,N
\end{align*}
where 
\begin{align*}
 Z^{h-}_{t,s}(p) = \displaystyle\sum_{i_1=s+1}^{[h/2]}\displaystyle\sum_{j_1 = 0}^{i_1} d_{j_1}^{t+p}(\theta_0)d_{i_1-j_1}^{t}(\theta_0)\eta_{t+p-j_1}\eta_{t+j_1-i_1} 
\end{align*}
\text{and} 
\begin{align*}
 Z^{h+}_{t,s}(p) =\displaystyle\sum_{i_1=[h/2]+1}^{\infty}\displaystyle\sum_{j_1 = 0}^{i_1} d_{j_1}^{t+p}(\theta_0)d_{i_1-j_1}^{t}(\theta_0)\eta_{t+p-j_1}\eta_{t+j_1-i_1}.  
\end{align*}
Using the Cauchy-Schwarz inequality, the assumptions $(\mathbf{A_3})$, $(\mathbf{A_4})$, Lemma~\ref{lem:6} and thanks to the stationarity of $(\eta_t)_{t\in\mathbb{Z}}$, for $\nu>0$ we get 
\begin{align}
\label{eq:1010}
   \|d_{j}^{t+p}(\theta_0)d_{i-j}^{t}(\theta_0)\eta_{t+p-j}\eta_{t+j-i}\|_{2+\nu} \leq & \|d_{j}^{t-p}(\theta_0)\|_{4+2\nu} \|d_{i-j}^{t}(\theta_0)\|_{4+2\nu}\|\eta_{t+p-j}\|_{4+2\nu}\|\eta_{t+j-i}\|_{4+2\nu} \notag \\ \leq & \ C^{(0)}\rho^i \|\eta_t\|_{4+2\nu}^2
\end{align}
where $C^{(0)}$ is a positive constant. Applying once again the Cauchy-Schwarz inequality, using Lemma~\ref{lem:5} and Equation (\ref{eq:1010}), we have on the one hand
\begin{align}
\label{eq:15}
    |\text{Cov}({Z^{h^-}_{t,s}(p)},Z_{t-h,s}(p))| \leq & \displaystyle\sum_{i_1=s+1}^{[h/2]}\displaystyle\sum_{j_1=0}^{i_1}\displaystyle\sum_{i_2=s+1}^{\infty}\displaystyle\sum_{j_2=0}^{i_2}|\text{Cov}(d_{j_1}^{t+p}(\theta_0)d_{i_1-j_1}^{t}(\theta_0)\eta_{t+p-j_1}\eta_{t+j_1-i_1},
   \notag \\ & \hspace{9.5em} d_{j_2}^{t-h+p}(\theta_0)d_{i_2-j_2}^{t-h}({\theta_0})\eta_{t-h+p-j_2}\eta_{t-h+j_2-i_2})| \notag \\ \leq & \ C^{(1)}\displaystyle\sum_{i_1=0}^{[h/2]}\displaystyle\sum_{j_1=0}^{i_1}\displaystyle\sum_{i_2=s+1}^{\infty}\displaystyle\sum_{j_2=0}^{i_2}\|d_{j_1}^{t+p}(\theta_0)d_{i_1-j_1}^{t}(\theta_0)\eta_{t+p-j_1}\eta_{t+j_1-i_1}\|_{2+\nu} \notag \\ & \ \hspace{1.5cm} \times \|d_{j_2}^{t-h+p}(\theta_0)d_{i_2-j_2}^{t-h}(\theta_0)\eta_{t-h+p-j_2}\eta_{t-h+j_2-i_2}\|_{2+\nu} \times \alpha^{\frac{\nu}{2+\nu}}_{\Delta,\eta}\left([h/2]\right) \notag \\ \leq & \ C^{(2)}\displaystyle\sum_{i_1=s+1}^{[h/2]}\displaystyle\sum_{i_2=s+1}^{\infty}\rho^{i_1}\rho^{i_2}\left(\alpha^{\frac{\nu}{2+\nu}}_{\eta}([h/2])+\alpha^{\frac{\nu}{2+\nu}}_{\Delta}(\lfloor h/2\rfloor)\right) \notag \\ \leq & \ C^{(3)}\rho^{s}\left(\alpha^{\frac{\nu}{2+\nu}}_{\eta}([h/2])+\alpha^{\frac{\nu}{2+\nu}}_{\Delta}(\lfloor h/2\rfloor)\right).
\end{align}
On the other hand we have
\begin{align}
\label{eq:16}
   |\text{Cov}({Z^{h^+}_{t,s}(p)},Z_{t-h,s}(p))| \leq & \ \mathbb{E}|Z^{h^+}_{t,s}(p)Z_{t-h,s}(p)| + \mathbb{E}|Z^{h^{+}}_{t,s}(p)|\mathbb{E}|Z_{t-h,s}(p)| \notag \\ \leq & \ \left(\mathbb{E}|Z^{h^+}_{t,s}(p)|^{2}\right)^{\frac{1}{2}}\left(\mathbb{E}|Z_{t-h,s}(p)|^{2}\right)^{\frac{1}{2}} + \mathbb{E}|Z^{h^{+}}_{t,s}(p)|\mathbb{E}|Z_{t-h,s}(p)|\notag\\ \leq & \ C^{(4)}\rho^{[h/2]}\rho^{s},
\end{align}
 where $C^{(1)},C^{(2)}, C^{(3)}$ \text{and} $C^{(4)}$ represent arbitrary positive constants.
\subparagraph{$\diamond$ Case 2 : suppose that \( h \geq 0 \) and \( \lfloor h/2 \rfloor \leq s \).}\ \\ 
  Then, from the Cauchy-Schwarz inequality and the Lemma~\ref{lem:6}, it follows that there exists a positive constant \( C^{(5)} \) such that
\begin{align}
\label{eq:17}
    |\text{Cov}(Z_{t,s}(p),Z_{t-h,s}(p)| & \leq  \ \left(\mathbb{E}|Z_{t,s}(p)|^{2}\right)^{\frac{1}{2}}\left(\mathbb{E}|Z_{t-h,s}(p)|^{2}\right)^{\frac{1}{2}} + \mathbb{E}|Z_{t,s}(p)|\mathbb{E}|Z_{t-h,s}(p)| \notag \\
    & \leq C^{(5)} \rho^{s}
\end{align}
Thus, by combining Equations (\ref{eq:15}), (\ref{eq:16}) and (\ref{eq:17}), we obtain
\begin{align}
\label{eq:18}
     \displaystyle\sum_{h=0}^{\infty}|c_s^Z(h)| = &\displaystyle\sum_{h=0}^{\infty}|c_s^Z(h)|\mathds{1}_{[h/2]\leq s} + \displaystyle\sum_{h=0}^{\infty}|c_s^Z(h)|\mathds{1}_{[h/2]>s}\notag \\
     \leq & \displaystyle\sum_{h=0}^{2s+1}|c_s^Z(h)|+\displaystyle\sum_{h=2(s+1)}^{\infty}|c_s^Z(h)| \notag\\ \leq &
     \ M^{(1)}(2s+2)\rho^{s} + C^{3}\displaystyle\sum_{h=2s+2}^{\infty} \ \rho^{s}\left(\alpha^{\frac{\nu}{2+\nu}}_{\eta}([h/2])+\alpha^{\frac{\nu}{2+\nu}}_{\Delta}(\lfloor h/2\rfloor)\right)\notag \\
     \leq & \ M^{(1)}(2s+2)\rho^s + M^{(2)}\rho^s \displaystyle\sum_{h=2(s+2)}^{\infty}\left(\alpha^{\frac{\nu}{2+\nu}}_{\eta}([h/2])+\alpha^{\frac{\nu}{2+\nu}}_{\Delta}(\lfloor h/2\rfloor)\right) \xrightarrow[s\rightarrow \infty]{} 0.
\end{align}
By a similar argument, one also shows that $\sum_{h=-\infty}^0|c_s^Z(h)|\xrightarrow[s\rightarrow \infty]{} 0$, so that 
\begin{align}
\label{eq:19}
    \displaystyle\sum_{h=-\infty}^{\infty} |c_s^Z(h)|\xrightarrow[s\rightarrow \infty]{} 0.
\end{align}
Therefore, from the combination of Equations (\ref{eq:1800}), (\ref{eq:18}) and (\ref{eq:19}), we deduce that
\begin{align*}
  \displaystyle\sup_{n}\text{Var}\left(\frac{1}{\sqrt{n}}\displaystyle\sum_{t=1}^{n}(Z_{t,s}(p)-\mathbb{E}_{\theta_0}(Z_{t,s}(p)))\right) \xrightarrow[s\rightarrow \infty]{} 0.
\end{align*}
And the proof is complete using \citet[Corollary 7.7.1, p. 426]{anderson1971family}.
\qed

\subsection{Proof of the convergence of the variance matrix estimator}\label{proof_thm4}
 We proceed to demonstrate the proof of Theorem \ref{thm:4} by employing a series of Lemmas.

 We consider the regression of $\mathcal{Y}_t$ on the family $\{\mathcal{Y}_{t-1}, \ldots, \mathcal{Y}_{t-r}\}$ defined by
\begin{align}
    \label{eq:22}
    \mathcal{Y}_t = \sum_{i=1}^{r}\varphi_{r,i}\mathcal{Y}_{t-i}+u_{r,t}.    
\end{align}
 Denote 
\begin{align*}
    \underline{\bm{\varphi}}_r^{\star} & := (\varphi_1,\ldots,\varphi_r)\in\mathbb{R}^{N\times rN}, & \underline{\bm{\varphi}}_r & := (\varphi_{r,1}, \ldots, \varphi_{r,r})\in\mathbb{R}^{N\times rN}, \\
    \underline{\mathcal{Y}}_{r,t} & := \left(\mathcal{Y}_{t-1}',\ldots,\mathcal{Y}_{t-r}^{\prime}\right)'\in\mathbb{R}^{rN}, & \underline{\hat{\mathcal{Y}}}_{r,t} & := (\hat{\mathcal{Y}}_{t-1}',\ldots,\hat{\mathcal{Y}}_{t-r}')'\in\mathbb{R}^{rN}.
\end{align*}
We recall that $ \hat{\mathcal{Y}}_t=Y_t-(\hat{c}_{1,0}, \dots, \hat{c}_{N,0})'$ and additionally, we maintain the convention $\hat{\mathcal{Y}}_t = \mathcal{Y}_t = 0$ for $t\leq 0$ or $t>n$.
We also denote $\hat{\underline{\bm{\varphi}}}_r : = (\hat{\varphi}_{r,1},\dots,\hat{\varphi}_{r,r})$.
 When the values of $\mathcal{Y}_1, \ldots, \mathcal{Y}_n$ are known, the expressions for the least squares estimators of 
$\underline{\bm{\varphi}}_r$ 
and $\Sigma_{u_r}  := \Var(u_{r,t})$ are given by
\begin{align*}
 \underline{\check{\bm{\varphi}}}_r = \hat{\Sigma}_{\mathcal{Y},\underline{\mathcal{Y}}_r}\hat{\Sigma}_{\underline{\mathcal{Y}}_r}^{-1} \quad  \text{ and }  \quad \hat{\Sigma}_{\hat{u}_r} = \frac{1}{n}\displaystyle\sum_{t=1}^n\left(\mathcal{Y}_t-\check{\underline{\bm{\varphi}}}_r\underline{\mathcal{Y}}_{r,t}\right)\left(\mathcal{Y}_t-\check{\underline{\bm{\varphi}}}_r\underline{\mathcal{Y}}_{r,t}\right)', 
\end{align*} where
\begin{align*}
\hat{\Sigma}_{\mathcal{Y},\underline{\mathcal{Y}}_r} = \frac{1}{n}\displaystyle\sum_{t=1}^{n}\mathcal{Y}_t\underline{\mathcal{Y}}_{r,t}'  \quad \text{ and }  \quad \hat{\Sigma}_{\underline{\mathcal{Y}}_r} = \frac{1}{n}\displaystyle\sum_{t=1}^{n}\underline{\mathcal{Y}}_{r,t}\underline{\mathcal{Y}}_{r,t}'.
\end{align*}
 Here, when the values of $\mathcal{Y}_1, \ldots, \mathcal{Y}_n$ are unobserved, which is the case for us, the least squares estimators for $\underline{\bm{\varphi}}_r$ and $\Sigma_{u_r}$ are defined by
\begin{align*}
 \underline{\hat{\bm{\varphi}}}_r = \hat{\Sigma}_{\hat{\mathcal{Y}},\underline{\hat{\mathcal{Y}}}_r}\hat{\Sigma}_{\underline{\hat{\mathcal{Y}}}_r}^{-1} \quad \text{ and } \quad \hat{\Sigma}_{\hat{u}_r} = \frac{1}{n}\sum_{t=1}^n\left(\hat{\mathcal{Y}}_t-\hat{\underline{\bm{\varphi}}}_r\underline{\hat{\mathcal{Y}}}_{r,t}\right)\left(\hat{\mathcal{Y}}_t-\hat{\underline{\bm{\varphi}}}_r\underline{\hat{\mathcal{Y}}}_{r,t}\right)' ,
\end{align*} where
\begin{align*}
   \hat{\Sigma}_{\hat{\mathcal{Y}},\underline{\hat{\mathcal{Y}}}_r} = \frac{1}{n}\displaystyle\sum_{t=1}^{n}\hat{\mathcal{Y}}_t\underline{\hat{\mathcal{Y}}}_{r,t}'  \quad  \text{ and }  \quad \hat{\Sigma}_{\underline{\hat{\mathcal{Y}}}_r} = \frac{1}{n}\displaystyle\sum_{t=1}^{n}\underline{\hat{\mathcal{Y}}}_{r,t}\underline{\hat{\mathcal{Y}}}_{r,t}'.
\end{align*}
 Let 
\begin{align*}
    \Sigma_{\mathcal{Y}_t,{\underline{\mathcal{Y}}_{r}}} = \mathbb{E}\mathcal{Y}_t{\underline{\mathcal{Y}'}_{r,t}}, \  \Sigma_{\mathcal{Y}} = \mathbb{E}\mathcal{Y}_t{\mathcal{Y}_t}', \  \Sigma_{\underline{\mathcal{Y}}_r} = \mathbb{E}\underline{\mathcal{Y}}_{r,t}{\underline{\mathcal{Y}'}_{r,t}}, \, \  \hat{\Sigma}_{{\mathcal{Y}}} = \frac{1}{n}\displaystyle\sum_{t=1}^{n}{\mathcal{Y}}_{t}{\mathcal{Y}}_{t}' \ \text{and} \ \hat{\Sigma}_{\hat{\mathcal{Y}}} = \frac{1}{n}\displaystyle\sum_{t=1}^{n}{\hat{\mathcal{Y}}_{t}\hat{\mathcal{Y}}}_{t}'.
\end{align*}
In what follows, we will adopt the multiplicative matrix norm given by $$\|A\| = \sup_{\|x\|\leq 1} \|Ax\| = \rho^{1/2}(A^{\prime}A),$$ where, here $A:=(a_{i,j})$ is a matrix of arbitrary dimensions, $\|x\|$ denotes the Euclidean norm for vectors, and $\rho(\cdot)$ represents the spectral radius. This particular norm is chosen because it satisfies the inequality
\begin{align}
\label{eq:505}
\|A\|^2 \leq \sum_{i,j} a_{i,j}^2. 
\end{align}
The choice of this norm is critical for proving the upcoming Lemmas.

\begin{l1}\label{lem:9}
Under the assumptions of Theorem~\ref{thm:4},  we have
\begin{align*}
\sup_{r\geq 1} \max\left\{\left\|\Sigma_{\mathcal{Y},\underline{\mathcal{Y}}_r}\right\|,\left\|\Sigma_{\underline{\mathcal{Y}}_r}\right\|,\left\|\Sigma_{\underline{\mathcal{Y}}_r}^{-1}\right\|\right\}<\infty.
\end{align*}
\end{l1}
\subsubsection*{Proof.}
We will initiate by demonstrating that $\sup_{r\geq 1} \left\|\Sigma_{\mathcal{Y}_r}\right\|<\infty$.
With the previous notations, we have
 \begin{align*} 
 \underline{\mathcal{Y}}_{r,t}\underline{\mathcal{Y}}'_{r,t} =\left( \begin{array}{cccc}
\mathcal{Y}_{t-1}\mathcal{Y}'_{t-1} & \mathcal{Y}_{t-1}\mathcal{Y}'_{t-2} & \cdots & \mathcal{Y}_{t-1}\mathcal{Y}'_{t-r} \\ 
\mathcal{Y}_{t-2}\mathcal{Y}'_{t-1} & \mathcal{Y}_{t-2}\mathcal{Y}'_{t-2} & \cdots & \mathcal{Y}_{t-2}\mathcal{Y}'_{t-r} \\ 
\vdots & \vdots & \ddots & \vdots \\
\mathcal{Y}_{t-r}\mathcal{Y}'_{t-1} & \mathcal{Y}_{t-r}\mathcal{Y}'_{t-2} & \cdots & \mathcal{Y}_{t-r}\mathcal{Y}'_{t-r}
\end{array} \right)\in\mathbb{R}^{rN \times rN}.
\end{align*}
Hence, by stationarity we obtain
\begin{align*}
    \Sigma_{\underline{\mathcal{Y}}_r} = \left[\mathbb{E}\left(\mathcal{Y}_{t-u}\mathcal{Y}'_{t-v}\right)\right]_{u,v = 1,\dots,r} = \left[\mathbb{E}\left(\mathcal{Y}_{0}\mathcal{Y}'_{u-v}\right)\right]_{u,v = 1,\dots,r}  = \left[C(u-v)\right]_{u,v = 1,\dots,r}
\end{align*}
where $C(k) = \mathbb{E}(\mathcal{Y}_{0}\mathcal{Y}'_{k})\in\mathbb{R}^{N\times N}, k\in\mathbb{Z}$.  Subsequently, we introduce the spectral density of the stationary process $(\mathcal{Y}_t)_{t\in\mathbb{Z}}$ defined by
\begin{align*}
    f(\omega) = \frac{1}{2\pi}\sum_{k=-\infty}^{\infty} C(k)e^{i\omega k}, \ \forall \ \omega\in \mathbb{R}.
\end{align*}
Consider $\rho(\Sigma_{\underline{\mathcal{Y}}_r})$ as the spectral radius of the matrix $\Sigma_{\underline{\mathcal{Y}}_r}$, which corresponds to the eigenvector $\gamma^{(r)} := \left({\gamma_1^{(r)}}',\cdots,{\gamma_r^{(r)}}'\right)\in\mathbb{R}^{rN}$. Here, each $\gamma_j^{(r)}$ belongs to $\mathbb{R}^{N}$ for $j = 1,\dots,r$, and $\|\gamma_j^{(r)}\| =  1$.  

\noindent We have on one hand
\begin{align}
\label{eq:23}
   {\gamma^{(r)}}'\Sigma_{\underline{\mathcal{Y}}_r}\gamma^{(r)} =  {\gamma^{(r)}}'\rho(\Sigma_{\underline{\mathcal{Y}}_r}) \gamma^{(r)} = \rho(\Sigma_{\underline{\mathcal{Y}}_r}) \|\gamma^{(r)}\|^{2} = \rho(\Sigma_{\underline{\mathcal{Y}}_r}) = \|\Sigma_{\underline{\mathcal{Y}}_r}\|,
\end{align}
\noindent and on the other hand using the inversion formula,
\begin{align}
\label{eq:24}
 {\gamma^{(r)}}'\Sigma_{\underline{\mathcal{Y}}_r}\gamma^{(r)} = & \sum_{m,n = 1}^{r}{\gamma_m^{(r)}}'C(m-n)\gamma_n^{(r)} \notag \\
 = & \sum_{m,n = 1}^{r}{\gamma_m^{(r)}}'\left(\int_{[-\pi,\pi]}f(\omega)e^{-i\omega(m-n)}\mathrm{d}\omega\right)\gamma_{n}^{(r)} \notag \\ = & \int_{[-\pi,\pi]}\left(\sum_{m,n = 1}^{r}{\gamma_{m}^{(r)}}' f(\omega)e^{i\omega(n-m)}\gamma_{n}^{(r)}\right)\mathrm{d}\omega \notag\\
 = & \int_{[-\pi,\pi]}\left(\sum_{n=1}^r \gamma_{n}^{(r)}e^{i\omega n}\right)'f(\omega)\left(\sum_{m=1}^r \gamma_{m}^{(r)}e^{-i\omega m}\right)\mathrm{d}\omega \notag\\
 = & \int_{[-\pi,\pi]}\left(\sum_{n=1}^r \gamma_{n}^{(r)}e^{i\omega n}\right)'f(\omega)\left(\overline{\sum_{m=1}^r \gamma_{m}^{(r)}e^{i\omega m}}\right)\mathrm{d}\omega  .
\end{align}
Let \(\omega \in \mathbb{R}\). Consider the bilinear mapping, \(F : (X,Y) \in \mathbb{C}^{N} \times \mathbb{C}^{N} \mapsto F(X,Y) = X'f(\omega)\overline{Y}\). Given $(C(k))_{k\in\mathbb{Z}}$ is symmetric with respect to zero and Hermitian, it follows that $f(\omega)$ is also Hermitian, and consequently, its eigenvalues are real. Therefore, there exists a diagonal matrix \(D_\omega\) containing the eigenvalues of \(f(\omega)\) and a unitary matrix \(U_\omega\) such that \(f(\omega) = U_\omega D_\omega \overline{U_\omega}'\).

\noindent Moreover, since \(U\) is unitary, we also have for all $X\in\mathbb{C}^N$,
\begin{align}
\label{eq:25}
X'f(\omega)\overline{X} = X'\left(U_\omega D_\omega\overline{U_\omega}'\right)\overline{X} = \left(U_\omega^{'} X\right)'D_\omega\left(\overline{U_\omega^{'}X}\right)\leq \sup_{\omega}\|f(\omega)\|\|X\|^2.
\end{align}
Furthermore,
\begin{align}
\label{eq:26}
    \frac{1}{2\pi}\int_{[-\pi,\pi]}\left(\sum_{m=1}^r \gamma_m^{(r)}e^{i\omega m}\right)'\left(\overline{\sum_{n=1}^r \gamma_n^{(r)}e^{i\omega n}}\right)\mathrm{d}\omega = & \frac{1}{2\pi}\int_{[-\pi,\pi]}\left(\sum_{m=1}^r \sum_{n=1}^r \gamma_m^{(r)}\overline{\gamma_n^{(r)}}e^{i\omega(m-n)}\right)\mathrm{d}\omega \notag \\ = & \ \frac{1}{2\pi}\sum_{m=1}^r \sum_{n=1}^r \gamma_{m}^{(r)}\overline{\gamma_{n}^{(r)}}\int_{[-\pi,\pi]}e^{i\omega(m-n)}\mathrm{d}\omega \notag \\ = & \ \sum_{m=1}^r \sum_{n=1}^r \gamma_{m}^{(r)}\overline{\gamma_{n}^{(r)}}\delta_{mn} \notag \\  = & \ \|\gamma^{(r)}\|^2.
\end{align}
By combining Equations (\ref{eq:23}), (\ref{eq:24}), (\ref{eq:25}), and (\ref{eq:26}), along with the fact that the series $f(\cdot)$ is convergent, we arrive at
\begin{align*}
   \|\Sigma_{\underline{\mathcal{Y}}_r}\|\leq 2\pi\sup_{\omega\in\mathbb{R}}\|f(\omega)\|<\infty.
\end{align*}

\noindent Knowing that the eigenvalues of $\Sigma_{\mathcal{Y}_r}^{-1}$ are the inverses of the eigenvalues of $\Sigma_{\mathcal{Y}_r}$, it follows that $\rho(\Sigma_{\mathcal{Y}_r}^{-1})$ is equal to the inverse of the smallest eigenvalue of $\Sigma_{\mathcal{Y}_r}$ which is non-zero by hypothesis. Following the same reasoning, we also show that 
\begin{align*}
    \sup_{r\geq 1}\left\|\Sigma_{\mathcal{Y}_r}^{-1}\right\|< \infty.
\end{align*}
Furthermore, by noting that $\|\Sigma_{\mathcal{Y},\underline{\mathcal{Y}}_r}\|\leq \left\|\Sigma_{\underline{\mathcal{Y}}_{r+1}}\right\|$ (see in \cite{francq2003goodness}, p. 23,  for more details), allows us to conclude the proof.
\qed


\begin{l1}\label{lem:11}
We assume that the condition $\mathbb{E}|\eta_t|^{8+4\nu}<\infty$ holds and that Assumption $(\mathbf{A_4})$ is satisfied for some $\nu>0.$ Then, there exists a positive constant $C$ such that 
\begin{align*}
\sup_{s,\ell,e\in\mathbb{N}}\sum_{h=-\infty}^{\infty}\left|\Cov\left(\mathcal{Y}_1(m_1)\mathcal{Y}_{1+s}(m_2), \mathcal{Y}_{\ell+h}(m_1)\mathcal{Y}_{1+e+h}(m_2)\right)\right| < C,
\end{align*}
where $m_1, m_2 \in \{1, \ldots, N\}$ and $\mathcal{Y}_t(k)$ denotes the $k$-th element of the vector $\mathcal{Y}_t$.
\end{l1}

\subsubsection*{Proof.}
Consider the parameters \( h \in \mathbb{Z}\) and \( m_1, m_2 \in \{1, \ldots, N\} \). We have
\begin{align*}
\displaystyle\sum_{h=-\infty}^{\infty}\left|\Cov\left(\mathcal{Y}_1(m_1)\mathcal{Y}_{1+s}(m_2), \mathcal{Y}_{\ell+h}(m_1)\mathcal{Y}_{1+e+h}(m_2)\right)\right| \leq \displaystyle\sum_{h=-\infty}^{\infty}\sum_{j=1}^3|w_j(s,h,\ell,e,m_1,m_2)|
\end{align*}
where
\begin{align*}
w_1(s,h,\ell,e,m_1,m_2): = & \ + \text{Cov}(X_1 X_{1+m_1} X_{1+s} X_{1+s+m_2}, X_{\ell+h} X_{\ell+h+m_1} X_{1+s+h} X_{1+e+h+m_2}) \\ 
& - \mathbb{E}(X_{1+e+h} X_{1+e+h+m_2}) \text{Cov}(X_1 X_{1+m_1} X_{1+s} X_{1+s+m_2}, X_{\ell+h} X_{\ell+h+m_1})\\ 
& - \mathbb{E}(X_{\ell+h} X_{\ell+h+m_1}) \text{Cov}(X_1 X_{1+m_1} X_{1+s} X_{1+s+m_2}, X_{1+e+h} X_{1+e+h+m_2}),\\ 
w_2(s,h,\ell,e,m_1,m_2) : = & \ -\mathbb{E}(X_{1+s} X_{1+s+m_2}) \text{Cov}(X_1 X_{1+m_1}, X_{\ell+h} X_{\ell+h+m_1} X_{1+s+h} X_{1+s+h+m_2})\\ 
& + \mathbb{E}(X_{1+s} X_{1+s+m_2}) \mathbb{E}(X_{1+e+h} X_{1+e+h+m_2}) \text{Cov}(X_1 X_{1+m_1}, X_{\ell+h} X_{\ell+h+m_1})\\ 
& + \mathbb{E}(X_{1+s} X_{1+s+m_2}) \mathbb{E}(X_{\ell+h} X_{\ell+h+m_1}) \text{Cov}(X_1 X_{1+m_1}, X_{1+s+h} X_{1+s+h+m_2}),\\ 
w_3(s,h,\ell,e,m_1,m_2) : = & \ -\mathbb{E}(X_1 X_{1+m_1}) \text{Cov}(X_{1+s} X_{1+s+m_2}, X_{\ell+h} X_{\ell+h+m_1} X_{1+e+h} X_{1+e+h+m_2})\\ 
& + \mathbb{E}(X_1 X_{1+m_1}) \mathbb{E}(X_{1+e+h} X_{1+e+h+m_2}) \text{Cov}(X_{1+s} X_{1+s+m_2}, X_{\ell+h} X_{\ell+h+m_1})\\
& + \mathbb{E}(X_1 X_{1+m_1}) \mathbb{E}(X_{\ell+h} X_{\ell+h+m_1}) \text{Cov}(X_{1+s} X_{1+s+m_2}, X_{1+e+h} X_{1+e+h+m_2}).
\end{align*}
In what follows, we will focus on bounding $w_1(s,h,\ell,e,m_1,m_2)$. Similarly, $w_2(s,h,\ell,e,m_1,m_2)$ and $w_3(s,h,\ell,e,m_1,m_2)$ can be bounded in the same manner.
Let us put
\begin{align}
\label{eq:27}
\zeta_{i_1,\cdots,i_8}^{s,h,\ell,e}(m_1,m_2): = & \ \text{Cov}(d_{i_1}^{1} d_{i_2}^{1+m_1} d_{i_3}^{1+s} d_{i_4}^{1+s+m_2} \eta_{1-i_1} \eta_{1+m_1-i_2} \eta_{1+s-i_3} \eta_{1+s+m_2-i_4}, \notag \\ 
& \hspace{0.2cm} d_{i_5}^{\ell+h} d_{i_6}^{\ell+h+m_1} d_{i_7}^{1+e+h} d_{i_8}^{1+e+h+m_2} \eta_{\ell+h-i_5} \eta_{\ell+h+m-i_6} \eta_{1+e+h-i_7} \eta_{1+e+h+m_2-i_8}) \notag \\
= & \ \mathbb{E}(d_{i_1}^{1} d_{i_2}^{1+m_1} d_{i_3}^{1+s} d_{i_4}^{1+s+m_2} d_{i_5}^{\ell+h} d_{i_6}^{\ell+h+m_1} d_{i_7}^{1+e+h} d_{i_8}^{1+h+e+m_2}) \notag \\
 & \times\ \text{Cov}(\eta_{1-i_1} \eta_{1+m_1-i_2} \eta_{1+s-i_3} \eta_{1+s+m_2-i_4}, \eta_{\ell+h-i_5} \eta_{\ell+h+m_1-i_6} \eta_{1+e+h-i_7} \eta_{1+e+h+m_2-i_8}) \notag \\
 &+ \ \text{Cov}(d_{i_1}^{1} d_{i_2}^{1+m_1} d_{i_3}^{1+s} d_{i_4}^{1+s+m_2}, d_{i_5}^{\ell+h} d_{i_6}^{\ell+h+m_1} d_{i_7}^{1+e+h} d_{i_8}^{1+e+h+m_2}) \notag \\
&\times  \ \mathbb{E}(\eta_{1-i_1} \eta_{1+m-i_2} \eta_{1+s-i_3} \eta_{1+s+m_2-i_4}) \mathbb{E}(\eta_{\ell+h-i_5} \eta_{\ell+h+m_1-i_6} \eta_{1+e+h-i_7} \eta_{1+e+h+m_2-i_8}),
\end{align}
where we recall that $d_{k}^{t}:= c(\theta_0,\Delta_{t},\dots,\Delta_{t-k})$, for all $k,i\in\mathbb{Z}$ with $c(\theta_0,\Delta_{t},\dots,\Delta_{t-k})$ defined in Notation \ref{not:1} (See Equation (\ref{eq:not})).

Considering Equations (\ref{eq:4}), (\ref{eq:113}) and (\ref{eq:27}) we have
\begin{align*}
\text{Cov}(X_1X_{1+m_1}X_{1+s}X_{1+s+m_2},X_{\ell+h}X_{\ell+h+m_1}X_{1+e+h}X_{1+e+h+m_2}) = \sum_{0 \leq i_1, \cdots, i_8 \leq \infty} \zeta_{i_1, \cdots, i_8}^{s,h,\ell,e}(m_1,m_2)(\theta_0). 
\end{align*}
 In the remainder of this proof, we set $\epsilon_{i_1,t,m,i_2}^{(2)} := \eta_{t-i_1}\eta_{t+m-i_2}$ for all $t, i_1, i_2, m \in \mathbb{Z}$.

Using Lemma \ref{lem:6} and applying H{\"o}lder's inequality, there exists a positive constant $K$ and $\rho\in(0,1)$ such that 
\begin{align}
\label{eq:30}
\left|\mathbb{E}(d_{i_1}^{1}d_{i_2}^{1+m_1}d_{i_3}^{1+s}d_{i_4}^{1+s+m_1}d_{i_5}^{\ell+h}d_{i_6}^{\ell+h+m_1}d_{i_7}^{1+e+h}d_{i_8}^{1+h+e+m_2})\right| & \leq K\rho^{\sum_{k=1}^{8}i_k}
\end{align}
and by stationarity
\begin{align}
\label{eq:31}
\max \Big\{\left|\mathbb{E}(\epsilon_{i_1,1,m_1,i_2}^{(2)}  \epsilon_{i_3,1+s,m_2,i_4}^{(2)})\right|, \left|\mathbb{E}(\epsilon_{i_5,\ell+h,m_1,i_6}^{(2)}\epsilon_{i_7,1+e+h,m_2,i_8}^{(2)})\right|\Big\}  \ \leq \mathbb{E}|\eta_{0}|^4.
 <  \infty.
\end{align}
In view of Equations (\ref{eq:30}) and (\ref{eq:31}), there exists a positive constant $C$ such that
\begin{align}
\label{eq:pp1}
& \left|\text{Cov} (X_1X_{1+m_1}X_{1+s}X_{1+s+m_2},X_{\ell+h}X_{\ell+h+m_1}X_{1+e+h}X_{1+e+h+m_2})\right|\notag\\ & \leq \sum_{0 \leq i_1, \cdots, i_8 \leq \infty}\Bigg\{K\rho^{\sum_{k=1}^{8}i_k} \left|\text{Cov}\left(\epsilon_{i_1,1,m_1,i_2}^{(2)}\epsilon_{i_3,1+s,m_2,i_4}^{(2)},\epsilon_{i_5,\ell+h,m_1,i_6}^{(2)}\epsilon_{i_7,1+e+h,m_2,i_8}^{(2)}\right)\right| \notag\\ &  \hspace{2.5cm} + C\left|\text{Cov}(d_{i_1}^{1}d_{i_2}^{1+m_1}d_{i_3}^{1+s}d_{i_4}^{1+s+m_2},d_{i_5}^{\ell+h}d_{i_6}^{\ell+h+m_1}d_{i_7}^{1+e+h}d_{i_8}^{1+e+h+m_2})\right|\Bigg\}.
\end{align}
We will now bound the two terms in Equation \eqref{eq:pp1}.
\subparagraph{$\diamond$ Assuming $h \geq 0$, let us define $m_0 := m_1 \lor m_2$.}\ \\

  Firstly, we have 
\begin{align*}
 & \sum_{0 \leq i_1, \cdots, i_8 \leq \infty} C\left|\text{Cov}(d_{i_1}^{1}d_{i_2}^{1+m_1}d_{i_3}^{1+s}d_{i_4}^{1+s+m_1},d_{i_5}^{\ell+h}d_{i_6}^{\ell+h+m_1}d_{i_7}^{1+e+h}d_{i_8}^{1+e+h+m_2}\right| \leq C\sum_{k=1}^{9}v_k(s,h,\ell,e,m_1,m_2),
\end{align*}
 where
\begin{align*}
 v_1(s,h,\ell,e,m_1,m_2): = & \sum_{i_1 > \lfloor h/2 \rfloor} \sum_{0 \leq i_2, \dots, i_8 \leq \infty} \Big|\text{Cov}(d_{i_1}^{1} d_{i_2}^{1+m_1} d_{i_3}^{1+s} d_{i_4}^{1+s+m_2}, d_{i_5}^{\ell+h} d_{i_6}^{\ell+h+m_1} d_{i_7}^{1+e+h} d_{i_8}^{1+e+h+m_2})\Big|,\\
 v_2(s,h,\ell,em_1,m_2) : = & \sum_{i_2 > \lfloor h/2 \rfloor} \sum_{0 \leq i_1, i_3, \dots, i_8 \leq \infty} \Big|\text{Cov}(d_{i_1}^{1} d_{i_2}^{1+m_1} d_{i_3}^{1+s} d_{i_4}^{1+s+m_2}, d_{i_5}^{\ell+h} d_{i_6}^{\ell+h+m_1} d_{i_7}^{1+e+h} d_{i_8}^{1+e+h+m_2})\Big|, \\
  v_3(s,h,\ell,e,m_1,m_2) : =  & \sum_{i_3 > \lfloor h/2 \rfloor} \sum_{0 \leq i_1, i_2, i_4, \dots, i_8 \leq \infty} \Big|\text{Cov}(d_{i_1}^{1} d_{i_2}^{1+m_1} d_{i_3}^{1+s} d_{i_4}^{1+s+m_2}, d_{i_5}^{\ell+h} d_{i_6}^{\ell+h+m_1} d_{i_7}^{1+e+h} d_{i_8}^{1+e+h+m_2})\Big| ,\\
   v_4(s,h,\ell,em_1,m_2) : = & \sum_{i_4 > \lfloor h/2 \rfloor} \sum_{0 \leq i_1, i_2, i_3, i_5, \dots, i_8 \leq \infty} \Big|\text{Cov}(d_{i_1}^{1} d_{i_2}^{1+m_1} d_{i_3}^{1+s} d_{i_4}^{1+s+m_2}, d_{i_5}^{\ell+h} d_{i_6}^{\ell+h+m_1} d_{i_7}^{1+e+h} d_{i_8}^{1+e+h+m_2})\Big|, \\
   v_5(s,h,\ell,e,m_1,m_2) : = & \sum_{i_5 > \lfloor h/2 \rfloor} \sum_{0 \leq i_1, \dots, i_4, i_6, i_7, i_8 \leq \infty} \Big|\text{Cov}(d_{i_1}^{1} d_{i_2}^{1+m_1} d_{i_3}^{1+s} d_{i_4}^{1+s+m_2}, d_{i_5}^{\ell+h} d_{i_6}^{\ell+h+m_1} d_{i_7}^{1+e+h} d_{i_8}^{1+e+h+m_2})\Big|,\\
   v_6(s,h,\ell,e,m_1,m_2) : = & \sum_{i_6 > \lfloor h/2 \rfloor} \sum_{0 \leq i_1, \dots, i_5, i_7, i_8 \leq \infty} \Big|\text{Cov}(d_{i_1}^{1} d_{i_2}^{1+m_1} d_{i_3}^{1+s} d_{i_4}^{1+s+m_2}, d_{i_5}^{\ell+h} d_{i_6}^{\ell+h+m_1} d_{i_7}^{1+e+h} d_{i_8}^{1+e+h+m_2})\Big|,\\
v_7(s,h,\ell,e,m_1,m_2) : = & \sum_{i_7 > \lfloor h/2 \rfloor} \sum_{0 \leq i_1, \dots, i_6, i_8 \leq \infty} \Big|\text{Cov}(d_{i_1}^{1} d_{i_2}^{1+m_1} d_{i_3}^{1+s} d_{i_4}^{1+s+m_2}, d_{i_5}^{\ell+h} d_{i_6}^{\ell+h+m_1} d_{i_7}^{1+e+h} d_{i_8}^{1+e+h+m_2})\Big|,\\
v_8(s,h,\ell,e,m_1,m_2) : = & \sum_{i_8 > \lfloor h/2 \rfloor} \sum_{0 \leq i_1, \dots, i_7 \leq \infty} \Big|\text{Cov}(d_{i_1}^{1} d_{i_2}^{1+m_1} d_{i_3}^{1+s} d_{i_4}^{1+s+m_2}, d_{i_5}^{\ell+h} d_{i_6}^{\ell+h+m_1} d_{i_7}^{1+e+h} d_{i_8}^{1+e+h+m_2})\Big|,\\
v_9(s,h,\ell,e,m_1,m_2) : = & \sum_{0 \leq i_1, \dots, i_8 \leq \lfloor h/2 \rfloor} \Big|\text{Cov}(d_{i_1}^{1} d_{i_2}^{1+m_1} d_{i_3}^{1+s} d_{i_4}^{1+s+m_2}, d_{i_5}^{\ell+h} d_{i_6}^{\ell+h+m_1} d_{i_7}^{1+e+h} d_{i_8}^{1+e+h+m_2})\Big|.
\end{align*}
 Using an argument similar to that in Equation (\ref{eq:14}), there exists a family of positive constants \((\mathcal{C}_i)_{1\leq i\leq 8}\) such that
\begin{align}
\label{eq:333}
\sup_{s,\ell,e\in\mathbb{N}}\displaystyle\sum_{h=0}^{\infty}v_i(s,h,\ell,e,m_1,m_2) \leq \displaystyle\sum_{h=0}^{\infty}\mathcal{C}_i\rho^{h/2} <\infty. \
\end{align}
 To bound the term \(v_9(s,h,\ell,e,m_1,m_2)\), we will apply the inequality from \cite{davydov1968convergence}, as presented in Lemma \ref{lem:5}. Additionally, we introduce the positive constants \(K_0\) and \(C_0\), which both represent the universal constant specified in Lemma \ref{lem:5}, respectively.
\subparagraph{$\diamond$ Let us suppose $\lfloor h/2 \rfloor\geq s+m_0$.}\ \\
 We have, using the Lemma \ref{lem:5},
\begin{align}
\label{eq:p1}
    & \sum_{\lfloor h/2 \rfloor=s+m_0}^{\infty} \displaystyle\sum_{i_1,\cdots,i_8=0}^{\lfloor h/2 \rfloor} \Big|\text{Cov}(d_{i_1}^{1}d_{i_2}^{1+m_1}d_{i_3}^{1+s}d_{i_4}^{1+s+m_2},d_{i_5}^{\ell+h}d_{i_6}^{\ell+h+m_1}d_{i_7}^{1+e+h}d_{i_8}^{1+e+h+m_2})\Big| \notag \\ & \leq C_0\sum_{\lfloor h/2 \rfloor=s+m_0}^{\infty}\displaystyle\sum_{i_1,\cdots,i_8=0}^{\lfloor h/2 \rfloor}\|d_{i_1}^{1}d_{i_2}^{1+m_1}d_{i_3}^{1+s}d_{i_4}^{1+s+m_2}\|_{2+\nu}\|d_{i_5}^{\ell+h}d_{i_6}^{\ell+h+m_1}d_{i_7}^{1+e+h}d_{i_8}^{1+e+h+m_2}\|_{2+\nu} \notag \\ & \hspace{2.5cm} \times \alpha_{\Delta}^{\frac{\nu}{2+\nu}}(\min\{\ell,1+e\}+h-m_2-s-\max\{i_5,i_7\}-1) \notag \\ & \leq \left(\sum_{i_1,\cdots,i_8=0}^{\lfloor h/2 \rfloor}\rho^{\sum_{k=1}^{8}i_k}\right)\left(\sum_{\lfloor h/2 \rfloor=s+m_0}^{\infty}\alpha_{\Delta}^{\frac{\nu}{2+\nu}}(\lfloor h/2 \rfloor +\min\{\ell,1+e\}-m_0-s-1)\right).
\end{align}
To deal with the terms obtained for $\lfloor h/2 \rfloor<s+m_0$, we can write the following decomposition by applying the equality \(\Cov(UV,WZ) = \Cov(UW,VZ)+\mathbb{E}(UW)\mathbb{E}(VZ)-\mathbb{E}(UV)\mathbb{E}(WZ)\) for real random variables \(U\), \(V\), \(W\) and \(Z\) so that
\begin{align}
\label{eq:p999}
 & \text{Cov}(d_{i_1}^{1} d_{i_2}^{1+m_1} d_{i_3}^{1+s} d_{i_4}^{1+s+m_2}, d_{i_5}^{\ell+h} d_{i_6}^{\ell+h+m_1} d_{i_7}^{1+e+h} d_{i_8}^{1+e+h+m_2}) \notag\\ 
 = & -\text{Cov}\left(d_{i_1}^{1} d_{i_2}^{1+m_1}, d_{i_3}^{1+s} d_{i_4}^{1+s+m_2}\right) \text{Cov}\left(d_{i_5}^{\ell+h} d_{i_6}^{\ell+h+m_1}, d_{i_7}^{1+e+h} d_{i_8}^{1+e+h+m_2}\right) \notag\\ 
 & -\text{Cov}\left(d_{i_1}^{1} d_{i_2}^{1+m_1}, d_{i_3}^{1+s} d_{i_4}^{1+s+m_2}\right) \mathbb{E}\left(d_{i_5}^{\ell+h} d_{i_6}^{\ell+h+m_1}\right) \mathbb{E}\left(d_{i_7}^{1+e+h} d_{i_8}^{1+e+h+m_2}\right) \notag\\ 
 & - \text{Cov}\left(d_{i_5}^{\ell+h} d_{i_6}^{\ell+h+m_1}, d_{i_7}^{1+e+h} d_{i_8}^{1+e+h+m_2}\right) \mathbb{E}\left(d_{i_1}^{1} d_{i_2}^{1+m_1}\right) \mathbb{E}\left(d_{i_3}^{1+s} d_{i_4}^{1+s+m_2}\right) \notag\\ 
 & + \text{Cov}\left(d_{i_1}^{1+m_1} d_{i_2}^{1+m_1}, d_{i_5}^{\ell+h} d_{i_6}^{\ell+h+m_1}\right) \text{Cov}\left(d_{i_3}^{1+s} d_{i_4}^{1+s+m_2}, d_{i_7}^{1+e+h} d_{i_8}^{1+e+h+m_2}\right) \notag\\ 
 & + \text{Cov}\left(d_{i_1}^1 d_{i_2}^{1+m_1}, d_{i_7}^{1+s+h} d_{i_8}^{1+s+h+m_2}\right) \mathbb{E}\left(d_{i_3}^{1+s} d_{i_4}^{1+s+m_2}\right) \mathbb{E}\left(d_{i_5}^{\ell+h} d_{i_6}^{\ell+h+m_1}\right) \notag\\ 
 & + \text{Cov}(d_{i_1}^{1} d_{i_2}^{1+m_1} d_{i_5}^{\ell+h} d_{i_6}^{\ell+h+m_1}, d_{i_3}^{1+s} d_{i_4}^{1+s+m_2} d_{i_7}^{1+e+h} d_{i_8}^{1+e+h+m_2}).
\end{align}
Noting that in the previous decomposition of Equation (\ref{eq:p999}), we also assume that $e > \ell$ and $e > s$. The others case can be handled in a similar manner. In what follows, we will bound each term in Equation (\ref{eq:p999}) by applying the inequality from \cite{davydov1968convergence}, as stated in Lemma \ref{lem:5}.

By Lemma \ref{lem:5},note that for $\lfloor h/2 \rfloor < 1+e-\ell-m_1$ we have
\begin{align}
\label{eq:bad1}
\left|\text{Cov}\left(d_{i_5}^{\ell+h} d_{i_6}^{\ell+h+m_1}, d_{i_7}^{1+e+h} d_{i_8}^{1+e+h+m_2}\right)\right| &\leq C_0 \|d_{i_5}^{1+h} d_{i_6}^{\ell+h+m_1}\|_{2+\nu} \|d_{i_7}^{1+e+h} d_{i_8}^{1+e+h+m_2}\|_{2+\nu} \notag \\
&\hspace{1.5cm} \times \alpha_{\Delta}^{\frac{\nu}{2+\nu}}\{1 + e - \ell - i_7 - m_1\},
\end{align}
and for  $\lfloor h/2 \rfloor \geq 1+e-\ell-m_1$,,
\begin{align}
\label{eq:bad2}
\left|\text{Cov}\left(d_{i_5}^{\ell+h} d_{i_6}^{\ell+h+m_1}, d_{i_7}^{1+e+h} d_{i_8}^{1+e+h+m_2}\right)\right| &= \left|\text{Cov}\left(d_{i_7}^{1+e+h} d_{i_8}^{1+e+h+m_2}, d_{i_5}^{\ell+h} d_{i_6}^{\ell+h+m_1}\right)\right| \notag\\ 
&\leq C_0 \|d_{i_5}^{\ell+h} d_{i_6}^{\ell+h+m_1}\|_{2+\nu} \|d_{i_7}^{1+e+h} d_{i_8}^{1+e+h+m_2}\|_{2+\nu} \notag \\
&\hspace{1.5cm} \times \alpha_{\Delta}^{\frac{\nu}{2+\nu}} \{\ell - 1 - e - i_5 - m_2\}.
\end{align}
 Thus, from Equations (\ref{eq:bad1}), (\ref{eq:bad2}) and the fact that  \(0 \leq i_5, i_7 \leq \lfloor h/2 \rfloor\), it follows 
\begin{align}
\label{eq:min1}
\left|\text{Cov}\left(d_{i_5}^{\ell+h} d_{i_6}^{\ell+h+m_1}, d_{i_7}^{1+e+h} d_{i_8}^{1+e+h+m_2}\right)\right| &\leq 2C_0 \|d_{i_5}^{\ell+h} d_{i_6}^{\ell+h+m_1}\|_{2+\nu} \|d_{i_7}^{1+s+h} d_{i_8}^{1+s+h+m_2}\|_{2+\nu} \notag \\
&\hspace{0.5cm} \times\alpha_{\Delta}^{\frac{\nu}{2+\nu}}\{\min\{\ell-e-1-m_0,1+e-\ell-m_0\}-\lfloor h/2 \rfloor\}.
\end{align}
Using Lemma \ref{lem:5}, Equation (\ref{eq:min1}) and the fact that \(\alpha_{\Delta}(k) = 1/4\) for \(k \leq 0\), we obtain:
\subparagraph{$\diamond$ First term in Equation (\ref{eq:p999})} 
\begin{align}
\label{eq:p2}
& \sum_{\lfloor h/2 \rfloor =0}^{s+m_0-1}\sum_{i_1,\cdots,i_8=0}^{\lfloor h/2 \rfloor}\left|\text{Cov}\left(d_{i_1}^{1}d_{i_2}^{1+m_1},d_{i_3}^{1+s}d_{i_4}^{1+s+m_2}\right)\text{Cov}\left(d_{i_5}^{\ell+h}d_{i_6}^{\ell+h+m_1},d_{i_7}^{1+e+h}d_{i_8}^{1+e+h+m_2}\right)\right| \notag\\  & \leq 2C_0^2\sum_{\lfloor h/2 \rfloor=0}^{s+m_0-1}\sum_{i_1,\cdots,i_8=0}^{\lfloor h/2 \rfloor}\Big\{\|d_{i_1}^{1}d_{i_2}^{1+m_1}\|_{2+\nu}\|d_{i_5}^{1+s}d_{i_4}^{1+s+m_2}\|_{2+\nu}\Big\}\alpha_{\Delta}^{\frac{\nu}{2+\nu}}\{s-m_1-i_3\} \notag\\ & \hspace{1cm} \times \Big\{\|d_{i_5}^{1+h}d_{i_6}^{\ell+h+m_1}\|_{2+\nu}\|d_{i_7}^{1+e+h}d_{i_8}^{1+e+h+m_2}\|_{2+\nu}\Big\}\alpha_{\Delta}^{\frac{\nu}{2+\nu}} \{\min\{\ell-e-1-m_0,1+e-\ell-m_0\}-\lfloor h/2 \rfloor\} \notag\\ & \leq 2C_0^2C\left(\sum_{i_1,\cdots,i_8=0}^{\infty}\rho^{\sum_{k=1}^{8}i_k}\right) \left(\sum_{\lfloor h/2 \rfloor=0}^{s+m_0-1} \alpha_{\Delta}^{\frac{\nu}{2+\nu}}\{\min\{\ell-e-1-m_0,1+e-\ell-m_0\}-\lfloor h/2 \rfloor\}\right). 
\end{align}
\subparagraph{$\diamond$ Second term in Equation (\ref{eq:p999})} 
\begin{align}
\label{eq:p3}
& \sum_{\lfloor h/2 \rfloor=0}^{s+m_0-1}\sum_{i_1,\cdots,i_8=0}^{\lfloor h/2 \rfloor}\left|\text{Cov}\left(d_{i_1}^{1}d_{i_2}^{1+m_1},d_{i_3}^{1+s}d_{i_4}^{1+s+m_2}\right)\mathbb{E}\left(d_{i_5}^{\ell+h}d_{i_6}^{\ell+h+m_1}\right)\mathbb{E}\left(d_{i_7}^{1+e+h}d_{i_8}^{1+e+h+m_2}\right)\right|\ \notag\\ & \leq  C_0\sum_{\lfloor h/2 \rfloor=0}^{s+m_0-1} \sum_{i_1,\cdots,i_8=0}^{\lfloor h/2 \rfloor}\Big\{\|d_{i_1}^{1}d_{i_2}^{1+m_1}\|_{2+\nu}\|d_{i_3}^{1+s}d_{i_4}^{1+s+m_2}\|_{2+\nu}\left\|d^{\ell+h}_{i_5}\right\|_4\left\|d^{\ell+h+m_1}_{i_4}\right\|_4 \left\|d^{1+e+h}_{i_7}\right\|_4\left\|d^{1+e+h+m_2}_{i_7}\right\|_4 \Big\} \notag\\ &\hspace{2.5cm} \times \alpha_{\Delta}^{\frac{\nu}{2+\nu}}\{s-m_1-i_3\} \notag\\ & \leq C_0C\left(\sum_{i_1,\cdots,i_8=0}^{\infty}\rho^{\sum_{k=1}^{8}i_k}\right)\left(\sum_{\lfloor h/2 \rfloor=0}^{s+m_0-1}\alpha_{\Delta}^{\frac{\nu}{2+\nu}}\{s-m_0-\lfloor h/2 \rfloor\}\right).
\end{align}
\subparagraph{$\diamond$ Third term in Equation (\ref{eq:p999})} 
\begin{align}
\label{eq:p4}
& \sum_{\lfloor h/2 \rfloor=0}^{s+m_0-1}\sum_{i_1,\cdots,i_8=0}^{\lfloor h/2 \rfloor}\left|\text{Cov}\left(d_{i_5}^{\ell+h}d_{i_6}^{\ell+h+m_1},d_{i_7}^{1+e+h}d_{i_8}^{1+e+h+m_2}\right) \mathbb{E}\left(d_{i_1}^{1}d_{i_2}^{1+m_1}\right)\mathbb{E}\left(d_{i_3}^{1+s}d_{i_4}^{1+s+m_2}\right)\right| \notag\\ & \leq 2C_0\sum_{\lfloor h/2 \rfloor=0}^{s+m_0-1}\sum_{i_1,\cdots,i_8=0}^{\lfloor h/2\rfloor} \Big\{\|d_{i_5}^{\ell+h}d_{i_6}^{\ell+h+m_1}\|_{2+\nu}\|d_{i_7}^{1+e+h}d_{i_8}^{1+e+h+m_2}\|_{2+\nu}\left\|d^{1}_{i_1}\right\|_4\left\|d^{1+m_1}_{i_2}\right\|_4 \left\|d^{1+s}_{i_3}\right\|_4\left\|d^{1+s+m_2}_{i_4}\right\|_4 \Big\} \notag\\ & \hspace{2.5cm} \times \alpha_{\Delta}^{\frac{\nu}{2+\nu}}\{\min\{\ell-e-1-m_0,1+e-\ell-m_0\}-\lfloor h/2 \rfloor\}\notag\\ & \leq 2C_0C\left(\sum_{i_1,\cdots,i_8=0}^{\infty}\rho^{\sum_{k=1}^{8}i_k}\right)\left(\sum_{\lfloor h/2 \rfloor=0}^{s+m_0-1} \alpha_{\Delta}^{\frac{\nu}{2+\nu}}\{\min\{\ell-e-1-m_0,1+e-\ell-m_0\}-\lfloor h/2 \rfloor\}\right).
\end{align}
\subparagraph{$\diamond$ Fourth term in Equation (\ref{eq:p999})} 
\begin{align}
\label{eq:p5}
&\sum_{\lfloor h/2 \rfloor=0}^{s+m_0-1}\sum_{i_1,\cdots,i_8=0}^{\lfloor h/2 \rfloor}\left|\text{Cov}\left(d_{i_1}^{1}d_{i_2}^{1+m_1},d_{i_5}^{\ell+h}d_{i_6}^{\ell+h+m_1}\right)\text{Cov}\left(d_{i_3}^{1+s}d_{i_4}^{1+s+m_2},d_{i_7}^{1+e+h}d_{i_8}^{1+e+h+m_2}\right)\right|\notag\\ & \leq \sum_{\lfloor h/2 \rfloor=0}^{s+m_0-1} C_0^2\sum_{i_1,\cdots,i_8=0}^{\lfloor h/2 \rfloor}\Big\{\|d_{i_1}^{1}d_{i_2}^{1+m_1}\|_{2+\nu}\|d_{i_7}^{1+e+h}d_{i_8}^{1+e+h+m_2}\|_{2+\nu}\Big\}\alpha_{\Delta}^{\frac{\nu}{2+\nu}}\{h+\ell-i_5-m_1-1\} \notag \\ &\hspace{2cm} \times \Big\{\|d_{i_3}^{1+s}d_{i_4}^{1+s+m_2}\|_{2+\nu}\|d_{i_7}^{1+e+h}d_{i_8}^{1+e+h+m_2}\|_{2+\nu}\Big\} \alpha_{\Delta}^{\frac{\nu}{2+\nu}}\{e-s+h-i_7-m_2\} \notag\\ & \leq 
C_0^2C\left(\sum_{i_1,\cdots,i_8=0}^{\infty}\rho^{\sum_{k=1}^{8}i_k}\right)\sum_{\lfloor h/2 \rfloor=0}^{s+m_0-1}\alpha_{\Delta}^{\frac{\nu}{2+\nu}}\{\lfloor h/2 \rfloor+\ell-m_0-1\} \alpha_{\Delta}^{\frac{\nu}{2+\nu}}\{e-s+\lfloor h/2 \rfloor-m_0\}\notag\\ & \leq 
C_0^2C\left(\sum_{i_1,\cdots,i_8=0}^{\infty}\rho^{\sum_{k=1}^{8}i_k}\right)\left(\sum_{\lfloor h/2 \rfloor=0}^{s+m_0-1}\alpha_{\Delta}^{\frac{\nu}{2+\nu}}\{e-s+\lfloor h/2 \rfloor-m_0-1\}\right).
\end{align}
\subparagraph{$\diamond$ Fifth term in Equation (\ref{eq:p999})}
\begin{align}
\label{eq:p6}
&\sum_{\lfloor h/2 \rfloor=0}^{s+m_0-1}\sum_{i_1,\cdots,i_8=0}^{\lfloor h/2 \rfloor}\left|\text{Cov}\left(d_{i_1}^1d_{i_2}^{1+m_1},d_{i_7}^{1+e+h}d_{i_8}^{1+e+h+m_2}\right)\mathbb{E}\left(d_{i_3}^{1+s}d_{i_4}^{1+s+m_2}\right)\mathbb{E}\left(d_{i_5}^{\ell+h}d_{i_6}^{\ell+h+m_1}\right)\right|\notag\\ & \leq \sum_{\lfloor h/2 \rfloor=0}^{s+m_0-1} C_0\sum_{i_1,\cdots,i_8=0}^{\lfloor h/2 \rfloor}\Big\{\|d_{i_1}^{1}d_{i_2}^{1+m_1}\|_{2+\nu}\|d_{i_7}^{1+e+h}d_{i_8}^{1+e+h+m_2}\|_{2+\nu}\left\|d^{1+s}_{i_3}\right\|_4\left\|d^{1+s+m_2}_{i_2}\right\|_4 \left\|d^{\ell+h}_{i_5}\right\|_4\left\|d^{\ell+h+m_1}_{i_8}\right\|_4 \Big\}\notag \\ & \hspace{2cm} \times \alpha_{\Delta}^{\frac{\nu}{2+\nu}}\{e+h-m_1-i_7\} \notag\\ & \leq C_0C\left(\sum_{i_1,\cdots,i_8=0}^{\infty}\rho^{\sum_{k=1}^{8}i_k}\right)\left(\sum_{\lfloor h/2 \rfloor=0}^{s+m_0-1} \alpha_{\Delta}^{\frac{\nu}{2+\nu}}\{e+\lfloor h/2\rfloor-m_0\}\right).
\end{align}

\subparagraph{$\diamond$ To control the sixth term, we distinguish between the subcases \(\lfloor 3h/2 \rfloor \leq s -\ell-m_1\) and \(\lfloor 3h/2 \rfloor > s - \ell-m_1\).}
\subparagraph{$\circ$ First subcase: \(\lfloor 3h/2 \rfloor \leq s - \ell-m_1\).}\ \\

 We have 
\begin{align}
\label{eq:p7}
& \sum_{\substack{\lfloor h/2 \rfloor=0 \\ 0 \leq \lfloor 3h/2 \rfloor \leq s - \ell-m_1}}^{s+m_0-1} \sum_{i_1,\cdots,i_8=0}^{\lfloor h/2 \rfloor} \left|\text{Cov}(d_{i_1}^{1} d_{i_2}^{1+m_1} d_{i_5}^{\ell+h} d_{i_6}^{\ell+h+m_1}, d_{i_3}^{1+s} d_{i_4}^{1+s+m_2} d_{i_7}^{1+e+h} d_{i_8}^{1+e+h+m_2})\right| \notag\\ 
& \leq C_0 \sum_{\substack{\lfloor h/2 \rfloor=0 \\ 0 \leq \lfloor 3h/2 \rfloor \leq s - \ell-m_1}}^{s+m_0-1} \sum_{i_1,\cdots,i_8=0}^{\lfloor h/2 \rfloor} \|d_{i_1}^{1} d_{i_2}^{1+m_1} d_{i_5}^{\ell+h} d_{i_6}^{\ell+h+m_2}\|_{2+\nu} \|d_{i_3}^{1+s} d_{i_4}^{1+s+m_2} d_{i_7}^{1+e+h} d_{i_8}^{1+e+h+m_2}\|_{2+\nu} \notag\\ 
& \hspace{2cm} \times \alpha_{\Delta}^{\frac{\nu}{2+\nu}} \{s - \ell - h - m_1 - i_3 + 1\} \notag\\ 
& \leq C_0C\left(\sum_{i_1,\cdots,i_8=0}^{\infty} \rho^{\sum_{k=1}^{\infty} i_k}\right) \left(\sum_{ \lfloor 3h/2 \rfloor = 0}^{s-\ell-m_1} \alpha_{\Delta}^{\frac{\nu}{2+\nu}} \{s - \ell - m_0 - \lfloor 3h/2 \rfloor\}\right).
\end{align}
\subparagraph{$\circ$ Second subcase: \(\lfloor 3h/2 \rfloor > s - \ell-m_1\).}\ \\

 In this second subcase, we use a decomposition similar to Equation (\ref{eq:p999}) to interchange the terms \(d^{\ell+h+m_1}_{i_6}\) and \(d^{1+s}_{i_3}\). 
\begin{align}
\label{eq:p9999}
& \Cov\left(d^{1}_{i_1} d^{1+m_1}_{i_2} d^{\ell+h}_{i_5}d^{\ell+h+m_1}_{i_6}, d^{1+s}_{i_3}d^{1+s+m_2}_{i_4} d^{1+s+h}_{i_7} d^{1+s+h+m_2}_{i_8}\right) \notag\\ 
& = \Cov\left(d_{i_1}^{1} d_{i_2}^{1+m_1} d_{i_5}^{\ell+h} d_{i_3}^{1+s}, d_{i_6}^{\ell+h+m_1} d_{i_4}^{1+s+m_2} d_{i_7}^{1+e+h} d_{i_8}^{1+e+h+m_2}\right) \notag \\
&  + \Cov(d^{1}_{i_1} d^{1+m_1}_{i_2} d^{\ell+h}_{i_5}, d^{1+s}_{i_3}) \Cov\left(d^{\ell+h+m_1}_{i_6}, d^{1+s+m_2}_{i_4} d^{1+e+h}_{i_7} d^{1+e+h+m_2}_{i_8}\right) \notag\\ 
& + \Cov\left(d^{1}_{i_1} d^{1+m_1}_{i_2} d^{\ell+h}_{i_5}, d^{1+s}_{i_3}\right) \mathbb{E}\left(d^{\ell+h+m_1}_{i_6}\right) \mathbb{E}\left(d^{1+s+m_2}_{i_4} d^{1+e+h}_{i_7} d^{1+e+h+m_2}_{i_8}\right) \notag\\ 
& + \Cov\left(d^{1+s+m_2}_{i_4} d^{1+e+h}_{i_7} d^{1+e+h+m_2}_{i_8}, d^{\ell+h+m_1}_{i_6}\right) \mathbb{E}\left(d^{1}_{i_1} d^{1+m_1}_{i_2} d^{\ell+h}_{i_5}\right) \mathbb{E}\left(d^{1+s}_{i_3}\right) \notag\\ 
& - \Cov\left(d^{1}_{i_1} d^{1+m_1}_{i_2} d^{\ell+h}_{i_5}, d^{\ell+h+m_1}_{i_6}\right) \Cov\left(d^{1+s}_{i_3}, d^{1+s+m_2}_{i_4} d^{1+e+h}_{i_7} d^{1+e+h+m_2}_{i_8}\right) \notag\\ 
& - \Cov\left(d^{1+s}_{i_3}, d^{1+s+m_2}_{i_4} d^{1+e+h}_{i_7} d^{1+e+h+m_2}_{i_8}\right) \mathbb{E}\left(d^1_{i_1} d^{1+m_1}_{i_2} d^{\ell+h}_{i_5}\right) \mathbb{E}\left(d^{\ell+h+m_1}_{i_6}\right) \notag\\ 
& - \Cov\left(d^1_{i_1} d^{1+m_1}_{i_2} d^{\ell+h}_{i_5}, d^{\ell+h+m_1}_{i_6}\right) \mathbb{E}\left(d^{1+s}_{i_3}\right) \mathbb{E}\left(d^{1+s+m_2}_{i_4} d^{1+e+h}_{i_7} d^{1+e+h+m_2}_{i_8}\right).
\end{align}
We focus on deriving an upper bound for the first term in Equation (\ref{eq:p9999}), as the remaining terms can be bounded similarly. To simplify the application of Lemma \ref{lem:5}, we define the following set:
\begin{align*}
    \mathcal{W}(h) := \bigg\{ 
    (i_3, i_4, i_7) \in \{0, \ldots, \lfloor h/2 \rfloor\}^3 : \ 
    & \lfloor 3h/2 \rfloor > s - \ell - m_1, \\
    & \ell + h < \min\big\{1 + s + m_2 - m_1 - i_4, 1 + s - i_3\big\}, \\
    & \ell + m_1 < 1 + e - i_7, \\
    & 0 \leq \lfloor h/2 \rfloor \leq s + m_0 - 1
    \bigg\}.
\end{align*}
The conditions in \(\mathcal{W}(h)\) ensure proper indices for identifying the $\sigma$-algebra generated when applying Lemma \ref{lem:5}. Alternatively, selecting appropriate indices (choosing the largest from the past and the smallest from the future) suffices to apply Davydov's inequality without requiring additional decomposition. Using Lemma \ref{lem:5} we therefore have
\begin{align}
\label{eq:p88}
& \sum_{(i_3,i_4,i_7)\in\mathcal{W}(h)} \sum_{i_1,\cdots,i_8=0}^{\lfloor h/2 \rfloor} \left|\Cov\left(d^{1}_{i_1}d^{1+m_1}_{i_2}d^{\ell+h}_{i_5}d^{1+s}_{i_3},d^{\ell+h+m_1}_{i_6}d^{1+s+m_2}_{i_4}d^{1+e+h}_{i_7}d^{1+e+h+m_2}_{i_8}\right)\right| \notag \\
& \leq C_0 \sum_{(i_3,i_4,i_7)\in\mathcal{W}(h)}\sum_{i_1,\cdots,i_8=0}^{\lfloor h/2 \rfloor} \left\|d^{1}_{i_1}d^{1+m_1}_{i_2}d^{\ell+h}_{i_5}d^{1+s}_{i_3}\right\|_{2+\nu}\left\|d^{\ell+h+m_1}_{i_6}d^{1+s+m_2}_{i_4}d^{1+s+h}_{i_7}d^{1+s+h+m_2}_{i_8}\right\|_{2+\nu} \notag\\ & \hspace{2cm} \times \alpha_{\Delta}^{\frac{\nu}{2+\nu}} \{\ell-s+\lfloor h/2 \rfloor +m_1\} \notag\\ & \leq C_0C \left(\sum_{i_1,\cdots,i_8=0}^{\infty} \rho^{\sum_{k=1}^{\infty} i_k}\right) \left(\sum_{(i_3,i_4,i_7)\in\mathcal{W}(h)}\alpha_{\Delta}^{\frac{\nu}{2+\nu}} \{\ell-s+m_1+\lfloor h/2\rfloor\}\right) .
\end{align}
Finally, we find that inequalities (\ref{eq:p1}) through (\ref{eq:p88}) are bounded by a constant independent of \(s\) and \(\ell\) under Assumption \((\mathbf{A_4})\), which proves that
\begin{align}
\label{eq:p8}
\sup_{s,\ell,e\in\mathbb{N}}\sum_{h=0}^{\infty}v_9(s,h,\ell,e,m_1,m_2) < \infty.
\end{align}
Consequently, in light of Equations (\ref{eq:333}) and (\ref{eq:p8}), we have
\begin{align}
\label{eq:p9}
\sup_{s,\ell,e\in\mathbb{N}}\sum_{h=0}^{\infty}\sum_{i_1,\cdots,i_8=0}^{\infty}\left|\text{Cov}(d_{i_1}^{1}d_{i_2}^{1+m_1}d_{i_3}^{1+s}d_{i_4}^{1+s+m_1},d_{i_5}^{\ell+h}d_{i_6}^{\ell+h+m_1}d_{i_7}^{1+e+h}d_{i_8}^{1+e+h+m_2})\right|<\infty.   
\end{align}
Secondly, we also have
\begin{align*}
\sum_{0\leq i_1,\ldots,i_8\leq \infty} K\rho^{\sum_{k=1}^{8}i_k}\left|\text{Cov}\left(\epsilon_{i_1,1,m_1,i_2}^{(2)}\epsilon_{i_3,1+s,m_2,i_4}^{(2)},\epsilon_{i_5,\ell+h,m_1,i_6}^{(2)}\epsilon_{i_7,1+e+h,m_2,i_8}^{(2)}\right)\right|\leq K\sum_{j=1}^{9}u_j(s,h,\ell,e,m_1,m_2)
\end{align*}
where 
\begin{align*}
    u_1(s,h,\ell,e,m_1,m_2) &:= \sum_{i_1 > \lfloor h/2 \rfloor} \sum_{0 \leq i_2, \cdots, i_8 \leq \infty} \rho^{\sum_{k=1}^{8}i_k}\left|\text{Cov}\left(\epsilon_{i_1,1,m_1,i_2}^{(2)}\epsilon_{i_3,1+s,m_2,i_4}^{(2)},\epsilon_{i_5,\ell+h,m_1,i_6}^{(2)}\epsilon_{i_7,1+e+h,m_2,i_8}^{(2)}\right)\right|, \\
    u_2(s,h,\ell,e,m_1,m_2) &:= \sum_{i_2 > \lfloor h/2 \rfloor} \sum_{0 \leq i_1, i_3, \cdots, i_8 \leq \infty} \rho^{\sum_{k=1}^{8}i_k}\left|\text{Cov}\left(\epsilon_{i_1,1,m_1,i_2}^{(2)}\epsilon_{i_3,1+s,m_2,i_4}^{(2)},\epsilon_{i_5,\ell+h,m_1,i_6}^{(2)}\epsilon_{i_7,1+e+h,m_2,i_8}^{(2)}\right)\right|, \\
    u_3(s,h,\ell,e,m_1,m_2) &:= \sum_{i_3 >\lfloor h/2 \rfloor} \sum_{0 \leq i_1, i_2, i_4, \cdots, i_8 \leq \infty} \rho^{\sum_{k=1}^{8}i_k}\left|\text{Cov}\left(\epsilon_{i_1,1,m_1,i_2}^{(2)}\epsilon_{i_3,1+s,m_2,i_4}^{(2)},\epsilon_{i_5,\ell+h,m_1,i_6}^{(2)}\epsilon_{i_7,1+e+h,m_2,i_8}^{(2)}\right)\right|,
\\
    u_4(s,h,\ell,e,m_1,m_2) &:= \sum_{i_4 >\lfloor h/2 \rfloor} \sum_{0 \leq i_1, i_2, i_3, i_5, \cdots, i_8 \leq \infty}\rho^{\sum_{k=1}^{8}i_k}\left|\text{Cov}\left(\epsilon_{i_1,1,m_1,i_2}^{(2)}\epsilon_{i_3,1+s,m_2,i_4}^{(2)},\epsilon_{i_5,\ell+h,m_1,i_6}^{(2)}\epsilon_{i_7,1+e+h,m_2,i_8}^{(2)}\right)\right|, \\
    u_5(s,h,\ell,e,m_1,m_2) &:= \sum_{i_5 >\lfloor h/2 \rfloor} \sum_{0 \leq i_1, \cdots, i_4, i_6, i_7, i_8 \leq \infty} \rho^{\sum_{k=1}^{8}i_k}\left|\text{Cov}\left(\epsilon_{i_1,1,m_1,i_2}^{(2)}\epsilon_{i_3,1+s,m_2,i_4}^{(2)},\epsilon_{i_5,\ell+h,m_1,i_6}^{(2)}\epsilon_{i_7,1+e+h,m_2,i_8}^{(2)}\right)\right|, \\
    u_6(s,h,\ell,e,m_1,m_2) &:= \sum_{i_6 > \lfloor h/2 \rfloor} \sum_{0 \leq i_1, \cdots, i_5, i_7, i_8 \leq \infty} \rho^{\sum_{k=1}^{8}i_k}\left|\text{Cov}\left(\epsilon_{i_1,1,m_1,i_2}^{(2)}\epsilon_{i_3,1+s,m_2,i_4}^{(2)},\epsilon_{i_5,\ell+h,m_1,i_6}^{(2)}\epsilon_{i_7,1+e+h,m_2,i_8}^{(2)}\right)\right|,\\
    u_7(s,h,\ell,e,m_1,m_2) &:= \sum_{i_7 > \lfloor h/2 \rfloor} \sum_{0 \leq i_1, \cdots, i_6, i_8 \leq \infty} \rho^{\sum_{k=1}^{8}i_k}\left|\text{Cov}\left(\epsilon_{i_1,1,m_1,i_2}^{(2)}\epsilon_{i_3,1+s,m_2,i_4}^{(2)},\epsilon_{i_5,\ell+h,m_1,i_6}^{(2)}\epsilon_{i_7,1+e+h,m_2,i_8}^{(2)}\right)\right|,\\
    u_8(s,h,\ell,e,m_1,m_2) &:= \sum_{i_8 >\lfloor h/2 \rfloor} \sum_{0 \leq i_1, \cdots, i_7 \leq \infty} \rho^{\sum_{k=1}^{8}i_k}\left|\text{Cov}\left(\epsilon_{i_1,1,m_1,i_2}^{(2)}\epsilon_{i_3,1+s,m_2,i_4}^{(2)},\epsilon_{i_5,\ell+h,m_1,i_6}^{(2)}\epsilon_{i_7,1+e+h,m_2,i_8}^{(2)}\right)\right|,\notag\\
    u_9(s,h,\ell,e,m_1,m_2) &:= \sum_{0 \leq i_1, \cdots, i_8 \leq \lfloor h/2 \rfloor} \rho^{\sum_{k=1}^{8}i_k}\left|\text{Cov}\left(\epsilon_{i_1,1,m_1,i_2}^{(2)}\epsilon_{i_3,1+s,m_2,i_4}^{(2)},\epsilon_{i_5,\ell+h,m_1,i_6}^{(2)}\epsilon_{i_7,1+e+h,m_2,i_8}^{(2)}\right)\right|.
\end{align*}
Using an argument similar to that in Equation (\ref{eq:14}), there also exists a family of positive constants \((\mathcal{D}_i)_{1\leq i\leq 8}\) such that
\begin{align}
\label{eq:3333}
\sup_{s,\ell,e\in\mathbb{N}}\displaystyle\sum_{h=0}^{\infty}u_i(s,h,\ell,e,m_1,m_2) \leq \displaystyle\sum_{h=0}^{\infty}\mathcal{D}_i\rho^{h/2} <\infty \ \forall \ 1\leq i\leq 8.
\end{align}
We now focus on the term \(u_9(s,h,\ell,e,m_1,m_2)\). Assume \(\lfloor h/2\rfloor \geq s + m_0\). By applying Lemma \ref{lem:5}, we have
\begin{align}
\label{eq:p10}
& \sum_{\lfloor h/2\rfloor=s+m_0}^{\infty}\sum_{i_1, \cdots, i_8 =0}^{\lfloor h/2 \rfloor} \rho^{\sum_{k=1}^{8}i_k}\left|\text{Cov}\left(\epsilon_{i_1,1,m_1,i_2}^{(2)}\epsilon_{i_3,1+s,m_2,i_4}^{(2)},\epsilon_{i_5,\ell+h,m_1,i_6}^{(2)}\epsilon_{i_7,1+e+h,m_2,i_8}^{(2)}\right)\right| \notag\\ & \leq K_0 \sum_{\lfloor h/2\rfloor=s+m_0}^{\infty} \sum_{i_1, \cdots, i_8 =0}^{\lfloor h/2 \rfloor}\rho^{\sum_{k=1}^{8}i_k}\left\|\eta_t\right\|^{8}_{8+4\nu}\alpha_{\eta}^{\frac{\nu}{2+\nu}}\{h-s-m_2-i_5+i_4\} \notag \\ & \leq K_0\left(\sum_{i_1, \cdots, i_8 =0}^{\infty}\rho^{\sum_{k=1}^{8}i_k}\right) \sum_{\lfloor h/2\rfloor=s+m_0}^{\infty} \left\|\eta_t\right\|^{8}_{8+4\nu}\alpha_{\eta}^{\frac{\nu}{2+\nu}}\{\lfloor h/2 \rfloor-s-m_0\}.
\end{align}
 As in (\ref{eq:p999}), to handle the case $\lfloor h/2\rfloor < s + m_0$, with $e>\max\{s,\ell\}$, we write
\begin{align}
\label{eq:fr72}
&\text{Cov}\left(\epsilon_{i_1,1,m_1,i_2}^{(2)}\epsilon_{i_3,1+s,m_2,i_4}^{(2)},\epsilon_{i_5,\ell+h,m_1,i_6}^{(2)}\epsilon_{i_7,1+e+h,m_2,i_8}^{(2)}\right) \notag\\ & = \text{Cov}\left(\epsilon_{i_1,1,m_1,i_2}^{(2)}\epsilon_{i_5,\ell+h,m_1,i_6}^{(2)},\epsilon_{i_3,1+s,m_2,i_4}^{(2)}\epsilon_{i_7,1+e+h,m_2,i_8}^{(2)}\right) \notag\\ & + \mathbb{E}\left\{\epsilon_{i_1,1,m_1,i_2}^{(2)}\epsilon_{i_5,\ell+h,m_1,i_6}^{(2)}\right\}\mathbb{E}\left\{\epsilon_{i_3,1+s,m_2,i_4}^{(2)}\epsilon_{i_7,1+e+h,m_2,i_8}^{(2)}\right\}\notag\\ & - \mathbb{E}\left\{\epsilon_{i_1,1,m_1,i_2}^{(2)}\epsilon_{i_3,1+s,m_2,i_4}^{(2)}\right\}\mathbb{E}\left\{\epsilon_{i_5,\ell+h,m_1,i_6}^{(2)}\epsilon_{i_7,1+e+h,m_2,i_8}^{(2)}\right\}.
\end{align}
 Using Lemma \ref{lem:5} once again, and the fact that $(\eta_t)_{t \in \mathbb{Z}}$ is weak white noise, we obtain :
\subparagraph{$\diamond$ First term in Equation (\ref{eq:fr72})}
\begin{align*}
& \sum_{\substack{\lfloor h/2 \rfloor =0 \\ \lfloor h/2 \rfloor \leq s-\ell-m_1}}^{s+m_0-1} \sum_{i_1, \cdots, i_8 =0}^{\lfloor h/2 \rfloor} \rho^{\sum_{k=1}^{8}i_k} \left|\text{Cov}\left(\epsilon_{i_1,1,m_1,i_2}^{(2)} \epsilon_{i_5,\ell+h,m_1,i_6}^{(2)}, \epsilon_{i_3,1+s,m_2,i_4}^{(2)} \epsilon_{i_7,1+e+h,m_2,i_8}^{(2)}\right)\right| \notag\\
& \leq K_0 \sum_{\substack{\lfloor h/2 \rfloor=0 \\ \lfloor h/2 \rfloor \leq s-\ell-m_1}}^{s+m_0-1} \sum_{i_1, \cdots, i_8 =0}^{\lfloor h/2 \rfloor} \rho^{\sum_{k=1}^{8}i_k} \left\|\eta_t\right\|^{8}_{8+4\nu} \alpha_{\eta}^{\frac{\nu}{2+\nu}}\{1+s-\ell-h-m_1-i_3+i_6\} \notag
\end{align*}
\begin{align}
\label{eq:p11}
& \leq K_0 \left(\sum_{i_1, \cdots, i_8 =0}^{\infty} \rho^{\sum_{k=1}^{8}i_k}\right) \sum_{\substack{\lfloor h/2 \rfloor=0 \\ \lfloor h/2 \rfloor \leq s-\ell-m_1}}^{s+m_0-1} \left\|\eta_t\right\|^{8}_{8+4\nu} \alpha_{\eta}^{\frac{\nu}{2+\nu}} \{1+s-\lfloor 3h/2 \rfloor - \ell - m_1\} \notag\\ 
& \leq K_0 \left(\sum_{i_1, \cdots, i_8 =0}^{\infty} \rho^{\sum_{k=1}^{8}i_k}\right)\left\|\eta_t\right\|^{8}_{8+4\nu} \left(\sum_{\lfloor h/2 \rfloor=0}^{s-\ell-m_1} \alpha_{\eta}^{\frac{\nu}{2+\nu}}\{s-\lfloor 3h/2 \rfloor - \ell - m_1\}\right).
\end{align}
\subparagraph{$\diamond$ Second term in Equation (\ref{eq:fr72})}
\begin{align}
\label{eq:p12}
& \sum_{\lfloor h/2 \rfloor=0}^{s+m_0-1}\sum_{i_1, \cdots, i_8 =0}^{\lfloor h/2 \rfloor}\rho^{\sum_{k=1}^{8}i_k}\mathbb{E}\left\{\epsilon_{i_1,1,m_1,i_2}^{(2)}\epsilon_{i_5,\ell+h,m_1,i_6}^{(2)}\right\} \notag\\ = & \sum_{\lfloor h/2 \rfloor=0}^{s+m_0-1}\sum_{i_1, \cdots, i_8 =0}^{\lfloor h/2 \rfloor}\rho^{\sum_{k=1}^{8}i_k} \Bigg\{\text{Cov}\left(\epsilon_{i_1,1,m_1,i_2}^{(2)},\epsilon_{i_5,\ell+h,m_1,i_2}^{(2)}\right)+\mathbb{E}\left(\epsilon_{i_1,1,m_1,i_2}^{(2)}\right)\mathbb{E}\left(\epsilon_{i_5,\ell+h,m_1,i_6}^{(2)}\right)\Bigg\}\notag\\  \leq & \sum_{\lfloor h/2 \rfloor =0}^{s+m_0-1}\sum_{i_1, \cdots, i_8 =0}^{\lfloor h/2 \rfloor}\rho^{\sum_{k=1}^{8}i_k}\Bigg\{K_0 \|\eta_t\|_{4+2\nu}^4\alpha_{\eta}^{\frac{\nu}{2+\nu}}\{\ell+h-m_1-i_5+i_2-1\} \notag\\ & \hspace{3cm} + K_0^2 \|\eta_{t}\|_{2+\nu}^{4}\alpha_{\eta}^{\frac{\nu}{2+\nu}}\{1+m_1-i_2+i_1\}\alpha_{\eta}^{\frac{\nu}{2+\nu}}\{m_1+i_6-i_5\}\Bigg\}\notag\\  \leq &\left(\sum_{i_1, \cdots, i_8 =0}^{\infty}\rho^{\sum_{k=1}^{8}i_k}\right)\Bigg\{\sum_{\lfloor h/2 \rfloor=0}^{s+m_1-1} K_0\|\eta_t|_{4+2\nu}^4\alpha_{\eta}^{\frac{\nu}{2+\nu}}\{\lfloor h/2 \rfloor-m_1-1\}\notag\\ & \hspace{3cm} + K_0^2\|\eta_{t}\|_{2+\nu}^{4}\alpha_{\eta}^{\frac{\nu}{2+\nu}}\{m_1-\lfloor h/2 \rfloor\}\alpha_{\eta}^{\frac{\nu}{2+\nu}}\{m_1-\lfloor h/2 \rfloor\}\Bigg\}.
\end{align}
\subparagraph{$\diamond$ Third term in Equation (\ref{eq:fr72})}
\begin{align}
\label{eq:p13}
& \sum_{\lfloor h/2 \rfloor=0}^{s+m_0-1}\sum_{i_1, \cdots, i_8 =0}^{\lfloor h/2 \rfloor}\rho^{\sum_{k=1}^{8}i_k}\mathbb{E}\left\{\epsilon_{i_1,1,m_1,i_2}^{(2)}\epsilon_{i_3,1+s,m_2,i_4}^{(2)}\right\} \notag\\ = & \sum_{\lfloor h/2 \rfloor=0}^{s+m_0-1}\sum_{i_1, \cdots, i_8 =0}^{\lfloor h/2 \rfloor}\rho^{\sum_{k=1}^{8}i_k}\Bigg\{\text{Cov}\left(\epsilon_{i_1,1,m_1,i_2}^{(2)},\epsilon_{i_3,1+s,m_2,i_4}^{(2)}\right)+\mathbb{E}\left(\epsilon_{i_1,1,m_1,i_2}^{(2)}\right)\mathbb{E}\left(\epsilon_{i_3,1+s,m_2,i_4}^{(2)}\right)\Bigg\}\notag\\  \leq & \sum_{\lfloor h/2 \rfloor =0}^{s+m_0-1}\sum_{i_1, \cdots, i_8 =0}^{\lfloor h/2 \rfloor}\rho^{\sum_{k=1}^{8}i_k}\Bigg\{K_0\|\eta_t\|_{4+2\nu}^4\alpha_{\eta}^{\frac{\nu}{2+\nu}}\{1+s-\max\{i_3,m_2-i_4\}\} \notag\\ & \hspace{3cm} + K_0^2 \|\eta_{t}\|_{2+\nu}^{4}\alpha_{\eta}^{\frac{\nu}{2+\nu}}\{m_1-i_2+i_1\}\alpha_{\eta}^{\frac{\nu}{2+\nu}}\{m_2-i_4+i_3\}\Bigg\}\notag\\ \leq & \left(\sum_{i_1, \cdots, i_8 =0}^{\infty}\rho^{\sum_{k=1}^{8}i_k}\right)\Bigg\{\sum_{\lfloor h/2 \rfloor=0}^{s+m_0-1} K_0\|\eta_t|_{4+2\nu}^4\alpha_{\eta}^{\frac{\nu}{2+\nu}}\{s-m_0-\lfloor h/2 \rfloor\}\notag\\ & \hspace{3cm} + K_0^2 \|\eta_{t}\|_{2+\nu}^{4}\alpha_{\eta}^{\frac{\nu}{2+\nu}}\{m_1-\lfloor h/2 \rfloor\}\alpha_{\eta}^{\frac{\nu}{2+\nu}}\{m_2-\lfloor h/2 \rfloor\}\Bigg\}.
\end{align}
 Ultimately, we conclude that the inequalities (\ref{eq:p10}), (\ref{eq:p12}), and (\ref{eq:p13}) are also bounded by a constant independent of \(s, e\), and \(\ell\). Similarly to the sixth term of (\ref{eq:p999}), we also bound (\ref{eq:p11}) by a constant independent of \(s, e\), and \(\ell\). Furthermore, using a decomposition similar to (\ref{eq:fr72}), we also show that 
\begin{align*}
\sum_{\substack{\lfloor h/2 \rfloor =0 \\ \lfloor h/2 \rfloor \geq s-\ell-m_1+1}}^{s+m_0-1} \sum_{i_1, \cdots, i_8 =0}^{\lfloor h/2 \rfloor} \rho^{\sum_{k=1}^{8}i_k} \left|\text{Cov}\left(\epsilon_{i_1,1,m_1,i_2}^{(2)} \epsilon_{i_5,\ell+h,m_1,i_6}^{(2)}, \epsilon_{i_3,1+s,m_2,i_4}^{(2)} \epsilon_{i_7,1+e+h,m_2,i_8}^{(2)}\right)\right|
\end{align*}
is bounded by a constant independent of \(s,e\) and \(\ell\). Consequently, we conclude that
\begin{align}
\label{eq:p14}
   \sup_{s,\ell,e\in\mathbb{N}} \sum_{h=0}^{\infty} u_9(s,h,\ell,e,m_1,m_2) < \infty.
\end{align}

\noindent Thus, by combining Equations (\ref{eq:3333}) and (\ref{eq:p14}), we arrive at
\begin{align}
\label{eq:p15}
\sup_{s,\ell,e\in\mathbb{N}}\sum_{h=0}^{\infty}\sum_{i_1,\cdots,i_8=0}^{\infty}\rho^{\sum_{k=1}^{8}i_k}\left|\text{Cov}\left(\epsilon_{i_1,1,m_1,i_2}^{(2)}\epsilon_{i_3,1+s,m_2,i_4}^{(2)},\epsilon_{i_5,\ell+h,m_1,i_6}^{(2)}\epsilon_{i_7,1+e+h,m_2,i_8}^{(2)}\right)\right|<\infty.
\end{align}

\medskip

\noindent \noindent Finally, in view of Equations (\ref{eq:p9}) and (\ref{eq:p15}), we reach the conclusion that
\begin{align}
\label{eq:344}
\sup_{s,\ell,e\in\mathbb{N}}\displaystyle\sum_{h=0}^{\infty}\Big|\text{Cov}(X_1X_{1+m_1}X_{1+s}X_{1+s+m_2},X_{\ell+h}X_{\ell+h+m_1}X_{1+e+h}X_{1+e+h+m_2}) \Big|< \infty.
\end{align}

 Furthermore, by applying H{\"o}lder inequality once again, it follows from Lemma \ref{lem:111} that
\begin{align*}
\left|\mathbb{E}\left(X_{1+e+h}X_{1+e+h+m_2}\right)\right| \leq & \left\|X_{1+e+h}\right\|_2\left\|X_{1+e+h+m_2}\right\|_2\leq \left\|X_0\right\|_2^2<\infty,
\end{align*}
 and 
\begin{align*}
\left|\mathbb{E}(X_{\ell+h}X_{\ell+h+m_1})\right|\leq \left\|X_{\ell+h}\right\|_2\left\|X_{\ell+h+m_1}\right\|_2\leq \|X_0\|_2^2<\infty.    
\end{align*}
 By  analogous arguments used to obtain Equation (\ref{eq:344}), we can readily demonstrate that
\begin{align*}
\sup_{\ell,s\in\mathbb{N}}\sum_{h=0}^{\infty}\Big|\text{Cov}(X_{1}X_{1+m_1}X_{1+s+m_2},X_{\ell+h}X_{\ell+h+m_1})\Big| < \infty,
\end{align*}
 and
\begin{align}
 \sup_{e,s\in\mathbb{N}} \sum_{h=0}^\infty\left|\text{Cov}(X_{1}X_{1+m_1}X_{1+s}X_{1+s+m_2},X_{1+e+h}X_{1+e+h+m_2})\right|<\infty.
\end{align}
 Consequently we have
\begin{align}
\label{eq:355}
\sup_{s,\ell,e\in\mathbb{N}}\sum_{h=0}^{\infty} \left| \mathbb{E}(X_{1+e+h}X_{1+e+h+m_2})\text{Cov}(X_{1}X_{1+m_1}X_{1+s+m_2},X_{\ell+h}X_{\ell+h+m_1})\right| < \infty,
\end{align}
and
\begin{align}
\label{eq:366}
\sup_{s,\ell,e\in\mathbb{N}}\sum_{h=0}^{\infty} \left|\mathbb{E}(X_{\ell+h}X_{\ell+h+m_1})\text{Cov}(X_{1}X_{1+m_1}X_{1+s}X_{1+s+m_2},X_{1+e+h}X_{1+e+h+m_2})\right| < \infty.
\end{align}
 Ultimately, by combining Equations (\ref{eq:344}), (\ref{eq:355}) and (\ref{eq:366}), we draw the conclusion that 
\begin{align}
\label{eq:3771}
\sup_{s,\ell,e\in\mathbb{N}}\sum_{h=0}^{\infty}|w_1(s,h,\ell,e,m_1,m_2)|< \infty.
\end{align}
Proceeding similarly as in Equation (\ref{eq:3771}), we also show that
\begin{align*}
\sup_{s,\ell,e\in\mathbb{N}}\sum_{h=0}^{\infty}|w_2(s,h,\ell,e,m_1,m_2)|< \infty,    
\ \text{and} \
\sup_{s,\ell,e\in\mathbb{N}}\sum_{h=0}^{\infty}|w_3(s,h,\ell,e,m_1,m_2)|< \infty.   
\end{align*}
 The same bounds clearly hold for \(h \leq 0\). Thus, we have demonstrated that
\begin{align*}
\sup_{s,\ell,e\in\mathbb{N}}\sum_{h=-\infty}^{\infty}\left|\text{Cov}\left(\mathcal{Y}_1(m_1)\mathcal{Y}_{1+s}(m_2), \mathcal{Y}_{\ell+h}(m_1)\mathcal{Y}_{1+e+h}(m_2)\right)\right| < \infty.
\end{align*}
This completes the proof.  
\qed

\begin{l1}\label{lem:12}
We assume that the condition $\mathbb{E}|\eta_t|^{8+4\nu}<\infty$ holds and that Assumption $(\mathbf{A_4})$ is satisfied for some $\nu>0.$ Then for any integer $s$, there exists a positive constant $C$ independent of $n$ such that
\begin{align*}
    \sum_{h=-\infty}^{\infty}\left| \Cov\left(\hat{\mathcal{Y}}_1(m_1)\hat{\mathcal{Y}}_{1+s}(m_2), \hat{\mathcal{Y}}_{1+h}(m_1)\hat{\mathcal{Y}}_{1+s+h}(m_2)\right) \right| < C,
\end{align*}
where $m_1, m_2 \in \{1, \ldots, N\}$ and $\hat{\mathcal{Y}}_t(k)$ denotes the $k$-th element of the vector $\hat{\mathcal{Y}}_t$. We also recall that $\hat{\mathcal{Y}}_t(m) = 0$ for all $m \in \{1, \ldots, N\}$ when $t > n$ or $t \leq 0$.
\end{l1}
\subsubsection*{Proof.}

 The proof of this lemma will closely follow the approach of Lemma \ref{lem:11}.

 Let \(h \in \mathbb{Z}\) and \(1 \leq m_1, m_2 \leq N.\) We have
\begin{align*}
\displaystyle\sum_{h=-\infty}^{\infty}\left|\text{Cov}\left(\hat{\mathcal{Y}}_1(m_1)\hat{\mathcal{Y}}_{1+s}(m_2), \hat{\mathcal{Y}}_{1+h}(m_1)\hat{\mathcal{Y}}_{1+s+h}(m_2)\right)\right| \leq \displaystyle\sum_{h=0}^{n-s-1} \sum_{j=1}^{4}g_j(s,h,m_1,m_2)\mathds{1}_{[0\leq s \leq n-1]}
\end{align*}
 with

\begin{align*}
g_{1}(s,h,m_1,m_2): = &  \text{Cov}(X_1X_{1+m_1}X_{1+s}X_{1+s+m_2},X_{1+h}X_{1+h+m_1}X_{1+s+h}X_{1+s+h+m_2}) \\ & -\text{Cov}(X_1X_{1+m_1}X_{1+s}X_{1+s+m_2},\hat{c}_{m_1,0}X_{1+s+h}X_{1+s+h+m_2}) \\ & -
\text{Cov}(X_1X_{1+m_1}X_{1+s}X_{1+s+m_2},\hat{c}_{m_2,0}X_{1+h}X_{1+h+m_1}) \\ &  + \text{Cov}(X_{1}X_{1+m_1}X_{1+s}X_{1+s+m_2},\hat{c}_{m_1,0}\hat{c}_{m_2,0})
,\\
g_2(s,h,m_1,m_2): = & - \text{Cov}(\hat{c}_{m_2,0}X_{1}X_{1+m_1},X_{1+h}X_{1+h+m_1}X_{1+h+m_1}X_{1+s+h}X_{1+s+h+m_2}) \\ & 
+ \text{Cov}(\hat{c}_{m_2,0}X_1X_{1+m_1},\hat{c}_{m_1,0}X_{1+s+h}X_{1+s+h+m_2})\\ & +
\text{Cov}(\hat{c}_{m_2,0}X_{1}X_{1+m_1},\hat{c}_{m_2,0}X_{1+h}X_{1+h+m_1})\\ & -\text{Cov}(\hat{c}_{m_2,0}X_{1}X_{1+m_1},\hat{c}_{m_1,0}\hat{c}_{m_2,0})
,\\
g_{3}(s,h,m_1,m_2): = & -\text{Cov}(\hat{c}_{m_1,0}X_{1+s}X_{1+s+m_2}, X_{1+h}X_{1+h+m_1}X_{1+s+h}X_{1+s+h+m_2})\\ & +\text{Cov}(\hat{c}_{m_1,0}X_{1+s}X_{1+s+m_2},\hat{c}_{m_2,0}X_{1+h}X_{1+h+m_1})\\ & + \text{Cov}(\hat{c}_{m_1,0}X_{1+s}X_{1+s+m_2},\hat{c}_{m_2,0}X_{1+s+h}X_{1+s+h+m_2})\\ & 
-\text{Cov}(\hat{c}_{m_1,0}X_{1+s}X_{1+s+m_2},\hat{c}_{m_1,0}\hat{c}_{m_2,0}),
\\
g_{4}(s,h,m_1,m_2):  = &   \text{Cov}(\hat{c}_{m_1,0}\hat{c}_{m_2,0},X_{1+h}X_{1+h+m_1}X_{1+s+h}X_{1+s+h+m_2})\\ & - \text{Cov}(\hat{c}_{m_1,0}\hat{c}_{m_2,0},\hat{c}_{m_1,0}X_{1+s+h}X_{1+s+h+m_2})\\ & -\text{Cov}(\hat{c}_{m_1,0}\hat{c}_{m_2,0},\hat{c}_{m_2,0}X_{1+h}X_{1+h+m_1})\\ & 
+\text{Cov}(\hat{c}_{m_1,0}\hat{c}_{m_2,0},\hat{c}_{m_1,0}\hat{c}_{m_2,0}),
\end{align*}
 where we recall that $\hat{c}_{k,0} = (n-k)^{-1}\sum_{t=1}^{n-k} X_tX_{t+k}$ for all $k\in\mathbb{N}^\star.$


  Based on Equation (\ref{eq:355}) of Lemma~\ref{lem:11}, there exists a positive constant $C$ independent of $t$ and $s$ such that
\begin{align}
\label{eq:4044}
\sum_{h=-\infty}^{\infty}\left|\text{Cov}(X_1X_{1+m_1}X_{1+s}X_{1+s+m_2},X_{t+h}X_{t+h+m_1}X_{1+s+h}X_{1+s+h+m_2}) \right| < C.
\end{align}
Thanks to Equation~(\ref{eq:4044}), we easily obtain that
\begin{align}
\label{eq:eq87}
&\sum_{h=0}^{n-s-1}\left|\text{Cov}\left(X_{1}X_{1+m_1}X_{1+s}X_{1+s+m_2},\hat{c}_{m_1,0}X_{1+s+h}X_{1+s+h+m_2}\right)\right|\notag\\&\leq \sum_{h=0}^{\infty}\left|\text{Cov}\left(X_{1}X_{1+m_1}X_{1+s}X_{1+s+m_2},\hat{c}_{m_1,0}X_{1+s+h}X_{1+s+h+m_2}\right)\right| \notag\\ 
&\leq \sum_{h=0}^{\infty}\left\{(n-m_1)^{-1}\sum_{t=1}^{n-m_1}\left|\text{Cov}\left(X_{1}X_{1+m_1}X_{1+s}X_{1+s+m_2},X_{t+h}X_{t+h+m_1}X_{1+s+h}X_{1+s+h+m_2}\right)\right|\right\}\notag\\
&= (n-m_1)^{-1}\sum_{t=1}^{n-m_1}\left\{\sum_{h=0}^{\infty}\left|\text{Cov}\left(X_{1}X_{1+m_1}X_{1+s}X_{1+s+m_2},X_{t+h}X_{t+h+m_1}X_{1+s+h}X_{1+s+h+m_2}\right)\right|\right\} < C,
\end{align}
Similarly, there exist positive constants, each denoted by \( C \) and independent of \( n \), such that 
\begin{align*}
 \sum_{h=0}^{n-s-1}\left|\text{Cov}(X_1X_{1+m_1}X_{1+s}X_{1+s+m_2},\hat{c}_{m_2,0}X_{1+h}X_{1+h+m_1})\right|< C.
\end{align*}
Let us now focus on the fourth term of $g_4(s,h,m_1,m_2)$, namely $\text{Cov}(\hat{c}_{m_1,0}\hat{c}_{m_2,0}, \hat{c}_{m_1,0}\hat{c}_{m_2,0})$, which represents the most challenging term to bound.

 By the stationarity of $(X_t)_{t \in \mathbb{Z}}$, we have
\begin{align}
\label{eq:l9}
 &\left|\Cov(\hat{c}_{m_1,0}\hat{c}_{m_2,0}, \hat{c}_{m_1,0}\hat{c}_{m_2,0})\right| \notag\\ & \leq  n^{-4}\left\{\left|\sum_{i_1,i_2=1}^{n}\sum_{i_3,i_4=1}^{n}\Cov\left(X_{i_1}X_{i_1+m_1}X_{i_2}X_{i_2+m_2},X_{i_3}X_{i_3+m_1}X_{i_4}X_{i_4+m_2}\right)\right|\right\}\notag\\ & \leq n^{-4}\sum_{h_1,h_2=1}^{n-1}\left\{\sum_{i_2,i_4=1}^{n-h_1\wedge h_2}\left|\Cov\left(X_{i_2+h_1}X_{i_2+h_1+m_1}X_{i_2}X_{i_2+m_2},X_{i_4+h_2}X_{i_4+h_2+m_1}X_{i_4}X_{i_4+m_2}\right)\right|\right\}\notag\\ &\leq n^{-3}\sum_{h_1,h_2=1}^{n-1}\left\{\sum_{k=-n+h_1\wedge h_2+1}^{n-h_1\wedge h_2-1}\frac{n-h_1\wedge h_2-|k|}{n}\left|\Cov\left(Y_{1+h_1}(m_1)Y_{1}(m_2),Y_{1+h_2-k}(m_1)Y_{1-k}(m_2)\right)\right|\right\}\notag\\ & \leq n^{-3}\sum_{h_1,h_2=1}^{n-1} \left\{\sum_{k= -\infty}^{\infty}\left|\Cov\left(Y_{1+h_1}(m_1)Y_{1}(m_2),Y_{1+h_2-k}(m_1)Y_{1-k}(m_2)\right)\right|\right\},
\end{align}
where we recall that for all $t \in \mathbb{Z}$, $Y_t(m) := X_t X_{t+m}$ for a given $m \in \{1, \dots, N\}$.

 Using Lemma \ref{lem:11}, we deduce that there exists a positive constant $\mathcal{K}$, independent of $n$, $h_1$ and $h_2$ such that
\begin{align}
\label{eq:l10}
 \sum_{k= -\infty}^{\infty}\left|\Cov\left(Y_{1+h_1}(m_1)Y_{1}(m_2),Y_{1+h_2-k}(m_1)Y_{1-k}(m_2)\right)\right| < \mathcal{K}. 
\end{align}
Thus, from Equations (\ref{eq:l9}) and (\ref{eq:l10}), it follows that
\begin{align}
\label{eq:l11}
    \sum_{h=0}^{n-s-1} \left|\text{Cov}(\hat{c}_{m_1,0}\hat{c}_{m_2,0}, \hat{c}_{m_1,0}\hat{c}_{m_2,0})\right| \leq  (n-s)(n-1)^2n^{-3} \mathcal{K}=\mathrm{O}(1).
\end{align}
 By reasoning similarly to Equation (\ref{eq:eq87}) or Equation (\ref{eq:l11}), the remaining terms involving $g_2$ and $g_3$ can easily be bounded by a constant independent of $n$ and $s$. Consequently, we arrive at
\begin{align*}
  \sum_{h=-\infty}^{\infty} |g_i(s,h,m_1,m_2)| < C, \quad \forall \ 1 \leq i \leq 4,
\end{align*}
where $C$ is an arbitrary constant independent of $n$.

This completes the proof.
\qed
\begin{l1}\label{lem:13}
We assume that the condition $\mathbb{E}|\eta_t|^{8+4\nu}<\infty$ holds and that Assumption $(\mathbf{A_4})$ is satisfied for some $\nu>0.$ Then for any integer $s$, there exists a positive constant $C$ such that
\begin{align*}
    \sum_{h=-\infty}^{\infty}\left| \text{Cov}\left(\mathcal{Y}_1(m_1)\mathcal{Y}_{1+s}(m_2), \hat{\mathcal{Y}}_{1+h}(m_1)\hat{\mathcal{Y}}_{1+s+h}(m_2)\right) \right| < C,
\end{align*}
where $m_1, m_2 \in \{1, \ldots, N\}$ and $\hat{\mathcal{Y}}_t(k)$ denotes the $k$-th element of the vector $\hat{\mathcal{Y}}_t$ with $\hat{\mathcal{Y}}_t(m) = 0$ for all $m \in \{1, \ldots, N\}$ whenever $t > n$ or $t \leq 0$.
\end{l1}

\subsubsection*{Proof.}
This proof follows a similar approach to that of Lemmas~\ref{lem:11} and ~\ref{lem:12}.  It involves expanding the covariance and appropriately bounding each term, as previously established and it is omitted.
\qed

\begin{l1}\label{lem:14}
 Under the assumptions of Theorem~\ref{thm:4}, the terms \(\sqrt{r}\|\hat{\Sigma}_{\underline{\mathcal{Y}}_r}-\Sigma_{\underline{\mathcal{Y}}_r}\|\), \(\sqrt{r}\|\hat{\Sigma}_{\mathcal{Y}}-\Sigma_{\mathcal{Y}}\|\) and \(\sqrt{r}\|\hat{\Sigma}_{\mathcal{Y},\underline{\mathcal{Y}}_r}-\Sigma_{\mathcal{Y},\underline{\mathcal{Y}}_r}\|\) tend towards 0 in probability as \(n \rightarrow \infty\) when \(r = \mathrm{o}(n^{1/3})\).
\end{l1}

\subsubsection*{Proof.} To demonstrate this lemma, we will focus solely on the proof of one term, while the other terms can be demonstrated similarly. 

 Using Markov inequality, we have
\begin{align}
\label{eq:42}
    \forall \ \epsilon > 0, \ \mathbb{P}\left(\sqrt{r}\|\hat{\Sigma}_{\underline{\mathcal{Y}}_r}-\Sigma_{\underline{\mathcal{Y}}_r}\| > \epsilon\right) \leq \frac{1}{\epsilon^2} \mathbb{E}\left(r\|\hat{\Sigma}_{\underline{\mathcal{Y}}_r}-\Sigma_{\underline{\mathcal{Y}}_r}\|^2\right).
\end{align} 
Let \(1 \leq m_1, m_2 \leq N\) and \(1 \leq r_1, r_2 \leq r\). The element of the $\{r_1N+m_1\}$ th-row and $\{r_2N+m_1\}$ th-column of $\left(\hat{\Sigma}_{\underline{\mathcal{Y}}_r}-\Sigma_{\underline{\mathcal{Y}}_r}\right)$ is of the form $n^{-1}\sum_{t=1}^{n}\mathcal{Y}_{t-r_1}(m_1)\mathcal{Y}_{t-r_2}(m_2) - \mathbb{E}\left(\mathcal{Y}_{t-r_1}(m_1)\mathcal{Y}_{t-r_2}(m_2)\right)$.

 Define \(\mathcal{X}_t(r_1,r_2,m_1,m_2) : = \mathcal{Y}_{t-r_1}(m_1)\mathcal{Y}_{t-r_2}(m_2)\). Using Equation (\ref{eq:505}),  Lemma~\ref{lem:11} and the stationarity of \((\mathcal{Y}_t)_{t\in\mathbb{Z}}\), it follows that
\begin{align}
\label{eq:43}
    \mathbb{E}\left(r\|\hat{\Sigma}_{\underline{\mathcal{Y}}_r}-\Sigma_{\underline{\mathcal{Y}}_r}\|^{2}\right) \leq & \sum_{m_1,m_2=1}^{N} \sum_{r_1,r_2=1}^{r} \mathbb{E}\left(r\left(    \left(\hat{\Sigma}_{\underline{\mathcal{Y}}_r}-\Sigma_{\underline{\mathcal{Y}}_r}\right)(r_1N+m_1, r_2N+m_2)\right)^2\right) \notag\\
     \leq  & \sum_{m_1,m_2=1}^{N} \sum_{r_1,r_2=1}^{r}\left\{\frac{1}{n^2}\mathbb{E}\left(r\left(\sum_{t=1}^n \left(\mathcal{Y}_{t-r_1}(m_1)\mathcal{Y}_{t-r_2}(m_2)-\mathbb{E}(\mathcal{Y}_{t-r_1}(m_1)\mathcal{Y}_{t-r_2}(m_2)\right)\right)^2\right)\right\} \notag \\
     \leq  & \sum_{m_1,m_2=1}^{N} \sum_{r_1,r_2=1}^{r}
     r \ \text{Var}\left(\frac{1}{n}\sum_{t=1}^n \left(\mathcal{Y}_{t-r_1}(m_1)\mathcal{Y}_{t-r_2}(m_2)\right)\right)\notag\\
      \leq  & \sum_{m_1,m_2=1}^{N} \sum_{r_1,r_2=1}^{r}
      \frac{r}{n^2}\sum_{k=-n+1}^{n-1}(n-|k|)\text{Cov}(\mathcal{X}_t(r_1,r_2,m_1,m_2), \mathcal{X}_{t-k}(r_1,r_2,m_1,m_2)) \notag\\
      \leq  & \sum_{m_1,m_2=1}^{N} \sum_{r_1,r_2=1}^{r}
      \frac{r}{n}\sum_{k=-\infty}^{\infty}\text{Cov}(\mathcal{X}_t(r_1,r_2,m_1,m_2), \mathcal{X}_{t-k}(r_1,r_2,m_1,m_2)) \notag\\
      \leq & \ \frac{\mathcal{C}_{12}N^2r^3}{n}\xrightarrow[n\rightarrow\infty]{r \ = \ \mathrm{o}(n^{1/3})}0, \notag\\
\end{align}
where \(\mathcal{C}_{12}\) is here a positive constant independent of \(r_1, r_2, m_1, m_2, r\), and \(n\). 

 Consequently using Equations (\ref{eq:42}) and (\ref{eq:43}), when \( r = \mathrm{o}(n^{1/3}) \), we finally obtain that 
\begin{align}
\label{eq:501}
    \|\hat{\Sigma}_{\underline{\mathcal{Y}}_r}-\Sigma_{\underline{\mathcal{Y}}_r}\| = \mathrm{o}_{\mathbb{P}}(1).
\end{align}

\noindent Similarly, when $r = \mathrm{o}(n^{1/3})$, we also prove that 
\begin{align*}
  \sqrt{r}\|\hat{\Sigma}_{\mathcal{Y}}-\Sigma_{\mathcal{Y}}\| = \mathrm{o}_{\mathbb{P}}(1) \ \text{and} \ \sqrt{r}\|\hat{\Sigma}_{\hat{\mathcal{Y}},\underline{\hat{\mathcal{Y}}}_r}-\Sigma_{\mathcal{Y},\underline{\mathcal{Y}}_r}\| = \mathrm{o}_{\mathbb{P}}(1).
\end{align*}
This completes the proof.
\qed
\begin{l1}\label{lem:15}
Let \(\hat{\Sigma}_{\mathcal{Y}}\) be the matrix obtained by replacing \(\hat{\mathcal{Y}}_t\) with \(\mathcal{Y}_t\) in \(\hat{\Sigma}_{\hat{\mathcal{Y}}}\). Under the assumptions of Theorem~\ref{thm:4}, the terms \(\sqrt{r}\|\hat{\Sigma}_{\underline{\hat{\mathcal{Y}}}_r}-\Sigma_{\underline{\mathcal{Y}}_r}\|\), \(\sqrt{r}\|\hat{\Sigma}_{\hat{\mathcal{Y}}}-\Sigma_{\mathcal{Y}}\|\) and \(\sqrt{r}\|\hat{\Sigma}_{\hat{\mathcal{Y}},\underline{\hat{\mathcal{Y}}}_r}-\Sigma_{\mathcal{Y},\underline{\mathcal{Y}}_r}\|\) tend towards 0 in probability as \(n \rightarrow \infty\) when \(r = \mathrm{o}(n^{1/3})\).
\end{l1}
\subsubsection*{Proof.}
Let \( r = r(n) \) be such that \( r = \mathrm{o}(n^{1/3}) \). Applying the Markov inequality once again, we have
\begin{align*}
   \forall \ \varepsilon > 0, \ \mathbb{P}\left(\sqrt{r}\|\hat{\Sigma}_{\underline{\hat{\mathcal{Y}}}_r} - \hat{\Sigma}_{\underline{\mathcal{Y}}_r}\| > \varepsilon\right) \leq \frac{1}{\varepsilon^2} \mathbb{E}\left(r\|\hat{\Sigma}_{\underline{\hat{\mathcal{Y}}}_r} - \hat{\Sigma}_{\underline{\mathcal{Y}}_r}\|^2\right).
\end{align*}
 Let $1 \leq m_1, m_2 \leq N$ and $1 \leq r_1, r_2 \leq r$. The component of $\hat{\Sigma}_{\underline{\hat{\mathcal{Y}}}_r} - \hat{\Sigma}_{\underline{\mathcal{Y}}_r}$ located in the $\{r_1N+m_1\}$-th row and $\{r_2N+m_2\}$-th column is of the form $n^{-1}\sum_{t=1}^{n}\hat{\mathcal{Y}}_{t-r_1}(m_1)\hat{\mathcal{Y}}_{t-r_2}(m_2) - \mathcal{Y}_{t-r_1}(m_1)\mathcal{Y}_{t-r_2}(m_2)$.

 Letting $\mathcal{Z}_t = \mathcal{Z}_t(r_1,r_2,m_1,m_2,n): = \hat{\mathcal{Y}}_{t-r_1}(m_1)\hat{\mathcal{Y}}_{t-r_2}(m_2) - \mathcal{Y}_{t-r_1}(m_1)\mathcal{Y}_{t-r_2}(m_2)$ and using once again the norm defined in Equation (\ref{eq:505}), we have
\begin{align}
\label{eq:44}
 \mathbb{E}\left(\|\hat{\Sigma}_{\underline{\hat{\mathcal{Y}}}_r}-\hat{\Sigma}_{\underline{\mathcal{Y}}_r}\|^{2}\right) \leq & \displaystyle\sum_{m_1,m_2=1}^{N}\sum_{r_1,r_2=1}^{r}\mathbb{E}\left(n^{-2}\displaystyle\left(\sum_{t=1}^n\mathcal{Z}_t\right)^2\right) \notag\\ \leq &  \sum_{m_1,m_2=1}^{N}\sum_{r_1,r_2=1}^{r}\left(\left(\mathbb{E}(\mathcal{Z}_1)\right)^2 + n^{-2}\text{Var}\left(\displaystyle\sum_{t=1}^n\mathcal{Z}_t\right)\right).
\end{align}
Let us consider \(1 \leq m \leq N\) and \(t \in \mathbb{Z}\). In view of Lemma \ref{lem:111}, using the Minkowski inequality and the stationarity of the process $(X_t)_{t\in\mathbb{Z}}$, we have
\begin{align}
\label{eq:45}
    \|\hat{\mathcal{Y}}_t(m)\|_2 = \left\|X_{t}X_{t+m}-\hat{c}_{m,0}\right\|_2\leq \left\|X_{t}X_{t+m}\right\|_2+\left\|(n-m)^{-1}\sum_{t=1}^{n-m}X_tX_{t+
    m}\right\|_2\leq 2\left\|X_0\right\|_4^2 < \infty.
\end{align}                                         
 Furthermore, due to stationarity, we also have
\begin{align}
\label{eq:46}
  \left\|\mathcal{Y}_t(m)\right\|_2^2 = \text{Var}(X_tX_{t+m})  = \text{Cov}(X_tX_{t+m},X_tX_{t+m})  \leq \|X_t\|_4^4+\|X_t\|_{2}^4 = \|X_0\|_4^4+\|X_0\|_{2}^4<\infty .
\end{align}
 Note also that the sequence \( (\hat{c}_{m,0} - c_{m,0})_{n\in\mathbb{N}} \) belongs to \( \mathcal{L}^{2} \). Moreover, by noting that \(\lim_{n\rightarrow\infty} \Var(\sqrt{n}F^{N,n})\) is an element of $\mathbb{R}^{N\times N}$ (see Lemma \ref{lem:7}), it follows that the $m$-th diagonal element \(\mathbb{E}(n|\hat{c}_{m,0} - c_{m,0}|^2)\) converges as \(n\) tends to infinity. In other words, there exists a positive constant $C(m)$ depending on $m$ such that 
\begin{align}
\label{eq:502}
 \left\|\hat{c}_{m,0} - c_{m,0}\right\|_{2} \underset{n \to \infty}{\sim} \frac{C(m)}{\sqrt{n}}.
\end{align}
 In view of Equations (\ref{eq:45}), (\ref{eq:46}) and (\ref{eq:502}), using the Cauchy Schwarz's inequality we obtain that
\begin{align}
\label{eq:47}
    \left|\mathbb{E}(\mathcal{Z}_t)\right| = & \ \left|\mathbb{E}\left(\hat{\mathcal{Y}}_{t-r_1}(m_1)\hat{\mathcal{Y}}_{t-r_2}(m_2)-\mathcal{Y}_{t-r_1}(m_1)\mathcal{Y}_{t-r_2}(m_2)\right)\right| \notag \\ 
     = & \ \left|\mathbb{E}\left\{\left(\hat{\mathcal{Y}}_{t-r_1}(m_1)-\mathcal{Y}_{t-r_1}(m_1)\right)\hat{\mathcal{Y}}_{t-r_2}(m_2)+\left(\hat{\mathcal{Y}}_{t-r_2}(m_2)-\mathcal{Y}_{t-r_2}(m_2)\right)\mathcal{Y}_{t-r_1}(m_1)\right\}\right| \notag\\
     \leq & \ \left\|\hat{\mathcal{Y}}_{t-r_1}(m_1)-\mathcal{Y}_{t-r_1}(m_1)\right\|_{2}\left\|\hat{\mathcal{Y}}_{t-r_2}(m_2)\right\|_2 + \left\|\hat{\mathcal{Y}}_{t-r_2}(m_2)-\mathcal{Y}_{t-r_2}(m_2)\right\|_{2}\left\|\mathcal{Y}_{t-r_1}(m_1)\right\|_2\notag\\
     = & \ \left\|\hat{c}_{m_1,0}-c_{m_1,0}\right\|_{2}\left\|\hat{\mathcal{Y}}_{t-r_2}(m_2)\right\|_2 + \left\|\hat{c}_{m_2,0}-c_{m_2,0}\right\|_{2}\left\|\mathcal{Y}_{t-r_1}(m_1)\right\|_2 \notag\\ 
     \leq & \ \frac{C(m_1,m_2)}{\sqrt{n}},
\end{align}
where $C(m_1, m_2)$ is a positive constant depending only on $m_1$ and $m_2$.

Additionally, owing to the stationarity property of the sequence $(X_t)_{t\in \mathbb{Z}}$, we obtain that 
\begin{align}
\label{eq:48}
  \text{Var}\left(\sum_{t=1}^{n}\mathcal{Z}_t\right) = \sum_{k=-n+1}^{n-1}(n-|k|)\text{Cov}(\mathcal{Z}_t,\mathcal{Z}_{t-k}).
\end{align}
Using Equations (\ref{eq:44}), (\ref{eq:47}) and (\ref{eq:48}), Lemmas~\ref{lem:11}, ~\ref{lem:12} and ~\ref{lem:13} we ultimately arrive at
\begin{align*}
 \mathbb{E}\left(r\|\hat{\Sigma}_{\underline{\hat{\mathcal{Y}}}_r}-\hat{\Sigma}_{\underline{\mathcal{Y}}_r}\|^{2}\right) \leq & \displaystyle\sum_{m_1,m_2=1}^{N}\sum_{r_1,r_2=1}^{r}r\left((\mathbb{E}\left(\mathcal{Z}_1)\right)^2 + n^{-2}\displaystyle\sum_{k=-n+1}^{n-1}(n-|k|)\text{Cov}(\mathcal{Z}_t,\mathcal{Z}_{t-k}) \right) \\ \leq
 & \displaystyle\sum_{m_1,m_2=1}^{N} \sum_{r_1,r_2=1}^{r} r\left((\mathbb{E}\left(\mathcal{Z}_1)\right)^2+ \displaystyle\sum_{k=-n+1}^{n-1}\frac{n-|k|}{n^2}\text{Cov}(\mathcal{Z}_t,\mathcal{Z}_{t-k})\right) \\ \leq & \displaystyle\sum_{m_1,m_2=1}^{N}\sum_{r_1,r_2=1}^{r} r \left((\mathbb{E}\left(\mathcal{Z}_1)\right)^2+n^{-1}\displaystyle\sum_{k=-\infty}^{\infty}\text{Cov}(\mathcal{Z}_t,\mathcal{Z}_{t-k})\right)\\
  \leq & \ N^{2}r^3\left((\mathbb{E}\left(\mathcal{Z}_1)\right)^2+n^{-1}\mathcal{C}_{13}\right)\xrightarrow[n\rightarrow\infty]{r \ = \ \mathrm{o}(n^{1/3})}0,
\end{align*}
where $\mathcal{C}_{13}$ is a strictly positive constant and independent of $r_1, r_2, m_1, m_2, r$, and $n$.

When $r = \mathrm{o}(n^{1/3})$ we thus have
\begin{align}
\label{eq:49}
  \sqrt{r}\|\hat{\Sigma}_{\underline{\hat{\mathcal{Y}}}_r}-\Sigma_{\underline{\mathcal{Y}}_r}\| = \mathrm{o}_{\mathbb{P}}(1).
\end{align}
Finally, noting that 
\begin{align*}
   \hat{\Sigma}_{\underline{\hat{\mathcal{Y}}_r}}-\Sigma_{\mathcal{Y}_r} = \left(\hat{\Sigma}_{\underline{\hat{\mathcal{Y}}}_r}-\hat{\Sigma}_{\underline{\mathcal{Y}}_r}\right) + \left(\hat{\Sigma}_{\underline{\mathcal{Y}}_r}-\Sigma_{\mathcal{Y}_r}\right),
\end{align*}
it follows from Equations (\ref{eq:501}) and (\ref{eq:49}) that when $r = \mathrm{o}(n^{1/3})$
\begin{align*}
\sqrt{r}\|\hat{\Sigma}_{\underline{\hat{\mathcal{Y}}}_r}-\Sigma_{\underline{\mathcal{Y}}_r}\| = \mathrm{o}_{\mathbb{P}}(1).
\end{align*}
The remaining results are similarly obtained, there by concluding the proof.
\qed

\begin{l1}\label{lem:16}
Under the assumptions of Theorem~\ref{thm:4}, we have
\begin{align*}
    \sqrt{r} \left\|\hat{\Sigma}_{\underline{\mathcal{Y}}_r}^{-1}-\Sigma_{\underline{\mathcal{Y}}_r}^{-1}\right\| = \mathrm{o}_{\mathbb{P}}(1)
\end{align*}
as $n\rightarrow\infty$ when $r = \mathrm{o}(n^{1/3})$ and $r\rightarrow\infty$.
\end{l1}

\subsubsection*{Proof.}
The proof is analogous to that given in \citet[Lemma A.6, p. 27]{francq2003goodness} (see also \citet[Lemma 6 of the supplementary material]{bmcf}) and it is omitted.
\qed

\begin{l1}\label{lem:17}
Under the assumptions of Theorem~\ref{thm:4}, we have
\begin{align*}
    \sqrt{r}\|\underline{\varphi}_r^{\star} - \underline{\varphi}_r\| \xrightarrow[n\rightarrow \infty]{r \ = \ \mathrm{o}(n^{1/3})} 0.
\end{align*}
\end{l1}

\subsubsection*{Proof.}
Recall that, in view of Equations (\ref{eq:20}) and (\ref{eq:22}), for $t\in\mathbb{Z}$ we have 
\begin{align}
\label{eq:50}
    \underline{\bm{\varphi}}_r\underline{\mathcal{Y}}_{r,t}+u_{r,t} = \mathcal{Y}_t = \underline{\bm\varphi}_r^{\star}\underline{\mathcal{Y}}_{r,t}+u_{r,t}^{\star},
\end{align}
where $u_{r,t}^{\star} := \displaystyle\sum_{i=r+1}^{\infty}\varphi_{i}\mathcal{Y}_{t-i}+u_t$.

Considering the orthogonality between $u_{r,t}$ and $\mathcal{Y}_t$, it follows from Equation (\ref{eq:50}) that
\begin{align}
\label{eq:51}
 \underline{\varphi}_r^{\star} - \underline{\varphi}_r = -\mathbb{E}(u_{r,t}^{\star}\underline{\mathcal{Y}}_{r,t}^{'})\Sigma_{\underline{\mathcal{Y}}_r}^{-1} := \Sigma_{\underline{\mathcal{Y}}_r}^{\star}\Sigma_{\underline{\mathcal{Y}}_r}^{-1}
\end{align}
where $\Sigma_{\underline{\mathcal{Y}}_r}^{\star} = \mathbb{E}(u_{r,t}^{\star}\underline{\mathcal{Y}}_{r,t}^{'})$. 

Furthermore, thanks to Minkowski's inequality, we obtain on the one hand
\begin{align}
\label{eq:52}
   \left\|\Sigma_{\underline{\mathcal{Y}}_r}^{\star}\right\| = & \left\|\mathbb{E}\left\{\left(\displaystyle\sum_{i=r+1}^{\infty}\varphi_{i}\mathcal{Y}_{t-i}+u_t\right)\underline{\mathcal{Y}}_{r,t}^{'}\right\}\right\| \notag\\
   = & \left\|\displaystyle\sum_{i=r+1}^{\infty}\varphi_{i}\mathbb{E}\left(\mathcal{Y}_{t-i}\underline{\mathcal{Y}}_{r,t}^{'}\right)\right\| \notag\\
   \leq & \displaystyle\sum_{i=r+1}^{\infty}\left\|\varphi_i\right\|\left\|\mathbb{E}(\mathcal{Y}_{t-i}\underline{\mathcal{Y}}_{r,t}^{'})\right\| \notag\\
   \leq & \displaystyle\sum_{i=1}^{\infty}\left\|\varphi_{r+i}\right\|\left\|\mathbb{E}(\mathcal{Y}_{t-i-r}\underline{\mathcal{Y}}_{r,t}^{'})\right\|.
\end{align}
 Using Equation (\ref{eq:505}) we obtain on the other hand
\begin{align}
\label{eq:53}
    \left\|\mathbb{E}(\mathcal{Y}_{t-i-r}\underline{\mathcal{Y}}_{r,t}^{'})\right\|^{2}\leq \displaystyle\sum_{k=1}^r\mathbb{E}\{(\mathcal{Y}_{t-i-r},\mathcal{Y}_{t-k}^{'})\}^{2}\leq r\displaystyle\sum_{1\leq m_1,m_2\leq N} \displaystyle\sum_{k=1}^{r}\left\|\mathcal{Y}_{t-i-r}(m_1)\right\|_{2}^{2}\left\|\mathcal{Y}_{t-k}(m_2)\right\|_{2}^{2}\leq \mathcal{C}_{14}rN^2.
\end{align}
In view of Equations (\ref{eq:52}), (\ref{eq:53}) and using the fact that $\|\varphi_{i}\| = \mathrm{o}(i^{-2})$ we arrive at
\begin{align*}
  \sqrt{r}\|\underline{\varphi}_r^{\star} - \underline{\varphi}_r\| \leq \displaystyle\sum_{i=1}^{\infty}\left(\mathrm{o}(i+r)^{-2}\right)\mathcal{C}_{14}^{1/2}rN\left(\sup_{s\geq 1}\left\|\Sigma_{\underline{\mathcal{Y}}_s}^{-1}\right\|\right) \xrightarrow[r\rightarrow \infty]{} 0
\end{align*}
where $\mathcal{C}_{14}$ is a positive constant that is independent of $n$ and $r$.

\noindent Hence, the proof is concluded. 
\qed


\begin{l1}\label{lem:18}
Under the assumptions of Theorem~\ref{thm:4}, we have
\begin{align*}
 \sqrt{r}\left\|\hat{\underline{\bm{\varphi}}_r}-\underline{\bm{\varphi}}_r\right\| = \mathrm{o}_{\mathbb{P}}(1)
\end{align*}
as $n\rightarrow\infty$ when $r=\mathrm{o}(n^{1/3})$.
\end{l1}
\subsubsection*{Proof.}
Taking into account the orthogonality condition between $u_{r,t}$ and $\underline{\mathcal{Y}}_{r,t}$ in the regression of $\mathcal{Y}_t$ on the family $(\mathcal{Y}_{t-i})_{1 \leq i \leq r}$ (see Equation (\ref{eq:22})), it follows that
\begin{align*}
  \mathbb{E}(\mathcal{Y}_t\mathcal{Y}_{r,t}^{'}) = \underline{\bm{\varphi}}_r\mathbb{E}(\underline{\mathcal{Y}}_{r,t}\mathcal{Y}_{r,t}^{'}).
\end{align*}
In other words, we consider the previously introduced notations
\begin{align*}
 \underline{\bm{\varphi}}_r = \Sigma_{\mathcal{Y},\underline{\mathcal{Y}}_{r}}\Sigma_{\underline{\mathcal{Y}}_r}^{-1}  \ \ \text{and} \ \ \hat{\underline{\bm{\varphi}}}_r = \hat{\Sigma}_{\hat{\mathcal{Y}},\underline{\hat{\mathcal{Y}}}_{r}}\hat{\Sigma}_{\underline{\hat{\mathcal{Y}}}_r}^{-1}.
\end{align*}
Consequently, by employing on the one hand the triangle inequality, and on the other hand Lemmas~\ref{lem:9}, ~\ref{lem:15} and ~\ref{lem:16}, when $n\rightarrow\infty$ and $r = \mathrm{o}(n^{1/3})$ we have 
\begin{align*}
 \sqrt{r}\left\|\hat{\underline{\bm{\varphi}}_r}-\underline{\bm{\varphi}}_r\right\| = &  \sqrt{r}\left\|\left(\hat{\Sigma}_{\hat{\mathcal{Y}},\underline{\hat{\mathcal{Y}}}_{r}}-\Sigma_{\mathcal{Y},\underline{\mathcal{Y}}_{r}}+\Sigma_{\mathcal{Y},\underline{\mathcal{Y}}_{r}}\right)\hat{\Sigma}_{\underline{\hat{\mathcal{Y}}}_r}^{-1}-\Sigma_{\mathcal{Y},\underline{\mathcal{Y}}_{r}}\left(\Sigma_{\underline{\mathcal{Y}}_r}^{-1}-\hat{\Sigma}_{\underline{\hat{\mathcal{Y}}}_r}^{-1}+\hat{\Sigma}_{\underline{\hat{\mathcal{Y}}}_r}^{-1}\right)\right\| \\
  = & \sqrt{r}\left\|\left(\hat{\Sigma}_{\hat{\mathcal{Y}},\underline{\hat{\mathcal{Y}}}_{r}}-\Sigma_{\mathcal{Y},\underline{\mathcal{Y}}_{r}}\right)\hat{\Sigma}_{\underline{\hat{\mathcal{Y}}}_r}^{-1}+\Sigma_{\mathcal{Y},\underline{\mathcal{Y}}_{r}}\left(\hat{\Sigma}_{\underline{\hat{\mathcal{Y}}}_r}^{-1}-\Sigma_{\underline{\mathcal{Y}}_r}^{-1}\right)\right\| \\ \leq & \sqrt{r}\left\|\hat{\Sigma}_{\hat{\mathcal{Y}},\underline{\hat{\mathcal{Y}}}_{r}}-\Sigma_{\mathcal{Y},\underline{\mathcal{Y}}_{r}}\right\|\left(\sup_{r\geq 1}\left\|\hat{\Sigma}_{\underline{\hat{\mathcal{Y}}}_r}^{-1}\right\|\right) + \left(\sup_{r\geq 1}\left\|\Sigma_{\mathcal{Y},\underline{\mathcal{Y}}_{r}}\right\|\right)\sqrt{r}\left\|\hat{\Sigma}_{\underline{\hat{\mathcal{Y}}}_r}^{-1}-\Sigma_{\underline{\mathcal{Y}}_r}^{-1}\right\| = \mathrm{o}_{\mathbb{P}}(1).
\end{align*}
This completes the proof.
\qed

\subsubsection*{Proof of Theorem~\ref{thm:4}.}
Given the expressions for $I$ and $\hat{I}^{SP}$, proving Theorem~\ref{thm:4} essentially involves demonstrating that $\hat{\bm{\varphi}}(1)$ converges in probability to $\bm{\varphi}(1)$ and $\hat{\Sigma}_{\hat{u}_r}$ converges in probability to $\Sigma_{u}$. By denoting $I_{N \times N}$ as the identity matrix of order $N$, it can be observed that
\begin{align*}
    \hat{\varphi}_{r,i} - \varphi_{r,i} = \left(\underline{\hat{\bm{\varphi}}}_r-\underline{\bm{\varphi}}_r\right)\mathcal{K}_i \quad \text{and} \quad \varphi_i-\varphi_{r,i} = \left(\underline{\bm{\varphi}}_r^{\star}-\underline{\bm{\varphi}}_r\right)\mathcal{K}_i \qquad \forall \ 1\leq i\leq r
\end{align*}
with
\begin{align*}
\mathcal{K}_i = \left( \begin{array}{c}
O_{N\times N} \\
\vdots \\
I_{N\times N} \\
\vdots \\
O_{ N\times N}
\end{array} \right)\in\mathbb{R}^{rN\times N},
\end{align*}
where $I_{N\times N}$ corresponds to the \(i\)-th block of $\mathcal{K}_{i}$.

Moreover, Equations (\ref{eq:20}) and (\ref{eq:21}) allow us to recall
\begin{align*}
   \hat{\bm{\varphi}}(L) = \sum_{i=1}^{r}\hat{\varphi}_{r,i}L^{i} \qquad \text{and} \qquad \bm{\varphi}(L) = \sum_{i=1}^{\infty}\varphi_iL^{i}.
\end{align*}
Using Lemmas~\ref{lem:17}, ~\ref{lem:18} and Equation (\ref{eq:505}), we obtain
\begin{align*}
    \left\|\hat{\bm{\varphi}}(1)-\bm{\varphi}(1)\right\| \leq & \ \left\|\displaystyle\sum_{i=1}^{r}\left(\hat{\varphi}_{i,r}-\varphi_{r,i}\right)\right\|+\left\|\displaystyle\sum_{i=1}^{r}\left(\varphi_{r,i}-\varphi_i\right)\right\| + \left\|\sum_{i=r+1}^{\infty}\varphi_i\right\| \\ \leq & \ 
    \left\|\left(\underline{\hat{\bm{\varphi}}}_r-\underline{\bm{\varphi}}_r\right)\sum_{i=1}^{r}\mathcal{K}_i\right\|+\left\|\left(\underline{\bm{\varphi}}_r^{\star}-\underline{\bm{\varphi}}_r\right)\sum_{i=1}^{r}\mathcal{K}_i\right\|+ \displaystyle\sum_{i=r+1}^{\infty}\left\|\varphi_i\right\|
    \\ \leq & \ \left\|\sum_{i=1}^{r}\mathcal{K}_i\right\|\left\{\left\|\underline{\hat{\bm{\varphi}}}_r-\underline{\bm{\varphi}}_r\right\| + \left\|\left(\underline{\bm{\varphi}}_r^{\star}-\underline{\bm{\varphi}}_r\right)\right\|\right\}+ \sum_{i=r+1}^{\infty}\left\|\varphi_i\right\| \\ \leq & \ \sqrt{N} \sqrt{r} \left\{\left\|\underline{\hat{\bm{\varphi}}}_r-\underline{\bm{\varphi}}_r\right\| + \left\|\left(\underline{\bm{\varphi}}_r^{\star}-\underline{\bm{\varphi}}_r\right)\right\|\right\}+\displaystyle\sum_{i=r+1}^{\infty} \left\|\varphi_i\right\| =\mathrm{o}_{\mathbb{P}}(1).
\end{align*}
We also have
\begin{align*}
   \hat{\Sigma}_{\hat{u}_r} = \hat{\Sigma}_{\hat{\mathcal{Y}}}-\underline{\hat{\bm{\varphi}}}_{r}\hat{\Sigma}_{\hat{\mathcal{Y}},\underline{\hat{\mathcal{Y}}}_r}^{'},
\end{align*}
  and
\begin{align*}
    \Sigma_{u} = & \ \mathbb{E}(u_tu_t^{'}) = \mathbb{E}\left(u_t\mathcal{Y}_t^{'}\right)=\mathbb{E}\left\{\left(\mathcal{Y}_t-\displaystyle\sum_{i=1}^{\infty}\varphi_i\mathcal{Y}_{t-i}\right)\mathcal{Y}_t^{'}\right\} \\
    = & \ \Sigma_{\mathcal{Y}}-\displaystyle\sum_{i=1}^{\infty}\varphi_i\mathbb{E}\left(\mathcal{Y}_{t-i}\mathcal{Y}_t^{'}\right) = \Sigma_{\mathcal{Y}}-\underline{\bm{\varphi}}_r^{\star}\Sigma_{\mathcal{Y},\underline{\mathcal{Y}}_r}^{'}-\displaystyle\sum_{i=r+1}^{\infty}\varphi_i\mathbb{E}\left(\mathcal{Y}_{t-i}\mathcal{Y}_{t}^{'}\right).
\end{align*}
Therefore, it follows that
\begin{align*}
    \left\|\hat{\Sigma}_{\hat{u}_r}-\Sigma_{u}\right\| = &\left\|\hat{\Sigma}_{\hat{\mathcal{Y}}}-\Sigma_{\mathcal{Y}}-\left(\underline{\hat{\bm{\varphi}}}_r-\underline{\bm{\varphi}}_r^{\star}\right)\hat{\Sigma}_{\hat{\mathcal{Y}},\underline{\hat{\mathcal{Y}}}_{r}}^{'}-\underline{\bm{\varphi}}_r^{\star}\left(\hat{\Sigma}_{\hat{\mathcal{Y}},\underline{\hat{\mathcal{Y}}}_{r}}^{'}-\Sigma_{\mathcal{Y},\underline{\mathcal{Y}}_r}^{'}\right)+\sum_{i=r+1}^{\infty}\varphi_i\mathbb{E}\left(\mathcal{Y}_{t-i}\mathcal{Y}_t^{'}\right)\right\| \\ 
    \leq & \left\|\hat{\Sigma}_{\hat{\mathcal{Y}}}-\Sigma_{\mathcal{Y}}\right\|+ \left\|\left(\underline{\hat{\bm{\varphi}}}_r-\underline{\bm{\varphi}}_r^{\star}\right)\left(\hat{\Sigma}_{\hat{\mathcal{Y}},\underline{\hat{\mathcal{Y}}}_{r}}^{'}-\Sigma_{\mathcal{Y},\underline{\mathcal{Y}}_r}^{'}\right)\right\| + \left\|\left(\underline{\hat{\bm{\varphi}}}_r-\underline{\bm{\varphi}}_r^{\star}\right)\Sigma_{\mathcal{Y},\underline{\mathcal{Y}}_r}^{'}\right\|\\ & + \left\|\underline{\bm{\varphi}}_r^{\star}\left(\hat{\Sigma}_{\hat{\mathcal{Y}},\underline{\hat{\mathcal{Y}}}_r}^{'}-\Sigma_{\mathcal{Y},\underline{\mathcal{Y}}_r}^{'}\right)\right\| + \left\|\sum_{i=r+1}^{\infty}\varphi_{i}\mathbb{E}\left(\mathcal{Y}_{t-i}\mathcal{Y}_t^{'}\right)\right\| \\ \leq & \left\|\hat{\Sigma}_{\hat{\mathcal{Y}}}-\Sigma_{\mathcal{Y}}\right\|+\left\|\underline{\hat{\bm{\varphi}}}_r-\underline{\bm{\varphi}}_r^{\star}\right\|\left\|\hat{\Sigma}_{\hat{\mathcal{Y}},\underline{\hat{\mathcal{Y}}}_{r}}^{'}-\Sigma_{\mathcal{Y},\underline{\mathcal{Y}}_r}^{'}\right\|+\left\|\underline{\hat{\bm{\varphi}}}_r-\underline{\bm{\varphi}}_r^{\star}\right\| \left(\sup_{s\geq 1}\left\|\Sigma_{\mathcal{Y},\underline{\mathcal{Y}}_s}^{'}\right\|\right)\\ + & \left\|\underline{\bm{\varphi}}_r^{\star}\right\|\left\|\hat{\Sigma}_{\hat{\mathcal{Y}},\underline{\hat{\mathcal{Y}}}_r}^{'}-\Sigma_{\mathcal{Y},\underline{\mathcal{Y}}_r}^{'}\right\| + \displaystyle\sum_{i=r+1}^{\infty}\left\|\varphi_i\right\|\left\|\mathbb{E}\left(\mathcal{Y}_{t-i}\mathcal{Y}_{t}^{'}\right)\right\|.
\end{align*}
By Lemma~\ref{lem:14}, we established that $\hat{\Sigma}_{\hat{\mathcal{Y}}}-\Sigma_{\mathcal{Y}} = \mathrm{o}_{\mathbb{P}}(r^{-1/2})$. Following Lemmas~\ref{lem:17} and ~\ref{lem:18}, it is observed that $\left\|\underline{\hat{\bm{\varphi}}}_r-\underline{\varphi}_r^{\star}\right\| = \mathrm{o}_{\mathbb{P}}(1)$. Lemma~\ref{lem:16} further indicates that $\left\|\hat{\Sigma}_{\hat{\mathcal{Y}},\underline{\hat{\mathcal{Y}}}_r}^{'}-\Sigma_{\mathcal{Y},\underline{\mathcal{Y}}_r}^{'}\right\| = \mathrm{o}_{\mathbb{P}}(r^{-1/2})$. Additionally, according to Lemma~\ref{lem:9}, we have $\sup_{r\geq 1}\left\|\Sigma_{\mathcal{Y},\underline{\mathcal{Y}}_r}^{'}\right\| = \mathrm{O}(1)$. Finally, using (\ref{eq:505}), we also derive that $\left\|\underline{\bm{\varphi}}_r^{\star}\right\|\leq \left(\sum_{i=1}^{\infty}\text{Tr}(\varphi_i\varphi_i^{'})\right)^{1/2}<\infty$ and $\left\|\mathbb{E}\left(\mathcal{Y}_{t-i}\mathcal{Y}_{t}^{'}\right)\right\| = \mathrm{O}(1)$. Together, these results permit us to complete the proof.
\qed

\bibliographystyle{apalike}
\bibliography{bibliography}


\end{document}